\documentclass[12pt,twoside]{amsart}
\usepackage{amssymb,amsmath,amsfonts,latexsym,mathtools,tikz,tikz-cd}
\usepackage{amscd}
\usepackage{xcolor}
\usepackage{array,multirow,tabularx,longtable}
\usepackage[normalem]{ulem}
\usepackage{comment}
\usepackage{hyperref}
\usepackage[abbrev,alphabetic]{amsrefs}
\usepackage[margin=1.25in]{geometry}

\usepackage{ascmac} 

\title[Discriminant divisors for conic bundles]
{Discriminant divisors for conic bundles} 
\author{Hiromu Tanaka} 

\address{Graduate School of Mathematical Sciences, 
The University of Tokyo, 
3-8-1 Komaba, Meguro-ku, Tokyo 153-8914, JAPAN} 
\email{tanaka@ms.u-tokyo.ac.jp}

\newcommand{\bdl}[0]{{\operatorname{bdl}}}

\newcommand{\Diff}[0]{{\operatorname{Diff}}}

\newcommand{\GL}[0]{{\operatorname{GL}}}
\newcommand{\Hilb}[0]{{\operatorname{Hilb}}}
\newcommand{\Univ}[0]{{\operatorname{Univ}}}
\newcommand{\univ}[0]{{\operatorname{univ}}}

\newcommand{\red}[0]{{\operatorname{red}}}

\newcommand{\Proj}[0]{{\operatorname{Proj}}}
\newcommand{\Spec}[0]{{\operatorname{Spec}}}

\newcommand{\Pic}[0]{{\operatorname{Pic}}}

\newcommand{\Ex}[0]{{\operatorname{Ex}}}

\newcommand{\Frac}[0]{{\operatorname{Frac}}}
\newtheorem{thm}{Theorem}[section]
\newtheorem{lem}[thm]{Lemma}
\newtheorem{cor}[thm]{Corollary}
\newtheorem{prop}[thm]{Proposition}

\newtheorem{step}{Step}

\theoremstyle{definition}
\newtheorem{ex}[thm]{Example}

\newtheorem{dfn}[thm]{Definition}

\newtheorem{rem}[thm]{Remark}
      
\newtheorem{nota}[thm]{Notation}         

\newtheorem{nasi}[thm]{}

\makeatletter
  
  \@addtoreset{equation}{thm}
  \makeatother

\newcommand{\cred}{\color{black}}

\newcommand{\MO}{\mathcal{O}}
\newcommand{\A}{\mathbb{A}}

\newcommand{\Q}{\mathbb{Q}}
\newcommand{\Z}{\mathbb{Z}}
\newcommand{\F}{\mathbb{F}}

\newcommand{\p}{\mathfrak{p}}
\newcommand{\m}{\mathfrak{m}}
\newcommand{\n}{\mathfrak{n}}
\renewcommand{\P}{\mathbb{P}}
\begin{document}

\begin{abstract}
We study some foundational properties on discriminant divisors 
for generically smooth conic bundles. In particular, we extend the formula $\Delta_f \equiv -f_*K_{X/T}^2$ to arbitrary characteristics. 
\end{abstract}

\maketitle



\tableofcontents


\section{Introduction}


Conic bundles played a crucial role in the classification of Fano threefolds \cite{MM81}, \cite{MM83}. 
Especially, Mori--Mukai establish the following formula 
\begin{equation}\label{e-MM-0}
\Delta_f \equiv -f_*K_{X/T}^2, 
\end{equation}
where {\cred $f: X \to T$ is a conic bundle from a smooth projective threefold to a smooth projective surface over an algebraically closed field of characteristic zero.} 
{\cred Originally,} 
$\Delta_f$ was defined as the reduced closed subscheme of $T$ 
that parametrises the singular fibres. 
Although this naive definition does not behave nicely in characteristic two, \cite{ABBB21} introduced a well-behaved definition. 
The main objective of this article is to establish the above formula (\ref{e-MM-0}) in arbitrary characteristic by using the definition of \cite{ABBB21}.

\begin{thm}[Remark \ref{r-disc-bdl-div}, Theorem \ref{t-cbf}]\label{intro-MM}
Let $k$ be an algebraically closed field. 
Let $f: X \to T$ be a generically smooth conic bundle, where $X$ and $T$ are smooth projective varieties over $k$. 
Then the numerical equivalence 
\[
\Delta_f \equiv -f_*(K_{X/T}^2)
\]
holds, where $K_{X/T} := K_X -f^*K_T$. 
\end{thm}

This theorem will be applied in forthcoming articles on the classification of Fano threefolds in positive characteristic \cite{Tan-Fano1}, \cite{Tan-Fano2}, \cite{AT}, \cite{Tan-Fano4}. 

\subsection{Overview of proofs and contents}\label{ss-intro-overview}


In Section \ref{s-disc-Sigma}, 
we shall introduce a discriminant scheme $\Delta_f$ for a conic bundle $f: X \to T$, 
which is a closed subscheme of $T$. 
The definition is designed to satisfy the following two axiomatic properties (I) and (II). 
\begin{enumerate}
    \item[(I)] Given a point $t \in T$, 
    $X_t$ is not smooth if and only if  $t \in \Delta_f$. 
    \item[(II)] An equality $\Delta_{f'} = \beta^{-1}\Delta_f$ of closed subschemes 
    holds for every carterian diagram of noetherian schemes
    \[
    \begin{CD}
    X' @>\alpha >> X\\
    @VVf' V @VVf V\\
    T' @>\beta >> T, 
    \end{CD}
    \]
    where $f$ and $f'$ are conic bundles.  
\end{enumerate}
The actual definition of $\Delta_f$ is given as follows. 
We first consider the case when the base scheme $T$ is affine and $X$ is a conic on $\P^2_T$. 
We have $T=\Spec\,A$, there is a closed embedding $X \subset \P^2_A$, and 
\[
X = \Proj\,A[x, y, z]/(Q), 
\quad Q := ax^2 + by^2+cz^2 + \alpha yz + \beta zx +\gamma xy, 
\quad a, b, c, \alpha, \beta, \gamma \in A.
\]
In this case, $\Delta_f$ is explicitly defined as follows: 
\begin{equation}\label{e-intro-delta}
\Delta_f = \Spec\,(A/\delta(Q)A), \quad \text{where}\quad 
\delta(Q) := 4abc + \alpha \beta \gamma - a \alpha^2 - b\beta^2 -c \gamma^2 \in A. 
\end{equation}
We shall check that 
the definition $\Delta_f$ is independent of the  choice of the closed embedding $X \subset \P^2_A$ and so on. 
The proof is done by using the fact that two embeddings are related by a linear transform. 
As mentioned already, this definition is based on \cite{ABBB21}. 
Although \cite{ABBB21} does not check the well-definedness, 
the equation (\ref{e-intro-delta}) is extracted from \cite{ABBB21}, 
which would be the most essential part to reach the definition explained above. 
Similarly, we also introduce a closed subscheme $\Sigma_f$ of $T$ such that, 
in addition to the property corresponding to (II), 
$\Sigma_f$ consists of the points $t$ of $T$ such that $X_t$ is not geometrically reduced. 




In Section \ref{s-sing-ambient}, we study the relation between 
$\Delta_f$ and the singularities of $X$. 
We here only treat the case when the base scheme $T$ is regular. 
We shall prove that the following  are equivalent (Theorem \ref{t-sm-red-fib}):
\begin{itemize}
    \item $\Delta_f$ is regular. 
    \item $X$ is regular and any fibre of $f$ is geometrically reduced. 
\end{itemize}
In particular, if both $X$ and $T$ are regular, 
then the singular locus of $\Delta_f$ is set-theoretically equal to $\Sigma_f$. 


In {\cred Section \ref{s-MM-formula}}, we prove the formula (Theorem \ref{intro-MM}): 
\begin{equation}\label{e-intro-MM2}
\Delta_f \equiv -f_*K_{X/T}^2. 
\end{equation}
Note that the argument in characteristic zero \cite[Proposition 6.8]{MM83} does not work (cf. Remark \ref{r-Bertini-fail}). 
We here overview some of the ideas of the proof. 
As a toy case, let us consider the case when $f : X \to T$ coincides with 
$\Univ_{\P^2_k/k}^{\theta} \to \Hilb_{\P^2_k/k}^{\theta}$, 
where 
\begin{itemize}
    \item $k$ is an algebraically closed field, 
    \item $\theta$ is the Hilbert polynomial of conics, 
    \item $\Hilb_{\P^2_k/k}^{\theta}$ denotes the Hilbert scheme parametrising 
    all the conics on $\P^2$, and 
    \item $\Univ_{\P^2_k/k}^{\theta}$ is its universal family. 
\end{itemize}
In particular, we have $\Hilb_{\P^2_k/k}^{\theta} \simeq \P^5_k$. 
In order to check the numerical equivalence (\ref{e-intro-MM2}), 
it is enough to find a line $L$ on $\Hilb_{\P^2_k/k}^{\theta} \simeq \P^5_k$ 
such that 
$\Delta_f \cdot L = -(f_*K_{X/T}^2) \cdot L$. 
This is carried out by taking a general line $L$ such that $\Delta_f \cap L$ is smooth. 
For the general case, the problem is reduced, by standard argument, 
to the case when $T$ is a smooth projective curve. 
We embed $T$ to the relative Hilbert scheme  $\Hilb_{\P_T(E)/T}^{\theta}$ for $E := f_*\omega_{X/T}^{-1}$, which is a locally free sheaf on $T$ of rank $3$. 
Then the problem is further reduced to the case when 
$f$ coincides with a relative version $\Univ_{\P_T(E)/T}^{\theta} \to \Hilb_{\P_T(E)/T}^{\theta}$ of 
$\Univ_{\P^2_k/k}^{\theta} \to \Hilb_{\P^2_k/k}^{\theta}$. 
By $\rho(\Hilb_{\P^2_k/k}^{\theta}) = \rho(\P^5)=1$, 
it was enough to find a line $L$ satisfying $\Delta_f \cdot L = -(f_*K_{X/T}^2) \cdot L$ in the above case. 
For our case, we can check that $\rho(\Hilb_{\P_T(E)/T}^{\theta})=2$. 
Then the proof is carried out by finding suitable two curves $C_1$ and $C_2$ such that 
$\Delta_f \cdot C_1 = -(f_*K_{X/T}^2) \cdot C_1$ and 
$\Delta_f \cdot C_2 = -(f_*K_{X/T}^2) \cdot C_2$. 
For more details, see Theorem \ref{t-cbf}.


In order to justify the above argument, 
we shall introduce an invertible sheaf $\Delta_f^{\bdl}$, called the discriminant bundle, which is a variant of the discriminant scheme $\Delta_f$ (Subsection \ref{s-disc-bdl}).  
If $f : X \to T$ is a generically smooth conic bundle such that $T$ is a noetherian integral scheme, 
then $\Delta_f$ is an effective Cartier divisor on $T$. 
However, $\Delta_f$ is useless when $f : X \to T$ is a wild conic bundle, i.e., no fibre is smooth (more explicitly, we have $\Delta_f =T$). 
On the other hand, the discriminant bundle $\Delta^{\bdl}_f$ is an invertible sheaf on $T$ even if $f$ is a wild conic bundle. 
Furthermore, the discriminant bundles $\Delta^{\bdl}_f$ satisfy the following properties (Theorem \ref{t-disc-bdl-bc}). 
\begin{enumerate}
    \item[(A)] If $f: X \to T$ is a generically smooth conic bundle such that $T$ is a noetherian integral scheme, then $\Delta^{\bdl}_f \simeq \MO_T(\Delta_f)$. 
    \item[(B)] 
 An isomorphism $\Delta_{f'}^{\bdl} \simeq \beta^{*}\Delta_f^{\bdl}$ of invertible sheaves  
    holds for every carterian diagram 
    \[
    \begin{CD}
    X' @>\alpha >> X\\
    @VVf' V @VVf V\\
    T' @>\beta >> T, 
    \end{CD}
    \]
    where $T$ and $T'$ are smooth varieties over a field.  
\end{enumerate}
The definition of $\Delta^{\bdl}_f$ is given as follows. 
For a conic bundle $f: X \to T$, we have the following cartesian diagram 
for the Hilbert polynomial $\theta$ of conics: 
\[
\begin{CD}
    X @>\varphi >> \widetilde{X} := \Univ^{\theta}_{\P_T(f_*\omega_{X/T}^{-1})/T} \\
    @VVf V @VV\widetilde f V\\
    T @>\psi >> \widetilde{T} := \Hilb^{\theta}_{\P_T(f_*\omega_{X/T}^{-1})/T}.  
    \end{CD}
\]
Then $\Delta^{\bdl}_f$ is defined by 
$\Delta^{\bdl}_f := \psi^*\MO_{\widetilde{T}}(\Delta_{\widetilde f})$. 
It is clear that (A) holds. 
The property (B) can be checked as follows. 
By standard argument, we may assume that 
$\beta : T' \to T$ is either smooth or 
a closed immersion such that  $\beta(T')$ is an effective Cartier divisor. 
For each case, it is not so hard to check (B). 

In Section \ref{s-surface}, we study $\Delta_f$ for surface conic bundles. 
We restrict ourselves to treating the case when 
$T$ is a smooth curve and $X$ is a surface having at worst canonical singularities. 
Since we are interested in the relation between $\Delta_f$ and 
the singularities of $X$, 
let us assume that there is a closed point $0 \in T$ such that 
$X_0$ is {\cred the} unique singular fibre. 
We shall first prove that 
\[
\deg \Delta_f = m-1, 
\]
where $m$ is the number of the irreducible components of the central fibre $Y_0$ for the minimal resolution $\varphi : Y \to X$ of $X$. 
Given such a conic bundle $f: X \to T$, 
we shall prove that 
the singularity of $X$ is determined by $\deg \Delta_f$ and whether the $X_0$ is reduced (Theorem \ref{t-disc-vs-sing}). 
We then exhibit several examples. 

In Section \ref{s-char2-ex}, we observe some phenomena 
which occur only in {\cred characteristic} two. 
For example, the following hold in characteristic $\neq 2$ 
(Proposition \ref{p-sm-nonred-fib}). 
\begin{enumerate}
\item If $f : X \to T$ is a conic bundle of smooth varieties, then $\Delta_f$ is a reduced divisor. 
\item If  $f : X \to T$ is a conic bundle of smooth varieties with $\dim T=2$, 
then $\Delta_f$ is normal crossing. 
\end{enumerate}
We shall see that both the properties fail in characteristic two (Example \ref{e-3fold-char2}). 
Furthermore, (2) does not hold even after replacing $\Delta_f$ by $(\Delta_f)_{\red}$. 
Roughly speaking, the singularities of $\Delta_f$ and $(\Delta_f)_{\red}$ can be arbitrarily bad even if $f : X \to T$ is a generically smooth conic bundle 
from a smooth threefold $X$ to a smooth surface $T$. 







\begin{rem}\label{r-Bea77-Sar82}
{\cred 
As explained above, we shall establish some foundational results on discriminant divisors $\Delta_f$ for conic bundles $f: X \to T$, 
where $X$ and $T$ are noetherian schemes. 
Most results in Section \ref{s-disc-Sigma} and 
Section \ref{s-sing-ambient} are essentially contained 
in \cite[Chapitre I]{Bea77} and \cite[Section 1]{Sar82} 
for the case when $X$ and $T$ are smooth varieties over an algebraically closed field of characteristic $\neq 2$. 
As far as the author knows, 
the contents in Section \ref{s-disc-Sigma} and 
Section \ref{s-sing-ambient} are new even when $X$ and $T$ are smooth varieties over an algebraically closed field of characteristic two. 
}
\end{rem}


\vspace{3mm}

{\cred 
\textbf{Acknowledgements:} 
The author 
thanks the referee for reading the manuscript carefully and for suggesting several  improvements.
The author was funded by JSPS KAKENHI Grant numbers JP22H01112 and JP23K03028.}

\section{Preliminaries}



\subsection{Notation}\label{ss:notation}



\begin{enumerate}
    \item We say that $X$ is a {\em variety} (over a field $k$) if $X$ is an integral scheme 
    that is separated and of finite type over $k$. We say that $X$ is a {\em curve} (resp. {\em surface}, resp. {\em threefold}) if $X$ is a variety of dimension one (resp. two, resp. three). 
\item We say that a noetherian scheme $X$ is {\em excellent} if 
all the stalks $\MO_{X, x}$ are excellent. 
\item Given a morphism $f:X \to S$ of schemes and a point $s \in S$, 
$X_s$ denotes the fibre of $f$ over $s$, i.e., $X_s := X \times_S \Spec\,\kappa(s)$, 
where $\kappa(s)$ denotes the residue field of $S$ at $s$. 
{\cred Unless otherwise specified, we consider $\kappa(s)$ as the base field of $X_s$. 
For example, we say that $X_s$ is {\em smooth} 
(resp. {\em geometrically reduced}) if $X_s$ is smooth (resp. geometrically reduced) over $\kappa(s)$. 
When $S$ is an integral scheme, the base field of the generic fibre 
is the function field $K(S)$ of $S$.  } 
\item 
For the definition of the relative dualising sheaf $\omega_{X/T}$, we refer to Remark \ref{r-dualising}. 
\item 
For the definition of conic bundles, see Subsection \ref{ss-conic-bdl-def}.
\item 
For the definition and basic properties on strictly henselian local rings, 
we refer to \cite[Subsection 2.8]{Fu15}. 
\item 
For a ring $A$, $\GL_n(A)$ denotes the group consisting of invertible matrices of size $n \times n$. 
Given a matrix $M$ of size $n \times n$, 
it is well known that $M$ is invertible if and only if $\det M \in A^{\times}$. 
\item 
Given $f \in A[x, y, z]$ and $M \in \GL_3(A)$, 
we define $f^M \in A[x, y, z]$ as follows. 
For $f = f(x, y, z)$, we set 
\[
f^M(x, y, z) := f( (M (x, y, z)^T)^T),
\]
where $(-)^T$ denotes the transposed matrix. 
More explicitly, for 
\[
M 
\begin{pmatrix}
x\\
y\\
z
\end{pmatrix}  = 
\begin{pmatrix}
m_{11} & m_{12} & m_{13}\\
m_{21} & m_{22} & m_{23}\\
m_{31} & m_{32} & m_{33}\\
\end{pmatrix} 
\begin{pmatrix}
x\\
y\\
z
\end{pmatrix} 
= 
\begin{pmatrix}
m_{11}x+ m_{12}y+ m_{13}z\\
m_{21}x+ m_{22}y+ m_{23}z\\
m_{31}x+ m_{32}y+ m_{33}z\\
\end{pmatrix},
\]
we have 
\[
f^M(x, y, z) = f(m_{11}x+ m_{12}y+ m_{13}z, m_{21}x+ m_{22}y+ m_{23}z, m_{31}x+ m_{32}y+ m_{33}z). 
\]
\item\label{ss-nota-conic} Set $\theta := 2m +1 \in \Z[m]$, which is the Hilbert polynomial of an arbitrary conic on $\P^2_k$ for a field $k$. 
\item 
Given a scheme $S$ and $S$-schemes $X$ and $Y$, 
we say that $X$ is {\em $S$-isomorphic} to $Y$ 
if there exists an isomorphism $\theta : X \to Y$ of schemes such that 
both $\theta : X \to Y$ and $\theta^{-1} : Y \to X$ are $S$-morphisms. 
When $S = \Spec\,A$ for a ring $A$, 
we say that $X$ is {\em $A$-isomorphic} to $Y$ if $X$ is $\Spec\,A$-isomorphic to $Y$. 
\end{enumerate}

\subsubsection{Singularities of minimal model program}\label{ss-mmp-sing}

We will freely use the standard notation in birational geometry, for which we refer to \cite{Kol13} and \cite{KM98}. 
Let $X$ be an integral normal excellent scheme admitting a dualising complex.  
We say that $X$ is {\em canonical} or {\em has at worst canonical singularities} if 
\begin{enumerate}
\item $K_X$ is $\Q$-Cartier, and 
\item all the coefficients $a_i$ are $\geq 0$ 
for every proper birational morphism $f: Y \to X$ from an integral normal excellent scheme  $Y$ 
and 
\[
K_Y =f^*K_X + \sum_{i=1}^r a_iE_i,
\]
where $E_1, ..., E_r$ are all the $f$-exceptional prime divisors. 
\end{enumerate}
Under assuming (1), 
if there exists a proper birational morphism $g : Z \to X$ from an integral regular excellent scheme $Z$, then (2) is known to be equivalent to (2)' below \cite[Lemma 2.30]{KM98}. 
\begin{enumerate}
\item[(2)'] all the coefficients $b_i$ are $\geq 0$ 
for some proper birational morphism $g: Z \to X$ from an integral regular excellent scheme  $Z$ 
and 
\[
K_Z =g^*K_X + \sum_{j=1}^s b_jF_j,
\]
where $F_1, ..., F_s$ are all the $g$-exceptional prime divisors. 
\end{enumerate}

\subsubsection{Linear algebra over local rings}


\begin{lem}\label{l-local-LA}
Let $A$ be a local ring. 
If $M \in \GL_n(A)$, then we can write $M = M_1 \cdots M_r$, where each $M_i$ is an elementary matrix. 
    Recall, e.g., that if $n=3$, then the list of elementary matrices are as follows: 
\begin{eqnarray*}
&&\begin{pmatrix}
\lambda & 0 & 0\\
0 & 1  & 0\\
0 & 0 & 1\\
\end{pmatrix} \qquad 
\begin{pmatrix}
1 & 0 & 0\\
0 & \lambda  & 0\\
0 & 0 & 1 \\
\end{pmatrix} \qquad 
\begin{pmatrix}
1 & 0 & 0\\
0 & 1  & 0\\
0 & 0 & \lambda \\
\end{pmatrix} \qquad (\lambda \in A^{\times})\\
&&\begin{pmatrix}
0 & 1 & 0\\
1 & 0  & 0\\
0 & 0 & 1\\
\end{pmatrix} \qquad 
\begin{pmatrix}
1 & 0 & 0\\
0 & 0  & 1\\
0 & 1 & 0 \\
\end{pmatrix} \qquad 
\begin{pmatrix}
0 & 0 & 1\\
0 & 1  & 0\\
1 & 0 & 0 \\
\end{pmatrix}\\
&&\begin{pmatrix}
1 & \mu & 0\\
0 & 1  & 0\\
0 & 0 & 1\\
\end{pmatrix} \qquad 
\begin{pmatrix}
1 & 0 & \mu\\
0 & 1  & 0\\
0 & 0 & 1\\
\end{pmatrix} \qquad 
\begin{pmatrix}
1 & 0 & 0\\
0 & 1  & \mu\\
0 & 0 & 1\\
\end{pmatrix}\\
&&\begin{pmatrix}
1 & 0 & 0\\
\mu & 1  & 0\\
0 & 0 & 1\\
\end{pmatrix} \qquad 
\begin{pmatrix}
1 & 0 & 0\\
0 & 1  & 0\\
\mu & 0 & 1\\
\end{pmatrix} \qquad 
\begin{pmatrix}
1 & 0 & 0\\
0 & 1  & 0\\
0 & \mu & 1\\
\end{pmatrix} \qquad (\mu \in A).
\end{eqnarray*}
\end{lem}

\begin{proof}
If $A$ is a field, then the assertion follows from linear algebra, 
i.e., $M$ becomes the identity matrix $E$ after applying elementary transformations 
finitely many times. 

The same argument works even when $A$ is a local ring. 
For example, we can find an entry $m_{ij}$ of $M$ which is a unit of $A$, 
where $m_{ij}$ denotes the $(i, j)$-entry of $M$. 
Then, applying suitable elementary transformations, we may assume that $m_{11} =1$. 
\end{proof}



\subsection{Definitions of conic bundles}\label{ss-conic-bdl-def}

In this subsection, we first introduce the definition of conic bundles 
(Definition \ref{d-conic-bdl}) and conics over a noetherian ring (Definition \ref{d-conicA}). 
We also summarise some basic properties, which should be well known to experts.

\subsubsection{Conic bundles}

\begin{dfn}\label{d-conic}
Let $\kappa$ be a field. 
\begin{enumerate}
\item 
We say that $C$ is {\em a conic on}  $\mathbb P^2_{\kappa}$ 
if  the equality  
\[
C = \Proj\,\kappa[x, y, z]/(ax^2 + by^2 + cz^2 + \alpha yz + \beta za + \gamma xy) 
\subset \P^2_{\kappa} = \Proj\,\kappa[x, y, z]
\]
holds for some $(a, b, c, \alpha, \beta, \gamma) \in \kappa^6 \setminus \{(0, 0, 0, 0, 0, 0)\}$. 
\item 
We say that $C$ is {\em a conic over} $\kappa$ 
if $C$ is $\kappa$-isomorphic to a conic on $\mathbb P^2_{\kappa}$. 
\end{enumerate}
\end{dfn}

\begin{dfn}\label{d-conic-bdl}
We say that $f: X \to T$ is a {\em conic bundle} if 
$f : X \to T$ is a flat proper morphism of noetherian schemes 
such that $X_t$ is a conic over $\kappa(t)$ for any point $t \in T$. 
\end{dfn}

\begin{rem}\label{r-conic-bdl}
If $f : X \to T$ is a conic bundle and $T' \to T$ is a morphism of noetherian schemes, then 
also the base change $X \times_T T' \to T'$ is a conic bundle. 
\end{rem}

\begin{lem}\label{l-conic-alg-fib}
Let $f :X \to T$ be a conic bundle. 
Then 
$f_*\MO_X = \MO_T$ and $R^if_*\MO_X =0$ for every $i>0$. 
\end{lem}

\begin{proof}
We have $H^i(X_t, \MO_{X_t})=0$ for every $i>0$. 
By \cite[Ch. III, Theorem 12.11]{Har77}, 
it holds that 
$R^if_*\MO_X=0$ for every $i>0$ and $f_*\MO_X$ is an invertible sheaf on $T$. 
Then $f^{\sharp}: \MO_T \to f_*\MO_X$ is an isomorphism, 
because 
the induced ring homomorphism 
\[
f^{\sharp} \otimes k(t) : \MO_T  \otimes k(t) \to (f_*\MO_X) \otimes k(t)
\]
is a nonzero $k(t)$-linear map of one-dimensional $k(t)$-vector spaces 
for every point $t \in T$. 
\end{proof}


Since a conic bundle $f: X \to T$ is a flat proper morphism of noetherian schemes whose fibres are Cohen--Macaulay of pure dimension one, 
there exists a dualising sheaf $\omega_{X/T}$ in the sense of \cite[page 157]{Con00}. 

\begin{rem}\label{r-dualising}
We here summarise some basic  properties on $\omega_{X/T}$ for later usage. 
\begin{enumerate}
\item 
Let $f: X \to T$ be a conic bundle of noetherian schemes. 
Since every fibre $X_t$ is Gorenstein of pure dimension one, 
 we have $\omega_{X/T} = \mathcal H^{-1}(f^!\MO_Y)$ and $\omega_{X/T}$ is an invertible sheaf on $X$. 
\item 
Let 
\[
\begin{CD}
X' @>\alpha >> X\\
@VVf'V @VVf V\\
T' @>\beta >> T
\end{CD}
\]
be a carterian diagram of noetherian schemes, where $f$ (and hence $f'$) is a conic bundle. 
We then obtain 
\[
\omega_{X'/T'} \simeq \alpha^*\omega_{X/T}. 
\]
In particular, $\omega_{X/T}|_{X_t} \simeq \omega_{X_t}$ for any point $t \in T$ and the fibre $X_t$ over $t$. 
\item 
Let $f: X \to T$ be a conic bundle, where $X$ and $T$ are Gorenstein normal varieties over an algebraically closed field. 
Then $\omega_{X/T} \simeq \MO_X(K_X -f^*K_T)$. 
\end{enumerate}
\end{rem}


\begin{prop}[{\cred cf. \cite{Bea77}*{Proposition 1.2}}]\label{p-conic-emb}
Let $f : X \to T$ be a conic bundle. 
Then the following hold. 
\begin{enumerate}
\item $R^if_*\omega_{X/T}^{-1} =0$ for every $i>0$. 
\item $f_*\omega_{X/T}^{-1} := f_*(\omega_{X/T}^{-1})$ is a locally free sheaf of rank $3$. 
\item $\omega_{X/T}^{-1}$ is very ample over $T$, and hence it defines a closed immersion 
$\iota: X \hookrightarrow \mathbb P( f_*\omega_{X/T}^{-1})$ over $T$. 
\end{enumerate}
\end{prop}

\begin{proof}
Let us show (1) and (2). 
Fix a point $t \in T$. 
By \cite[Ch. III, Corollary 12.9]{Har77}, it is enough to show that 
$\dim_{\kappa(t)} H^0(X_t, \omega^{-1}_{X/T}|_{X_t})=3$ 
and $H^i(X_t, \omega^{-1}_{X/T})=0$ for every $i>0$. 
This follows from 
$\omega^{-1}_{X/T}|_{X_t} \simeq \omega^{-1}_{X_t}$ and 
the fact that $X_t$ is a conic on $\mathbb P^2_{\kappa(t)}$. 
Thus (1) and (2) hold.

Let us show (3).  
First, we prove that $\omega_{X/T}^{-1}$ is $f$-free, 
  i.e., $f^*f_*\omega_{X/T}^{-1} \to \omega_{X/T}^{-1}$ is surjective. 
  Fix a closed point $x\in X$ and set $t:=f(x)$.
  Then $x\in X_t$, and there exists a section $\sigma\in H^0(X_t,\omega_{X/T}^{-1}|_{X_t})$ such that $\sigma|_{\{x\}}\neq 0$.
  We have two maps
  \begin{equation}\label{embed-1}
    f_*\omega_X^{-1}\to (f_*\omega_{X/T}^{-1})\otimes k(t) 
    \xrightarrow{\simeq}
    H^0(X_t,\omega_{X/T}^{-1}|_{X_t}).
  \end{equation}
  The first map is the projection, and the second map is obtained by \cite[Ch. III, Corollary 12.9]{Har77}.
  Let $T^\prime$ be an affine open subset containing $t$, and 
  set ${\cred X^\prime := f^{-1}(T')}$. 
  Then, by the above (\ref{embed-1}), we obtain 
  \begin{equation}\label{embed-2}
  H^0(X^\prime,\omega^{-1}_{X/T}|_{X^\prime})\to H^0(X^\prime,\omega^{-1}_{X/T}|_{X^\prime})\otimes k(t) \xrightarrow{\simeq} H^0(X_t,\omega_{X/T}^{-1}|_{X_t}).
  \end{equation}
  Since the first map is surjective and the second map is an isomorphism, 
   (\ref{embed-2}) is surjective. 
   Hence we can take $\tau\in H^0(X^\prime,\omega^{-1}_{X/T}|_{X^\prime})$ such that 
  $\tau|_{X_t} = \sigma$.
  Thus $\omega_{X/T}^{-1}$ is $f$-free. 
  Let $\iota : X \to \P(f_*\omega_{X/T}^{-1})$  be the induced morphism. 

  Next, we show that 
  $\iota : X \to \P(f_*\omega_{X/T}^{-1})$ is injective.
  Pick distinct points $x, x' \in X$ such that $\iota(x)=\iota(x') =:p$. 
Set $t:=f(x)=f(x')$.
  Then $x,x' \in X_t$, and hence we can take a $\sigma\in H^0(X_t,\omega_X^{-1}|_{X_t})$ such that $\sigma|_{\{x\}}=0$ and $\sigma|_{\{x'\}}\neq 0$. 
By the surjection (\ref{embed-2}), there exists $\tau \in H^0(X', \omega_{X/T}^{-1}|_{X'})$ 
such that $\tau|_{X_t} = \sigma$. 
This implies $\iota(x) \neq \iota(x')$, which is absurd. 

  Finally, we show that $\iota : X \to \P(f_*\omega_{X/T}^{-1}) =: {\cred P}$ is a closed immersion.
  Now, we have a commutative diagram 
  \begin{center}
    \begin{tikzcd}
      X_t \arrow[r,"\iota|_{X_t}"] \arrow[d,""] & P_t \arrow[d,""]\\
      X \arrow[r,"\iota"] & P,
    \end{tikzcd}
  \end{center}
  where $P_t$ is a fibre of $\pi\colon P\to T$.
  Since $X_t$ is a conic over $\kappa(t)$, we can embed this into $\mathbb{P}^2_{\kappa(t)}$, and this closed immersion is induced by $|-K_{X_t}|$. 
  By construction of $\iota$, 
  also $\iota|_{X_t}$ is a closed immersion. 
  Then we can check that $\iota$ is a closed immersion by applying Nakayama's lemma. 
  Thus (3) holds. 
\end{proof}

\subsubsection{Conics over rings}

\begin{dfn}\label{d-conicA}
Let $A$ be a noetherian ring. 
\begin{enumerate}
\item 
We say that $X$ is {\em a conic on}  $\mathbb P^2_A$ 
if  the equation 
\[
X = \Proj\,A[x, y, z]/(ax^2 + by^2 + cz^2 + \alpha yz + \beta za + \gamma xy) 
\subset \P^2_{\kappa} = \Proj\,A[x, y, z]
\]
holds for some $(a, b, c, \alpha, \beta, \gamma) \in A^6$ 
and the induced morphism $X \to \Spec\,A$ is a conic bundle. 
\item 
We say that $X$ is {\em a conic over} $A$ 
if $X$ is $A$-isomorphic to a conic on $\mathbb P^2_{A}$. 
\end{enumerate}
\end{dfn}

\begin{prop}\label{p-flat-local}
Let $A$ be a noetherian ring and take 
\[
Q := ax^2 + by^2 + cz^2 + \alpha yz + \beta zx + \gamma xy \in A[x, y, z]. 
\]
For $X := \Proj\,A[x, y, z]/(Q)$, the following are equivalent. 
\begin{enumerate}
\item $X$ is a conic on $\P^2_A$. 
\item The induced morphism $f: X \to \Spec\,A$ is a conic bundle. 
\item The induced morphism $f: X \to \Spec\,A$ is flat and any fibre of $f$ is one-dimensional. 
\end{enumerate}
Furthermore, if $(A, \m)$ is a noetherian local ring, 
then each of (1)--(3) is equivalent to (4). 
\begin{enumerate}
\item[(4)] {\cred At least} one of $a, b, c, \alpha, \beta, \gamma$ is not contained in $\m$. 
\end{enumerate}
\end{prop}

\begin{proof}
The implications $(1) \Rightarrow (2) \Rightarrow (3)$ are clear. 
Assume (3). 
Then $f: X \to \Spec\,A$ is a proper flat morphism. 
For any $t \in \Spec\,A$, its fibre $f^{-1}(t)$ is a conic over $k(t)$, as it is one-dimensional. 
Thus (1) holds. 

Assume that $(A, \m)$ is a local ring. 
Then it is obvious that (3) implies (4). 
Assume (4). 
We first reduce the problem to the case when $a \not\in \m$. 
By symmetry, we may assume that $\gamma \not\in \m$ and $a, b, c\in \m$. 
Applying the linear transform $(x, y, z) \mapsto (x, y+x, z)$, the problem is reduced to the case when $a \not\in \m$. 

Then the affine open subset $D_+(z)$ of $X$ can be written as 
\[
\Spec\,A[x, y]/(ax^2 +\varphi(y)x + \psi(y)), \qquad \varphi(y), \psi(y) \in A[y]. 
\]
We get an $A$-module isomorphism $A[x, y]/(ax^2 +\varphi(y)x + \psi(y)) \simeq A[y] \oplus x A[y]$, 
which is a free $A$-module. 
Therefore, $f: X \to \Spec\,A$ is flat. 
As the fibre $f^{-1}(\m)$ is one-dimensional, any fibre of $f$ is one-dimensional. 
Thus (3) holds. 
\end{proof}

\subsection{Local description of conic bundles}

Given a strictly henselian noetherian local ring $A$ and a conic $X = \Proj\,A[x, y, z]/(Q)$ on $\P^2_A$, 
the purpose of this subsection is to simplify the defining equation $Q$ via $A$-linear {\cred transformations}. 
In Subsection \ref{ss-local-desc-char0} and  Subsection \ref{ss-local-desc-char2}, 
we treat the case when the residue field of $A$ is of characteristic 
$\neq 2$ and $=2$, respectively.

\begin{lem}\label{l-local0}
Let $(A, \m, \kappa)$ be a noetherian local ring. 
Let 
\[
 X := {\rm Proj}\,A[x, y, z]/(Q).
\]
be a conic on $\P^2_A$ with $Q \in A[x, y, z]$. 
Then there exists $M \in \GL_3(A)$ such that 
\[
Q^M =a x^2 + by^2 + cz^2 + \alpha yz + \beta zx + \gamma xy 
\]
for some $a, b, c \in A^{\times}$ and $\alpha, \beta, \gamma \in A$. 
\end{lem}

\begin{proof}
By the same argument as in the second paragraph of the proof of Proposition \ref{p-flat-local}, 
we may assume that $a \in A^{\times}$. 
Applying a linear transform 
$(x, y, z) \mapsto {\cred (x + dx +ex, y, z)}$ for some $d, e \in \{0, 1\}$, 
we obtain $a \not\in \m, b \not\in \m$, and $c \not\in \m$. 
\end{proof}

\subsubsection{Local description in characteristic $\neq 2$}\label{ss-local-desc-char0}

\begin{prop}\label{p-char0-local0}
Let $\kappa$ be a field of characteristic $\neq 2$. 
Let 
\[
 X := {\rm Proj}\,\kappa[x, y, z]/(Q).
\]
be a conic on $\P^2_{\kappa}$ with $Q \in \kappa [x, y, z]$. 
Take the $3 \times 3$ symmetric matrix $S_Q$ such that 
$Q = (x\,y\,z)S_Q(x\,y\,z)^T$. 
Then the following hold. 
\begin{enumerate}
\item 
There exists $M \in \GL_3(A)$ such that 
\[
Q^M =a x^2 + by^2 + cz^2 
\]
for some $a, b, c \in \kappa$. 
\item 
The following are equivalent. 
\begin{enumerate}
\item $X$ is smooth. 
\item ${\rm rank}(S_Q) =3$. 
\end{enumerate}
\item 
The following are equivalent. 
\begin{enumerate}
\item $X$ is reduced. 
\item $X$ is geometrically reduced. 
\item ${\rm rank}(S_Q) =2$. 
\end{enumerate}
\item 
The following are equivalent. 
\begin{enumerate}
\item $X$ is not reduced. 
\item $X$ is not geometrically reduced. 
\item ${\rm rank}(S_Q) =1$. 
\end{enumerate}
\end{enumerate}
\end{prop}

\begin{proof}
All the assertions follow from linear algebra.  
\end{proof}

\begin{prop}\label{p-char0-local1}
Let $(A, \m, \kappa)$ be a noetherian local ring such that $\kappa$ is of characteristic $\neq 2$. 
Let 
\[
 X := {\rm Proj}\,A[x, y, z]/(Q).
\]
be a conic on $\P^2_A$ with $Q \in A[x, y, z]$. Let $0 \in \Spec\,A$ be the closed point. 
Then the following hold. 
\begin{enumerate}
\item 
There exists $M \in \GL_3(A)$ such that 
\[
Q^M =a x^2 + by^2 + cz^2 + \alpha yz. 
\]
for some $a \in A^{\times}$ and $b, c, \alpha \in A$. 
Furthermore, $X_0$ is not reduced if and only if $b, c, \alpha \in \m$. 
\item 
If $X_0$ is reduced, then 
there exists $M' \in \GL_3(A)$ such that 
\[
Q^{M'} =a' x^2 + b'y^2 + c'z^2 
\]
for some $a', b' \in A^{\times}$ and $c' \in A$. 
Furthermore, $X_0$ is not smooth if and only if $c' \in \m$. 
\end{enumerate}
\end{prop}

\begin{proof}
Note that  $2 \in A^{\times}$. 
Indeed, there exists $u \in A$ such that $2u \in 1+ \m$, which implies that $2 \in A^{\times}$.

Let us show (1).  
By Proposition \ref{p-flat-local},  
we can write 
\[
Q= ax^2 + by^2 + cz^2 + \alpha yz + \beta zx + \gamma xy, 
\]
where $a \in A^{\times}$ and $b, c, \alpha, \beta, \gamma \in A$. 
Completing the square  
\[
ax^2 +  \beta zx + \gamma xy = a \left(x+ \frac{\beta}{2a}z +\frac{\gamma }{2a}y\right)^2 - 
\left(\frac{\beta}{2a}z +\frac{\gamma }{2a}y\right)^2,
\]
we may assume that $\beta  = \gamma =0$: 
\[
Q =a x^2 + by^2 + cz^2 + \alpha yz. 
\]
By $a \in A^{\times}$ and Proposition \ref{p-char0-local0}, $X_0$ is not reduced if and only if $b, c, \alpha \in \m$. 
Thus (1) holds. 



Let us show (2). By (1), we may assume that 
\[
Q = a x^2 + by^2+cz^2 +\alpha yz
\]
for some $a \in A^{\times}$ and $b, c, \alpha \in A$. 
Since $X_0$ is reduced, one of $b, c, \alpha$ is contained in $A^{\times}$. 
We may assume that $b \in A^{\times}$ by applying $(x, y, z) \mapsto (x, y, z+y)$ for the case when 
$b, c \in \m$ and $\alpha \in A^{\times}$. 
Then we are done by completing the square again: 
\[
by^2+cz^2 +\alpha y z = b\left( y+\frac{\alpha}{2b}z\right)^2 + 
\left(c - \frac{\alpha^2}{4{\cred b^2}}\right) z^2.  
\]
Thus (2) holds. 
\end{proof}

\begin{cor}\label{c-char0-local2}
Let $(A, \m, \kappa)$ be a strictly henselian noetherian local ring such that $\kappa$ is of characteristic $\neq 2$. 
Let $0 \in  \Spec\,A$ be the closed point.  
Let 
\[
 X := {\rm Proj}\,A[x, y, z]/(Q).
\]
be a conic on $\P^2_A$ with $Q \in A[x, y, z]$. 
Then the following hold. 
\begin{enumerate}
\item 
There exists $M \in \GL_3(A)$ such that 
\[
Q^M =x^2 + by^2 + cz^2 + \alpha yz. 
\]
for some $b, c, \alpha \in A$. 
Furthermore, $X_0$ is not reduced if and only if $b, c, \alpha \in \m$. 
\item 
If $X_0$ is reduced, then there exists $M' \in \GL_3(A)$ such that 
\[
Q^{M'} =x^2 + y^2 + c'z^2 
\]
for some $c' \in A$. 
Furthermore, $X_0$ is not smooth over $\kappa$ if and only if $c' \in \m$. 
\item 
If $X_0$ is smooth over $\kappa$, then 
there exists $M'' \in \GL_3(A)$ such that 
\[
Q^{M''} = x^2 + y^2 + z^2. 
\]
\end{enumerate} 
\end{cor}

\begin{proof}
Since $A$ is a strictly henselian local ring, 
there exists $\sqrt{d} \in A$ for any element $d \in A^{\times}$. 
Then the assertions follow from Proposition \ref{p-char0-local1}. 
\end{proof}

\subsubsection{Local description in characteristic $2$}\label{ss-local-desc-char2}

\begin{prop}\label{p-char2-local}
Let $(A, \m, \kappa)$ be a noetherian local ring such that $\kappa$ is of characteristic two. 
Let 
\[
X := \Proj\,A[x, y, z]/(Q)
\]
be a conic on $\P^2_A$ with $Q \in A[x, y, z]$. 
Let $0 \in  \Spec\,A$ be the closed point.  
Then the following hold. 
\begin{enumerate}
\item 
If $X_0$ is geometrically reduced, then there exists $M \in \GL_3(A)$ such that 
\[
Q^M =  ax^2 + by^2 + cz^2 + yz + 2\beta zx + 2\gamma xy 
\]
for some 
$a, c, \beta, \gamma \in A$ and $b \in A^{\times}$. 
\item 
If $A$ is a strictly henselian local ring and  $X_0$ is geometrically reduced, 
then there exists $M' \in \GL_3(A)$ such that 
\[
Q^{M'} = a'x^2 +yz 
\]
for some $a' \in A$. 
In this case, $X_0$ is smooth if and only if $a' \in A^{\times}$. 
\end{enumerate}
\end{prop}

\begin{proof}
Let us show (1). 
We can write 
\[
Q =  ax^2 + by^2 + cz^2 + \alpha yz + \beta zx + \gamma xy 
\]
for some $(a, b, c, \alpha, \beta, \gamma) \in A^6 \setminus \{(0, ..., 0)\}$. 
We set $\overline{\alpha} := \alpha \mod \m, \overline{\beta} := \beta \mod \m, ...$ etc. 
As $X_0$ is geometrically reduced, we have $(\overline \alpha, \overline \beta, \overline \gamma) \neq (0, 0, 0)$. 
By symmetry, we may assume that $\overline{\alpha} \neq 0$, i.e., $\alpha \in A^{\times}$. 
Applying $(x, y, z) \mapsto (x, \alpha^{-1}y, z)$, 
the problem is reduced to the case when $\alpha =1$: 
 \[
 Q(x, y, z) = a x^2 +b y^2 +c z^2 +  yz + \beta zx + \gamma xy. 
 \]
We have 
\[
yz  + \beta zx + \gamma xy = (y + \beta x)(z +\gamma x) - \beta \gamma x^2. 
\]
Therefore, by applying $(x, y+\beta x, z + \gamma x) \mapsto (x, y, z)$, 
we may assume that $\alpha =1, \beta = 2\beta', \gamma=2\gamma'$ 
for some $\beta', \gamma' \in A$:  
\[
 Q(x, y, z) = a x^2 +b y^2 +c z^2 +  yz + 2\beta' zx + 2\gamma' xy. 
 \]
If $b \in A^{\times}$ or $c \in A^{\times}$, then 
we may assume $b \in A^{\times}$ by switching $y$ and $z$ if necessary. 
If $b \in \m$ and $c \in \m$. 
then we apply $(x, y, z) \mapsto (x, y, y+z)$: 
\begin{eqnarray*}
Q(x, y, y+z) 
&=& a x^2 +b y^2 +c (y+z)^2 +  y(y+z) + 2\beta' (y+z)x + 2\gamma' xy\\
&=& a x^2 +(b+c+1) y^2 +c z^2 +  (1+2c)yz + 2\beta' zx + 2(\beta'+\gamma') xy. 
\end{eqnarray*}
By applying $(1+2c)y \mapsto y$, 
we are done, because $1+b+c \in 1 + \m \subset A^{\times}$. 
Thus (1) holds. 

Let us show (2). 
By (1), we may assume that 
\[
Q(x, y, z) = ax^2 + by^2 + cz^2 + yz+ 2\beta zx + 2\gamma xy
\]
for some $a, c, \beta, \gamma  \in A$ and $b \in A^{\times}$. 
Since $\kappa$ is separably closed, we have 
\[
\overline{b}Y^2 + Y + \overline{c}  = \overline{b}(Y+\overline{\delta})(Y+ \overline{\epsilon})
\]
for some $\overline \delta, \overline \epsilon \in \kappa$ {\cred with $\overline \delta \neq \overline \epsilon$}. 
By Hensel's lemma, it holds that  
\[
bY^2 + Y + c = b(Y+ \delta)(Y+ \epsilon) \qquad \text{in} \qquad A[Y] 
\]
for some lifts $\delta, \epsilon \in A$ of $\overline \delta, \overline \epsilon \in \kappa$. 
Hence  we have $by^2 + yz + cz^2 = b(y+ \delta z)(y+ \epsilon z)$ and 
\[
Q(x, y, z) = ax^2 + by^2 + cz^2 + yz + 2\beta zx + 2\gamma xy= ax^2 + b(y+ \delta z)(y+ \epsilon z)+ 2\beta zx + 2\gamma xy. 
\]
Consider the matrix 
\[
\begin{pmatrix}
1 & \delta \\
1 & \epsilon
\end{pmatrix}. 
\]
This is an invertible matrix, because its determinant $\epsilon - \delta $ is contained in $A^{\times}$ 
(indeed, its reduction $\overline{\epsilon} -\overline{\delta} \in \kappa$ is nonzero). 
Applying $(x,  y+\delta z, y+\epsilon z) \mapsto (x, y, z)$, the problem is reduced to the case when 
$Q(x, y, z) = ax^2 +\alpha yz + 2\beta zx + 2\gamma xy$ for some $a, \beta, \gamma \in A$ and 
$\alpha \in A^{\times}$. 
By using $(x, y, z) \mapsto (x, y, \alpha^{-1}z)$, we may assume that $\alpha=1$: $Q(x, y, z) = ax^2 + yz + 2\beta zx + 2\gamma xy$. 
For $a' :=a-4\beta\gamma$, $y':= {\cred y+2\beta x}$ and $z':= z+2\gamma x$, we obtain 
\[
Q(x, y, z) = ax^2 + yz + 2\beta zx + 2\gamma xy 
= (a-4\beta\gamma)x^2+ (y+ 2\beta x)(z + 2\gamma x) = a'x^2 + y'z'. 
\]
Then it is clear that $X_0$ is smooth if and only if $a' \in A^{\times}$. 
Thus (2) holds. 
\end{proof}

\begin{prop}\label{p-char2-local2}
Let $\kappa$ be a separably closed field of characteristic two. 
Let 
\[
C := \Proj\,\kappa[x, y, z]/(Q)  
\]
be a conic on $\P^2_{\kappa}$. Then the following hold. 
\begin{enumerate}
 \item 
 There exists $M \in \GL_3(\kappa)$ such that one of the following holds.  
 \begin{enumerate}
 \item $Q^M = ax^2 + by^2 + cz^2$ for some $a, b, c \in \kappa$. 
 In this case, $C$ is not geometrically reduced. 
 \item $Q^M = ax^2 + yz$ for some $a \in \kappa$. 
 In this case, $C$ is geometrically reduced. Furthermore, $C$ is smooth if and only if $a\neq 0$. 
 \end{enumerate}
\item 
 If $\kappa$ is algebraically closed, then there exists $M' \in \GL_3(\kappa)$ such that one of the following holds.  
 \begin{enumerate}
 \item $Q^{M'} = x^2$. In this case, $C$ is not reduced. 
 \item $Q^{M'} = yz$. In this case, $C$ is not smooth but {\cred is} reduced. 
 \item $Q^{M'} = x^2 + yz$. In this case, $C$ is smooth. 
 \end{enumerate}
\end{enumerate}
\end{prop}

\begin{proof}
Let us show (1). 
If $X$ is geometrically reduced, then the assertion follows from Proposition \ref{p-char2-local}(2). 
We may assume that $C$ is not geometrically reduced. Then $\alpha = \beta = \gamma =0$. 
Thus (1) holds. 
The assertion (2) immediately follows from (1). 
\qedhere
 
\end{proof}

\subsection{Generic smoothness}

\begin{dfn}
Let $f: X \to T$ be a conic bundle such that $T$ is a noetherian integral scheme. 
We say that $f$ is {\em generically smooth} if the generic fibre of $f$ is smooth. 
We say that $f$ is {\em wild} if $X$ is normal and  $f$ is not generically smooth. 
\end{dfn}

{\cred 
For a conic bundle $f: X \to T$ such that $T$ is a noetherian integral scheme, 
$f$ is generically smooth if and only if there exists a non-empty open subset $T'$ of $T$ such that $f^{-1}(t)$ is smooth for every point $t \in T'$ 
\cite[Th\'eor\`eme 12.2.4]{EGAIV3}.} 

\begin{lem}\label{l-char0-gene-sm}
Let $f: X \to T$ be a conic bundle such that $T$ is a noetherian integral scheme. 
 If $X$ is normal and the function field $K(T)$ of $T$ is of characteristic $\neq 2$, 
then $f$ is generically smooth. 
\end{lem}

\begin{proof}
Taking the base change by $\Spec\,K(T) \to T$, we may assume that $T = \Spec\,k$ for a field $k$. 
Then $X$ is a conic on $\P^2_{k}$ {\cred and $X$ is a regular scheme}. 
Taking the base change to the separable closure of $k$, we may assume that $k$ is separably closed. 
By Proposition \ref{p-char0-local0}(3) or Corollary \ref{c-char0-local2}(2), we can write 
$X \simeq \{ x^2 +y^2 + cz^2 =0\} \subset \P^2_{k}$ for some $c \in k$. 
If $c=0$, then $X$ is not regular, as the $k$-rational point $[0:0:1] \in \P^2_k$ is a non-regular point. 
Hence $c \neq 0$. In this case, $X$ is smooth over $T=\Spec\,k$. 
\end{proof}

\begin{ex}\label{e-char2-wild-cb}
We work over an algebraically closed field of characteristic two. 
Set 
\[
X := \{ t_0x_0^2+t_1 x_1^2+ t_2x_2^2  =0\} \subset \P^2_{x_0, x_1, x_2} \times \P^2_{t_0, t_1, t_2}. 
\]
Then the induced morphism $f: X \to T := \P^2_{t_0, t_1, t_2} = \Proj\,k[t_0, t_1, t_2]$ 
is a wild conic bundle, because 
we have $a_0x_0^2+a_1 x_1^2+ a_2x_2^2 = (\sqrt{a_0} x_0 + \sqrt{a_1} x_1 + \sqrt{a_2}x_2)^2$ 
for all $a_0, a_1, a_2 \in k$. 
\end{ex}

\section{Discriminant divisors}\label{s-disc-Sigma}

Let $f: X \to T$ be a conic bundle. 
The purpose  of this section is to introduce two closed subschemes $\Delta_f$ and $\Sigma_f$ of $T$ 
which satisfy the following properties (I) and (II). 
\begin{enumerate}
    \item[(I)] Given a point $t \in T$, 
    $X_t$ is not smooth (resp. not geometrically reduced) if and only if  $t \in \Delta_f$ 
    (resp. $t \in \Sigma_f$) (Theorem \ref{t-disc-non-sm}). 
    \item[(II)] Equalities $\Delta_{f'} = \beta^{-1}\Delta_f$ and $\Sigma_{f'} = \beta^{-1}\Sigma_f$ of closed subschemes 
    hold for a carterian diagram of noetherian schemes
    \[
    \begin{CD}
    X' @>\alpha >> X\\
    @VVf' V @VVf V\\
    T' @>\beta >> T, 
    \end{CD}
    \]
    where $f$ and $f'$ are conic bundles (Remark \ref{r-disc-bc}, Remark \ref{r-Sigma-bc}).   
\end{enumerate}
In Subsection \ref{ss-def-disc} (resp. Subsection \ref{ss-def-nonred-locus}), 
we give the definition of $\Delta_f$ (resp. $\Sigma_f$). 
The main technical difficulty is their well-definedness. 
The property (I) will be established in Subsection \ref{ss-Delta-Sigma-loci}. 
In Subsection \ref{ss-Sigma-char2}, we introduce a simpler version $\Sigma'_f$ of $\Sigma_f$ which can be defined only in characteristic two.

\subsection{Definition of discriminant divisors}\label{ss-def-disc}

In this subsection, we introduce a closed subscheme $\Delta_f$ of $T$, called the discriminant scheme of $f$, associated with a conic bundle $f : X \to T$. 
We first define $\Delta_f$ for the case when $T$ is affine, say $T = \Spec\,A$, 
and a closed embedding 
\[
X =\Proj\,A[x, y, z]/(Q) \subset \P^2_A
\]
is fixed (Definition \ref{d-delta}). 
In this case, we set $\Delta_f := \Spec\,(A/\delta(Q)A)$, 
where $\delta(Q) \in A$ is determined by $Q$. 
Then we shall check that $\Delta_f = \Spec\,(A/\delta(Q)A)$ is independent of the choice of 
the embedding $X \subset \P^2_A$ (Proposition \ref{p-delta-welldef}), 
which enables us to define $\Delta_f$ for the general case, i.e., 
we may glue $\{\Delta_{f_i}\}_{i \in I}$ together 
for a suitable affine open cover $T = \bigcup_{i \in I} T_i$ and the induced morphisms $f_i : f^{-1}(T_i) \to T_i$ (Definition \ref{d-disc}).

\begin{dfn}[{\cred cf.\cite{ABBB21}*{Section 2.b}}]\label{d-delta}
Let $A$ be a noetherian ring. 
For 
\[
Q:= ax^2 + by^2 +cz^2 + \alpha yz + \beta zx + \gamma xy \in A[x, y, z], 
\]
we set 
\[
\delta(Q) := 4abc + \alpha \beta \gamma - a \alpha^2 - b\beta^2 -c \gamma^2 \in A. 
\] 
\end{dfn}

{\cred 
\begin{rem}
Note that if $2 \in A^{\times}$, then 
\[
\delta(Q) = \frac{1}{2} \det 
\begin{pmatrix}
2a & \gamma & \beta\\
\gamma & 2b & \alpha\\
\beta & \alpha & 2c
\end{pmatrix}. 
\]
\end{rem}
}

\begin{prop}\label{p-delta-welldef}
Let $A$ be a noetherian ring. 
For  $M \in \GL_3(A)$ and 
\[
Q := ax^2 + by^2 +cz^2 + \alpha yz + \beta zx + \gamma xy \in A[x, y, z], 
\]
the equality 
\[
\delta(Q)A = \delta(Q^M)A
\]
of ideals hold. 
\end{prop}

\begin{proof}
We first reduce the problem to the case when $A$ is a noetherian local ring. 
For $J := \delta(Q)A \cap \delta(Q^M) A$, we have the inclusion 
\[
\varphi : J = \delta(Q)A \cap \delta(Q^M) A \hookrightarrow \delta(Q) A, 
\]
which is an $A$-module homomorphism. 
For a prime ideal $\p \in \Spec\,A$, we get the following inclusion \cite[Corollary 3.4, ii)]{AM69}: 
\[
\varphi_{\p} : J_{\p} = 
\delta(Q) A_{\p} \cap \delta(Q^M) A_{\p} \hookrightarrow \delta(Q) A_{\p}. 
\]
As we are assuming that $\varphi_{\p}$ is an isomorphism for every $\p \in \Spec\,A$, 
also $\varphi$ is an isomorphism, i.e., $\delta(Q) A \subset \delta(Q^M)A$. 
By symmetry, we obtain the opposite inclusion $\delta(Q) A \supset \delta(Q^M)A$, 
which implies $\delta(Q) A = \delta(Q^M)A$. 
Thus the problem is reduced to the case when  $A$ is a noetherian local ring.

Since $A$ is a local ring, we have $M=M_1 \cdots M_r$ for some elementary matrices $M_1, ..., M_r$ 
(Lemma \ref{l-local-LA}). 
By symmetry, we may assume that one of (1)--(3) holds.  
\begin{enumerate}
\item 
$M = 
\begin{pmatrix}
\lambda & 0 & 0\\
0 & 1 & 0\\
0 & 0 & 1\\
\end{pmatrix}$ for some $\lambda \in A^{\times}$. 
\item 
$M = 
\begin{pmatrix}
0 & 1 & 0\\
1 & 0 & 0\\
0 & 0 & 1\\
\end{pmatrix}$. 
\item 
$M = 
\begin{pmatrix}
1 & \mu & 0\\
0 & 1 & 0\\
0 & 0 & 1\\
\end{pmatrix}$ for some $\mu \in A$. 
\end{enumerate}
Set $Q(x, y, z) := Q$. 

(1) In this case, we have 
\[
Q^M(x, y, z) = Q(\lambda x, y, z) = 
(\lambda^2a) x^2 + by^2 +cz^2 + \alpha yz + (\lambda\beta) zx + (\lambda\gamma) xy
\]
and 
\[
\delta(Q^M) = 4(\lambda^2a) bc + \alpha (\lambda \beta)(\lambda \gamma) - (\lambda ^2a) \alpha^2 - b(\lambda \beta)^2 -c (\lambda \gamma)^2 = \lambda^2 \delta(Q), 
\]
as required.

(2) In this case, 
we obtain $Q^M(x, y, z) = Q(y, x, z)$, which implies $\delta(Q^M) =\delta(Q)$ by symmetry 
(Definition \ref{d-delta}). 

(3) In this case, we have 
\begin{eqnarray*}
Q^M(x, y, z) &=& Q(x+\mu y, y, z)\\
&=& a(x+\mu y)^2 + by^2 +cz^2 + \alpha yz + \beta z(x+\mu y) + \gamma (x+\mu y)y\\
&=& ax^2  + (a\mu^2 + b + \gamma \mu)y^2 +cz^2
 + (\alpha + \beta \mu) yz + \beta  zx+(2a\mu + \gamma )xy.
\end{eqnarray*}
Therefore, the assertion holds by the following computation: 
\begin{eqnarray*}
\delta(Q^M) &=& 
4a (a\mu^2 + b + \gamma \mu) c + 
(\alpha + \beta \mu) \cdot \beta  \cdot (2a\mu + \gamma )\\
&-&a (\alpha + \beta \mu)^2 - (a\mu^2 + b + \gamma \mu) \beta^2 
-c (2a\mu + \gamma )^2\\
&=& (4a^2c \mu^2 + 4ac \gamma \mu + 4abc) \\
&+& (2a \beta^2 \mu^2 + (2a\alpha\beta + \beta^2 \gamma)\mu + a\beta \gamma)\\
&+& (-a\beta^2 \mu^2 - 2a\alpha\beta\mu - a\alpha^2)\\
&+& (-a\beta^2 \mu^2 -  \beta^2\gamma \mu - b\beta^2)\\
&+& (-4a^2 c\mu^2 -4ac\gamma\mu -c\gamma^2)\\
&=&  4abc + \alpha \beta \gamma - a \alpha^2 - b\beta^2 -c \gamma^2\\
&=&\delta(Q). 
\end{eqnarray*}
\end{proof}

\begin{dfn}\label{d-disc}
Let $f: X \to T$ be a conic bundle. 
We define the closed subscheme $\Delta_f$ of $T$, called {\em the discriminant locus} or  
{\em discriminant scheme} of $f$, as follows. 
\begin{enumerate}
\item If $T = \Spec\,A$ and $X =\Proj\,A[x, y, z]/(Q)$, then 
we set 
\[
\Delta_f := \Spec\,(A/\delta(Q)A). 
\]
By Proposition \ref{p-delta-welldef}, this definition is independent of the choice of $Q$. 
\item 
Fix an affine open cover $T = \bigcup_{i \in I} T_i$ 
such that 
$f^{-1}(T_i)$ is $T_i$-isomorphic to a conic on $\mathbb P^2_{T_i}$ for every $i \in I$ 
(the existence of such an open cover is guaranteed by Proposition \ref{p-conic-emb}). 
We can write 
\[
X_i := f^{-1}(T_i) = \Proj\, A_i[x, y, z]/(Q_i), \qquad Q_i \in A_i[x, y, z], \qquad 
A_i := \Gamma(T_i, \MO_T). 
\]
By (1), we have a closed subscheme $\Delta_{f_i}$ of $T_i$ for every $i \in I$. 
Again by (1), we obtain $\Delta_{f_i}|_{T_i \cap T_j} = \Delta_{f_j}|_{T_i \cap T_j}$, 
so that 
there exists a closed subscheme $\Delta_f$ on $T$ such that $\Delta_f|_{T_i} = \Delta_{f_i}$ 
for every $i \in I$. 
It follows from (1) that $\Delta_f$ does not depend on the choice of the affine open cover $T = \bigcup_{i \in I} T_i$. 
\end{enumerate}
We also call $\Delta_f$ the {\em discriminant divisor} of $f$ 
when $\Delta_f$ is an effective Cartier divisor on $T$. 
\end{dfn}

\begin{rem}\label{r-disc-bc}
Let 
\[
\begin{CD}
X' @>\alpha >> X\\
@VVf'V @VVfV\\
T' @> \beta >> T
\end{CD}
\]
be a cartesian diagram of noetherian schemes, 
where $f$ is a conic bundle. 
In particular, also $f'$ is a conic bundle. 
In this case, the equality 
\[
\beta^{-1}(\Delta_f) = \Delta_{f'}
\]
holds by Definition \ref{d-disc}, where $\beta^{-1}(\Delta_f)$ is the scheme-theoretic inverse image of $\Delta_f$. 
\end{rem}

\subsection{Locus of non-reduced fibres}\label{ss-def-nonred-locus}


In this subsection, we introduce a closed subscheme $\Sigma_f$ of $T$ associated with a conic bundle $f : X \to T$. 
The outline is similar to that of Subsection \ref{ss-def-nonred-locus}, whilst its proof is more involved, because the defining ideal $\sigma(Q)$ is no longer principal. 
The proof is carried out by reducing the problem to the case of characteristic zero.

\begin{dfn}\label{d-sigma}
Let $A$ be a noetherian ring. 
For 
\[
Q(x, y, z) := ax^2 + by^2 +cz^2 + \alpha yz + \beta zx + \gamma xy \in A[x, y, z], 
\]
we set 
\[
\sigma(Q) 
:= (4ab -\gamma^2, 4bc - \alpha^2, 4ca - \beta^2, 
2a \alpha - \beta \gamma, 2b\beta - \gamma \alpha, 2c\gamma - \alpha \beta) 
\subset A, 
\] 
which is an ideal of $A$. 
\end{dfn}

{\cred 
\begin{rem}
Note that $\sigma(Q)$ is the ideal of $A$ generated by the $2 \times 2$ minors of 
\[
\begin{pmatrix}
2a & \gamma & \beta\\
\gamma & 2b & \alpha\\
\beta & \alpha & 2c
\end{pmatrix}. 
\]
\end{rem}
}

The following lemma is the key to prove that $\sigma(Q)$ does not depend on the choice of 
the closed embedding into $\P^2_A$.

\begin{lem}\label{l-sigma-welldef}
Let $A$ be a noetherian ring and fix $\mu \in A$.  
Set 
\[
Q(x, y, z) := ax^2 + by^2 +cz^2 + \alpha yz + \beta zx + \gamma xy \in A[x, y, z] 
\]
and 
\[
M := 
\begin{pmatrix}
1 & 0 & \mu \\
0 & 1 & 0\\
0 & 0 & 1
\end{pmatrix}. 
\]
Then the equality  
\begin{equation}\label{e1-sigma-welldef}
\sigma(Q) = \sigma(Q^M). 
\end{equation}
of ideals of $A$ holds. 
\end{lem}

\begin{proof}
We introduce the following notation. 
\begin{itemize}
\item $\sigma_{11}(Q) := 4bc -\alpha^2$. 
\item $\sigma_{22}(Q) := 4ca -\beta^2$. 
\item $\sigma_{33}(Q) := 4ab -\gamma^2$. 
\item $\sigma_{12}(Q) := 2c\gamma - \alpha \beta$. 
\item $\sigma_{23}(Q) := 2a\alpha - \beta \gamma$. 
\item $\sigma_{31}(Q) := 2b\beta - \gamma \alpha$. 
\end{itemize}

\setcounter{step}{0}

\begin{step}\label{s1-sigma-welldef}
The equality (\ref{e1-sigma-welldef}) holds if the following equalities (1)--(6) hold. 
\begin{enumerate}
\item $\sigma_{11}(Q^M) = \sigma_{11}(Q) - 2\mu \sigma_{31}(Q) +\mu^2 \sigma_{33}(Q)$. 
\item $\sigma_{22}(Q^M) = \sigma_{22}(Q)$. 
\item $\sigma_{33}(Q^M) = \sigma_{33}(Q)$.  
\item $\sigma_{12}(Q^M) = \sigma_{12}(Q) -\mu \sigma_{23}(Q)$. 
\item $\sigma_{23}(Q^M) = \sigma_{23}(Q)$. 
\item $\sigma_{31}(Q^M) = \sigma_{31}(Q) + \mu \sigma_{33}(Q)$. 
\end{enumerate}
\end{step}

\begin{proof}[Proof of Step \ref{s1-sigma-welldef}]
Note that $\{ \sigma_{ij}(Q) \}_{i, j}$ and $\{ \sigma_{ij}(Q^M) \}_{i, j}$ 
are  generators of $\sigma(Q)$ and $\sigma(Q^M)$, respectively: 
\begin{eqnarray*}
\sigma(Q) &=& (\sigma_{11}(Q), \sigma_{22}(Q), \sigma_{33}(Q), 
\sigma_{12}(Q), \sigma_{23}(Q), \sigma_{31}(Q))\\
\sigma(Q^M) &=& (\sigma_{11}(Q^M), \sigma_{22}(Q^M), \sigma_{33}(Q^M), 
\sigma_{12}(Q^M), \sigma_{23}(Q^M), \sigma_{31}(Q^M)). 
\end{eqnarray*}
Therefore, (1)--(6) imply the equality (\ref{e1-sigma-welldef}). 
This completes the proof of Step \ref{s1-sigma-welldef}. 
\end{proof}

\begin{step}\label{s2-sigma-welldef}
If $A$ is a field of characteristic zero, then the equalities (1)--(6) of Step \ref{s1-sigma-welldef} hold. 
\end{step}

\begin{proof}[Proof of Step \ref{s2-sigma-welldef}]
For the symmetric matrix
\[
S := \begin{pmatrix}
2a & \gamma & \beta \\
\gamma & 2b & \alpha \\
\beta & \alpha  & 2c
\end{pmatrix}, 
\]
we have 
\[
Q(x, y, z) = ax^2 + by^2 +cz^2 + \alpha yz + \beta zx + \gamma xy 
= 
\frac{1}{2}
\begin{pmatrix}
x & y & z
\end{pmatrix} 
S
\begin{pmatrix}
x \\
y \\
z
\end{pmatrix}. 
\]
Furthermore, the following equalities (i)--(vi) hold if 
each $S_{ij}$ denotes the $2 \times 2$-matrix 
obtained from $S$ by excluding the $i$-th row and the $j$-th column. 
\begin{enumerate}
\renewcommand{\labelenumi}{(\roman{enumi})}
\item $\sigma_{11}(Q) = 4bc -\alpha^2 = \det S_{11}$. 
\item $\sigma_{22}(Q) = 4ca -\beta^2  = \det S_{22}$. 
\item $\sigma_{33}(Q) = 4ab -\gamma^2 = \det S_{33}$. 
\item $\sigma_{12}(Q) = 2c\gamma - \alpha \beta = \det S_{12} =\det S_{21}$. 
\item $\sigma_{23}(Q) = 2a\alpha - \beta \gamma = \det S_{23} =\det S_{32}$. 
\item $\sigma_{31}(Q) = 2b\beta - \gamma \alpha = -\det S_{31}  =-\det S_{13}$. 
\end{enumerate}
It holds that 
\[
Q^M(x, y, z) = \frac{1}{2}
\begin{pmatrix}
x & y & z
\end{pmatrix} 
(M^TSM)
\begin{pmatrix}
x \\
y \\
z
\end{pmatrix}, 
\]
and 
\[
M^TSM = 
\begin{pmatrix}
2a & \gamma & \beta +\mu \cdot (2a) \\
\gamma & 2b & \alpha + \mu \gamma \\
\beta +\mu \cdot (2a) & \alpha + \mu \gamma  & 2c +2\mu \beta + \mu^2 \cdot (2a)
\end{pmatrix}
\]
where $M^T$ denotes the transposed matrix of $M$. 

(1) 
\begin{eqnarray*}
\sigma_{11}(Q^M) 
&=& \det (M^TSM)_{11} \\
&=& \det 
\begin{pmatrix}
2b & \alpha + \mu \gamma \\
\alpha + \mu \gamma  & 2c +2\mu \beta + \mu^2 \cdot (2a)
\end{pmatrix}\\
&=&  \det 
\begin{pmatrix}
2b & \alpha + \mu \gamma \\
\alpha   & 2c +\mu \beta
\end{pmatrix}
+ \mu \det 
\begin{pmatrix}
2b & \alpha + \mu \gamma \\
\gamma   & \beta +\mu \cdot (2a)
\end{pmatrix}\\
&=&  
\det 
\begin{pmatrix}
2b & \alpha\\
\alpha   & 2c
\end{pmatrix}
+\mu 
\det 
\begin{pmatrix}
2b & \gamma \\
\alpha   & \beta
\end{pmatrix}
+ \mu \det 
\begin{pmatrix}
2b & \alpha\\
\gamma   & \beta
\end{pmatrix}
+ \mu^2 \det 
\begin{pmatrix}
2b & \gamma \\
\gamma   & 2a
\end{pmatrix}\\
&=& \det S_{11} - \mu \det S_{13} -\mu \det S_{31} +\mu^2 \det S_{33}\\
&=& \sigma_{11}(Q) - 2\mu \sigma_{31}(Q) +\mu^2 \sigma_{33}(Q). 
\end{eqnarray*}

(2) 
\[
\sigma_{22}(Q^M) = \det (M^TSM)_{22} = \det 
\begin{pmatrix}
2a & \beta +\mu \cdot (2a) \\
\beta +\mu \cdot (2a) & 2c +2\mu \beta + \mu^2 \cdot (2a)
\end{pmatrix}
\]
\[
=
\det 
\begin{pmatrix}
2a & \beta +\mu \cdot (2a) \\
\beta  & 2c +\mu \beta 
\end{pmatrix} 
= 
\det 
\begin{pmatrix}
2a & \beta\\
\beta  & 2c
\end{pmatrix} = \det S_{22} = \sigma_{22}(Q). 
\]

(3) $\sigma_{33}(Q^M) = \det (M^TSM)_{33} =\det S_{33} = \sigma_{33}(Q)$. 

(4) 
\begin{eqnarray*}
\sigma_{12}(Q^M) &=&\det (M^TSM)_{12} \\
&=&\det \begin{pmatrix}
\gamma &  \alpha + \mu \gamma\\
\beta +\mu \cdot (2a)  & 2c +2\mu \beta + \mu^2 \cdot (2a)
\end{pmatrix}\\
&=&\det \begin{pmatrix}
\gamma &  \alpha \\
\beta +\mu \cdot (2a)  & 2c +\mu \beta 
\end{pmatrix}\\
&=& \det \begin{pmatrix}
\gamma & \alpha \\
\beta   & 2c 
\end{pmatrix}
+ \mu 
\det \begin{pmatrix}
\gamma & \alpha \\
2a  & \beta
\end{pmatrix}\\
&=& \det S_{12} - \mu \det S_{32}\\
&=& \sigma_{12}(Q) - \mu \sigma_{23}(Q). 
\end{eqnarray*}

(5) 
\[
\sigma_{23}(Q^M) =\det (MSM^T)_{23} =
\det 
\begin{pmatrix}
2a & \gamma \\
\beta +\mu \cdot (2a) & \alpha + \mu \gamma
\end{pmatrix} = \det \begin{pmatrix}
2a & \gamma \\
\beta & \alpha
\end{pmatrix}
\]
\[
=\det S_{23} =  \sigma_{23}(Q).
\]

(6) 
\[
\sigma_{31}(Q^M) = -\det (M^TSM)_{13} 
=-\det 
\begin{pmatrix}
\gamma & 2b \\
\beta +\mu \cdot (2a) & \alpha + \mu \gamma
\end{pmatrix}
\]
\[
= -\det \begin{pmatrix}
\gamma & 2b \\
\beta & \alpha
\end{pmatrix}
-\mu \det \begin{pmatrix}
\gamma & 2b \\
2a & \gamma
\end{pmatrix}
= -\det S_{13} + \mu \det S_{33} =\sigma_{31}(Q) +\mu \sigma_{33}(Q).
\]
This completes the proof of Step \ref{s2-sigma-welldef}. 
\end{proof}

\begin{step}\label{s3-sigma-welldef}
The equalities (1)--(6) of Step \ref{s1-sigma-welldef} hold without any additional assumptions. 
\end{step}

\begin{proof}[Proof of Step \ref{s3-sigma-welldef}]
Let $A' := \Z [a', b', c', \alpha', \beta', \gamma', \mu']$ 
be the polynomial ring over $\Z$ with 
seven variables $a', b', c', \alpha', \beta', \gamma', \mu'$. 
Let 
\[
\varphi :  A' = \Z [a', b', c', \alpha', \beta', \gamma', \mu'] \to A
\]
be the ring homomorphism such that 
\[
\varphi(a') = a,\quad \varphi(b')=b,\quad \varphi(c') = c, \quad \varphi(\alpha') = \alpha, 
\quad \varphi(\beta') = \beta, \quad \varphi(\gamma') = \gamma, \quad \varphi(\mu') = \mu. 
\]
Take the embedding $\psi : A' \hookrightarrow A'' := \Frac\,A'$ into its total ring of fractions 
$A'' = \Frac\,A'$. 
Set 
\[
Q'(x, y, z) := a'x^2 + b'y^2 +c'z^2 + \alpha' yz + \beta' zx + \gamma' xy \in A'[x, y, z] 
\subset A''[x, y, z]. 
\]
and \[
M' := 
\begin{pmatrix}
1 & 0 & \mu' \\
0 & 1 & 0\\
0 & 0 & 1
\end{pmatrix}. 
\]

Let us only show that (1) of Step \ref{s1-sigma-welldef} holds, as the other ones (2)--(6) are similar. 
Applying Step \ref{s2-sigma-welldef} to $A''$, we obtain an equality 
\[
\sigma_{11}(Q'^{M'}) = \sigma_{11}(Q') - 2\mu' \sigma_{31}(Q') +\mu'^2 \sigma_{33}(Q')
\quad {\rm in} \quad A''. 
\]
Since both sides are elements of $A'$, it follows from $A' \subset A''$ that 
\[
\sigma_{11}(Q'^{M'}) = \sigma_{11}(Q') - 2\mu' \sigma_{31}(Q') +\mu'^2 \sigma_{33}(Q')
\quad {\rm in} \quad A'.
\]
Take the image by  $\varphi: A' \to A$: 
\[
\varphi(\sigma_{11}(Q'^{M'})) = \varphi(\sigma_{11}(Q')) - 2\mu \varphi(\sigma_{31}(Q')) +
\mu^2 \varphi(\sigma_{33}(Q')).  
\]
We can check that 
$\varphi(\sigma_{11}(Q'^{M'})) = \sigma_{11}(Q^M)$ and $\varphi(\sigma_{ij}(Q')) = \sigma_{ij}(Q)$ 
for all $i, j$. 
This completes the proof of Step \ref{s3-sigma-welldef}. 
\end{proof}
Step \ref{s1-sigma-welldef} and Step \ref{s3-sigma-welldef} complete the proof of 
Lemma \ref{l-sigma-welldef}. 
\end{proof}

\begin{prop}\label{p-sigma-welldef}
Let $A$ be a noetherian ring. 
For  $M \in \GL_3(A)$ and 
\[
Q(x, y, z) := ax^2 + by^2 +cz^2 + \alpha yz + \beta zx + \gamma xy \in A[x, y, z], 
\]
the equality  
\[
\sigma(Q) = \sigma(Q^M)
\]
of ideals hold. 
\end{prop}

\begin{proof}
By the same argument as in  the proof of Proposition \ref{p-delta-welldef}, 
we may assume that $A$ is a noetherian local ring and one of (1)-(3) holds.  
\begin{enumerate}
\item 
$M = 
\begin{pmatrix}
\lambda & 0 & 0\\
0 & 1 & 0\\
0 & 0 & 1\\
\end{pmatrix}$ for some $\lambda \in A^{\times}$. 
\item 
$M = 
\begin{pmatrix}
0 & 1 & 0\\
1 & 0 & 0\\
0 & 0 & 1\\
\end{pmatrix}$. 
\item 
$M = 
\begin{pmatrix}
1 & 0 & \mu\\
0 & 1 & 0\\
0 & 0 & 1\\
\end{pmatrix}$ for some $\mu \in A$. 
\end{enumerate}

(1) In this case, we have 
\[
Q^M(x, y, z) = Q(\lambda x, y, z) = 
(\lambda^2a) x^2 + by^2 +cz^2 + \alpha yz + (\lambda\beta) zx + (\lambda\gamma) xy
\]
and  
\[
\sigma(Q^M) = (4\lambda^2 ab -(\lambda\gamma)^2, 4bc - \alpha^2, 4c(\lambda^2a) - (\lambda \beta)^2, 
\]
\[
2(\lambda^2a)\alpha -(\lambda \beta) (\lambda \gamma), 
2b(\lambda\beta) - (\lambda \gamma) \alpha, 
2c(\lambda\gamma) - \alpha(\lambda \beta)) 
=\sigma(Q), 
\]
as required.

(2) In this case, 
we have $\sigma(Q^M) =\sigma(Q)$ by symmetry.

(3) In this case, 
the equality $\sigma(Q^M) =\sigma(Q)$
follows from Lemma \ref{l-sigma-welldef}. 
\end{proof}

\begin{dfn}\label{d-Sigma}
Let $f: X \to T$ be a conic bundle. 
We define the closed subscheme $\Sigma_f$ of $T$ as follows. 
\begin{enumerate}
\item If  $T=\Spec\,A$ and $X =\Proj\,A[x, y, z]/(Q)$, then 
we set 
\[
\Sigma_f := \Spec\,(A/\sigma(Q)). 
\]
By Proposition \ref{p-sigma-welldef}, this definition is independent of the choice of $Q$. 
\item 
Fix an affine open cover $T = \bigcup_{i \in I} T_i$ 
such that 
$f^{-1}(T_i)$ is $T_i$-isomorphic to a conic on $\mathbb P^2_{T_i}$ for every $i \in I$ 
(the existence of such an open cover is guaranteed by Proposition \ref{p-conic-emb}). 
We can write 
\[
X_i := f^{-1}(T_i) = \Proj\, A_i[x, y, z]/(Q_i), \qquad Q_i \in A_i[x, y, z], \qquad 
A_i := \Gamma(T_i, \MO_T). 
\]
By (1), we have a closed subscheme $\Sigma_{f_i}$ of $T_i$ for every $i \in I$. 
Again by (1), we obtain $\Sigma_{f_i}|_{T_i \cap T_j} = \Sigma_{f_j}|_{T_i \cap T_j}$, 
so that 
there exists a closed subscheme $\Sigma_f$ on $T$ such that $\Sigma_f|_{T_i} = \Sigma_{f_i}$ 
for every $i \in I$. 
It follows from (1) that $\Sigma_f$ does not depend on the choice of the affine open cover $T = \bigcup_{i \in I} T_i$. 
\end{enumerate}
\end{dfn}

\begin{rem}\label{r-Sigma-bc}
Let 
\[
\begin{CD}
X' @>\alpha >> X\\
@VVf'V @VVfV\\
T' @> \beta >> T
\end{CD}
\]
be a cartesian diagram of noetherian schemes, 
where $f$ is a conic bundle. 
In particular, also $f'$ is a conic bundle. 
In this case, the equation  
\[
\beta^{-1}(\Sigma_f) = \Sigma_{f'}
\]
holds by Definition \ref{d-Sigma}, where $\beta^{-1}(\Sigma_f)$ is the scheme-theoretic inverse image of $\Sigma_f$. 
\end{rem}

\subsection{Singular fibres}\label{ss-Delta-Sigma-loci}

The purpose of this subsection is to prove Theorem \ref{t-disc-non-sm}, 
which states that $\Delta_f$ (resp. $\Sigma_f$) 
is set-theoretically equal to the locus parametrising non-smooth fibres (resp. geometrically non-reduced fibres). 
Since the problem is reduced to the case when the base scheme is $\Spec\,\kappa$ for a field $\kappa$, 
we start with the following.

\begin{prop}\label{p-char0-field}\label{p-char2-field}
Let $\kappa$  be a field. 
Let 
\[
Q:= Q(x, y, z) := ax^2 + by^2 + cz^2 + \alpha yz + \beta zx + \gamma xy \in \kappa[x, y, z] 
\setminus \{0\} 
\]
and 
set 
\[
C := \Proj\,\kappa[x, y, z]/(Q), 
\]
which is a conic on $\P^2_{\kappa}$. Then the following hold. 
\begin{enumerate}
\item 
$C$ is not smooth over $\kappa$ if and only if $\delta(Q) =0$. 
\item 
$C$ is not geometrically reduced if and only if $\sigma(Q)=0$. 
\end{enumerate}
\end{prop}

\begin{proof}
{\cred By 
Proposition \ref{p-delta-welldef} 
and Proposition \ref{p-sigma-welldef}, 
we may perform coordinate changes.}

We first treat the case when ${\cred \kappa}$ is of characteristic $\neq 2$. 
By Proposition \ref{p-char0-local0}, 
we may assume that $\alpha = \beta = \gamma =0$: 
\[
Q = ax^2 + by^2+cz^2. 
\]
In this case, 
the following hold (Definition \ref{d-delta}, Definition \ref{d-sigma}): 
\[
\delta(Q) = 4abc, \qquad 
\sigma(Q) = (ab, bc, ca). 
\]
Thus (1) and (2) hold by Proposition \ref{p-char0-local0}. 

Hence we may assume that ${\cred \kappa}$ is an algebraically closed field of characteristic two. 
By Proposition \ref{p-char2-local2}(2), 
the problem is reduced to the case when $Q$ is equal to one of $x^2, x^2+yz, yz$. 
For each case, (1) and (2) hold by Definition \ref{d-delta} and Definition \ref{d-sigma}, respectively. 
\end{proof}

\begin{thm}\label{t-disc-non-sm}
Let $f : X \to T$ be a conic bundle and let $t$ be a point of $T$. 
Then the following hold. 
\begin{enumerate}
\item 
$X_t$ is not smooth 
if and only if $t \in \Delta_f$. 
\item 
$X_t$ is not geometrically reduced 
if and only if $t \in \Sigma_f$. 
\end{enumerate}
\end{thm}

\begin{proof}
By Remark \ref{r-disc-bc}, we may assume that $T=\Spec\,\kappa$, where $\kappa$ is a field. 
Then the assertion follows from Proposition \ref{p-char0-field}. 
\end{proof}

Theorem \ref{t-disc-non-sm} immediately deduces the set-theoretic inclusion 
$\Sigma_f \subset \Delta_f$. 
As the following proposition shows, this inclusion holds even as closed subschemes.

\begin{prop}\label{p-Sigma-Delta}
Let $f : X \to T$ be a conic bundle. 
Then the inclusion $\Sigma_f \subset \Delta_f$ of closed subschemes holds, i.e.,  
the inclusion $I_{\Sigma_f} \supset I_{\Delta_f}$ holds for the corresponding ideal sheaves.  
\end{prop}

\begin{proof}
By Proposition \ref{p-conic-emb}, we may assume that $T = \Spec\,A$ and 
\[
X = \Proj\,A[x, y, z]/(Q) \subset \P^2_A
\]
for 
\[ 
Q =ax^2 + by^2 + cz^2 + \alpha yz + \beta zx + \gamma xy \in A[x, y, z]. 
\]
It suffices to show that $\delta(Q) \in \sigma(Q)$. 
Recall that 
\begin{itemize}
\item $\delta(Q) = 4abc + \alpha \beta \gamma -a\alpha^2 -b\beta^2 - c\gamma^2$ and 
\item $\sigma(Q) 
= (4ab -\gamma^2, 4bc - \alpha^2, 4ca - \beta^2, 
2a \alpha - \beta \gamma, 2b\beta - \gamma \alpha, 2c\gamma - \alpha \beta)$. 
\end{itemize}
It holds that 
\[
4abc - a\alpha^2, \quad 4abc -b\beta^2, \quad 4abc -c\gamma^2 \in \sigma(Q), 
\]
and hence 
\[
12 abc -a\alpha^2 -b\beta^2 - c\gamma^2 \in \sigma(Q).  
\]
By $8abc -2a\alpha^2 = 2a(4bc-\alpha^2) \in \sigma(Q)$ and 
$2a\alpha^2 - \alpha \beta \gamma =\alpha( 2a\alpha -\beta\gamma) \in \sigma(Q)$, we get 
\[
8abc - \alpha \beta \gamma \in \sigma(Q), 
\]
which implies 
\[
\delta(Q) = (12 abc -a\alpha^2 -b\beta^2 - c\gamma^2) - (8abc - \alpha \beta \gamma) \in \sigma(Q).
\]
\end{proof}

\subsection{Case of characteristic two}\label{ss-Sigma-char2}

In this subsection, we introduce a simpler version $\Sigma'_f$ of $\Sigma_f$ which can be defined only in characteristic two. 
Let us start by recalling the defining equations of $\Delta_f$ and $\Sigma_f$. 

\begin{nasi}\label{n-delta-sigma-char2}
Let $A$ be  a noetherian $\F_2$-algebra. 
For 
\[
Q:= ax^2 + by^2 + cz^2 + \alpha yz + \beta za + \gamma xy \in A[x, y, z] 
\]
we have the following (Definition \ref{d-delta}, Definition \ref{d-sigma}).  
\begin{enumerate}
\item $\delta(Q) = \alpha \beta \gamma +a \alpha^2 + b\beta^2+ c \gamma^2$. 
\item $\sigma(Q) = (\alpha^2, \beta^2, \gamma^2, \alpha\beta, \beta\gamma, \gamma\alpha)$. 
\end{enumerate}
\end{nasi}


We now introduce the ring-theoretic counterpart $\sigma'$ to $\Sigma'_f$. 

\begin{dfn}\label{d-sigma'}
Let $A$ be a noetherian $\F_2$-algebra. 
For 
\[
Q := ax^2 + by^2 + cz^2 + \alpha yz + \beta za + \gamma xy \in A[x, y, z], 
\]
we define 
\[
\sigma'(Q) = (\alpha, \beta, \gamma), 
\] 
which is an ideal of $A$. 
\end{dfn}

\begin{prop}\label{p-sigma'-welldef}
Let $A$ be a noetherian $\F_2$-algebra.  
For  $M \in \GL_3(A)$ and 
\[
Q
:= ax^2 + by^2 +cz^2 + \alpha yz + \beta zx + \gamma xy \in A[x, y, z], 
\]
the equality 
\[
\sigma'(Q) = \sigma'(Q^M)
\]
of ideals hold. 
\end{prop}

\begin{proof}
By the same argument as in the proof of Proposition \ref{p-delta-welldef}, 
we may assume that $A$ is a noetherian local ring and one of the following holds.  
\begin{enumerate}
\item 
$M = 
\begin{pmatrix}
\lambda & 0 & 0\\
0 & 1 & 0\\
0 & 0 & 1\\
\end{pmatrix}$ for some $\lambda \in A^{\times}$. 
\item 
$M = 
\begin{pmatrix}
0 & 1 & 0\\
1 & 0 & 0\\
0 & 0 & 1\\
\end{pmatrix}$. 
\item 
$M = 
\begin{pmatrix}
1 & 0 & \mu\\
0 & 1 & 0\\
0 & 0 & 1\\
\end{pmatrix}$ for some $\mu \in A$. 
\end{enumerate}We only treat the case (3), i.e., we assume that $M = 
\begin{pmatrix}
1 & 0 & \mu\\
0 & 1 & 0\\
0 & 0 & 1\\
\end{pmatrix}$. 
Then there exist $
c' \in A$ such that the following holds: 
\begin{eqnarray*}
Q^M(x, y, z) 
&=& Q(x+\mu z, y, z) \\
&=& a(x+\mu z)^2 + by^2 +cz^2 + \alpha yz + \beta z(x+\mu z) + \gamma (x+\mu z)y\\
&=& ax^2 +by^2 +c'z^2 + (\alpha + \mu \gamma)yz + \beta zx + \gamma xy. 
\end{eqnarray*}
Hence we obtain 
$\sigma'(Q^M) = (\alpha + \mu \gamma, \beta, \gamma) = (\alpha, \beta, \gamma) = \sigma(Q)$. 
\end{proof}

\begin{dfn}\label{d-Sigma'}
Let $f: X \to T$ be a conic bundle. 
We define the closed subscheme $\Sigma'_f$ of $T$ as follows. 
\begin{enumerate}
\item If $T=\Spec\,A$ and $X =\Proj\,A[x, y, z]/(Q)$, then 
we set 
\[
\Sigma'_f := \Spec\,(A/\sigma'(Q)). 
\]
By Proposition \ref{p-sigma'-welldef}, this definition is independent of the choice of $Q$. 
\item 
Fix an affine open cover $T = \bigcup_{i \in I} T_i$ 
such that 
$f^{-1}(T_i)$ is $T_i$-isomorphic to a conic on $\mathbb P^2_{T_i}$ for every $i \in I$ 
(the existence of such an open cover is guaranteed by Proposition \ref{p-conic-emb}). 
We can write 
\[
X_i := f^{-1}(T_i) = \Proj\, A_i[x, y, z]/(Q_i), \qquad Q_i \in A_i[x, y, z], \qquad 
A_i := \Gamma(T_i, \MO_T). 
\]
By (1), we have a closed subscheme $\Sigma'_{f_i}$ of $T_i$ for every $i \in I$. 
Again by (1), we obtain $\Sigma'_{f_i}|_{T_i \cap T_j} = \Sigma'_{f_j}|_{T_i \cap T_j}$, 
so that 
there exists a closed subscheme $\Sigma'_f$ on $T$ such that $\Sigma'_f|_{T_i} = \Sigma'_{f_i}$ 
for every $i \in I$. 
It follows from (1) that $\Sigma'_f$ does not depend on the choice of the affine open cover $T = \bigcup_{i \in I} T_i$. 
\end{enumerate}
\end{dfn}

\begin{rem}\label{r-Sigma'-bc}
Let 
\[
\begin{CD}
X' @>\alpha >> X\\
@VVf'V @VVfV\\
T' @> \beta >> T
\end{CD}
\]
be a cartesian diagram of schemes, 
where $T$ and $T'$ are noetherian schemes and $f$ is a conic bundle. 
In particular, also $f'$ is a conic bundle. 
In this case, the equation  
\[
\beta^{-1}(\Sigma'_f) = \Sigma'_{f'}
\]
holds by Definition \ref{d-Sigma'}, 
where $\beta^{-1}(\Sigma'_f)$ is the scheme-theoretic inverse image of $\Sigma'_f$. 
\end{rem}

\begin{rem}\label{r-Sigma'-non-red}
Given a conic bundle $f: X \to T$, 
{\cred we have an inclusion $\Sigma_f \supset \Sigma'_f$ of closed subschemes 
of $T$ and} the following holds ((\ref{n-delta-sigma-char2}), Definition \ref{d-sigma'}):
\[
(\Sigma_f)_{\red} = (\Sigma'_f)_{\red}.
\]
\end{rem}

\section{Singularities of ambient spaces}\label{s-sing-ambient}

Given a conic bundle $f: X \to T$, 
it is natural to seek a relation between 
$\Delta_f$ and the singularities of $X$. 
The main result of this section is to prove that the following are equivalent under the assumption that $T$ is regular (Theorem \ref{t-sm-red-fib}).
\begin{enumerate}
    \item $\Delta_f$ is regular. 
    \item $X$ is regular and every fibre of $f$ is geometrically reduced. 
\end{enumerate}
In particular, if $f: X \to T$ is a conic bundle between regular noetherian schemes, then the equality $(\Delta_f)_{{\rm non}\text{-}{\rm reg}} = \Sigma_f$  holds, where 
$(\Delta_f)_{{\rm non}\text{-}{\rm reg}}$ denotes the non-regular locus of $\Delta_f$ (Theorem \ref{t-Delta-nonreg}). 
{\cred These results are known when $f : X \to T$ is a conic bundle of smooth varieties over an algebraically closed field of characteristic $\neq 2$ \cite[Proposition 1.8]{Sar82}.}

\begin{prop}\label{p-Delta-reg-red}
Let $f: X \to T$ be a generically smooth conic bundle, where $T$ is a regular noetherian integral scheme. 
Fix a point $0 \in \Delta_f$. 
If $\Delta_f$ is regular at $0$, then the fibre $X_0$ is geometrically reduced. 
\end{prop}

\begin{proof}
Assume that $X_0$ is not geometrically reduced. 
It suffices to show that $\Delta_f$ is not regular at $0$. 
By taking the strict henselisation of the local ring $\MO_{T, 0}$, 
we may assume that $T= \Spec\,A$, 
$(A, \m, \kappa)$ is a strictly henselian regular local ring, 
and $0$ is the closed point of $T$. 
We can write 
\[
X = \Proj\,A[x, y, z]/(Q) \qquad\text{for}\qquad 
Q = ax^2+by^2 +cz^2 + \alpha yz + \beta zx + \gamma xy \in A[x, y, z]. 
\]
By Definition \ref{d-delta} and Definition \ref{d-disc}, we get 
\[
\Delta_f = \Spec\,(A/\delta(Q)A)\qquad\text{and}\qquad
\delta(Q) = 4abc + \alpha \beta \gamma - a \alpha^2 - b\beta^2 -c \gamma^2 \in A.  
\]

If $\kappa$ is of characteristic two, then 
Theorem \ref{t-disc-non-sm} and 
Remark \ref{r-Sigma'-non-red} imply 
\[
2, \alpha, \beta, \gamma \in \m. 
\]
In particular, we get $\delta(Q) \in \m^2$, and hence $\Delta_f$ is not regular at $0$. 

Then the problem is reduced to the case when 
$\kappa$ is of characteristic $\neq 2$. 
By Corollary \ref{c-char0-local2}(1), we may assume that 
\[
X = \Proj\,A[x, y, z]/(x^2 +by^2 + cz^2+\alpha yz), 
\]
where $b, c, \alpha \in \m$. 
By $\delta(Q) = 4bc -\alpha^2 \in \m^2$, $\Delta_f$ is not regular at $0$. 
\end{proof}

\begin{prop}\label{p-sm-red-fib}
Let $f: X \to T$ be a  generically smooth conic bundle, where $X$ and $T$ are regular noetherian integral schemes. 
Fix a point $0 \in \Delta_f$. 
Then the following are equivalent. 
\begin{enumerate}
\item $X_0$ is geometrically reduced. 
\item $\Delta_f$ is regular at $0$. 
\end{enumerate}
\end{prop}

\begin{proof}
The implication (2) $\Rightarrow$ (1) has been settled in Proposition \ref{p-Delta-reg-red}. 
Let us show the opposite one: (1) $\Rightarrow$ (2). 
By taking the strict henselisation of the local ring $\MO_{T, 0}$, 
we may assume that $T= \Spec\,A$, 
$(A, \m, \kappa)$ is a strictly henselian regular local ring, 
and $0$ is the closed point of $T$. 


\setcounter{step}{0}

\begin{step}\label{s1-sm-red-fib}
The implication (1) $\Rightarrow$ (2) holds if $\kappa$ is of characteristic $\neq 2$. 
\end{step}

\begin{proof}[Proof of Step \ref{s1-sm-red-fib}]
Assume (1). 
By Corollary \ref{c-char0-local2}(2), we may assume that 
\[
X = \Proj\,A[x, y, z]/(x^2 + y^2+cz^2)
\]
for some $c \in A$. 
Then the following hold ({\cred Definition \ref{d-delta}}, Theorem \ref{t-disc-non-sm}): 
\[
\delta(Q) = 4c \qquad {\rm and}\qquad \sigma(Q) = A. 
\]
It suffices to show that $c \not\in \m^2$. 
Suppose $c \in \m^2$. 
Let $\n$ be the closed point of $X$ lying over $0 \in T$ 
corresponding to $[0:0:1] \in \Proj\,\kappa[x, y, z] =\P^2_{\kappa}$. 
{\cred For $\widetilde{x} := x/z$ and $\widetilde{y}:=y/z$}, 
we have the induced injective ring homomorphism 
\[
A \hookrightarrow  B:= A[\widetilde{x}, \widetilde{y}]/(\widetilde x^2+\widetilde y^2+c), 
\]
where $\Spec\,B$ is an open neighbourhood of $\n \in X$. 
Via this injection, we consider $A$ as a subring of $B$. 
We have $\m \subset \mathfrak n$, which implies $c \in \m^2 \subset \mathfrak n^2$. 
By $\widetilde x \in \n$ and $\widetilde y \in \mathfrak n$,  
we obtain $\widetilde x^2+\widetilde y^2+c \in \mathfrak n^2$, and hence $B$ is not regular at $\mathfrak n$, which is a contradiction. 
Thus (1) implies (2). 
This completes the proof of Step \ref{s1-sm-red-fib}. 
\end{proof}

\begin{step}\label{s2-sm-red-fib}
The implication (1) $\Rightarrow$ (2) holds if  $\kappa$ is of characteristic two. 
\end{step}

\begin{proof}[Proof of Step \ref{s2-sm-red-fib}]
Assume (1), i.e., $X_0$ is geometrically reduced. 
It follows from $0 \in \Delta_f$ that $X_0$ is not smooth. 
By Proposition \ref{p-char2-local}(2), we can write 
\[
X = \Proj\,A[x, y, z]/(ax^2 + yz) 
\]
for some $a \in \m$. 
Then we have 
\[
\delta(Q) = -a \qquad \text{and}\qquad 
\Delta_f = \Spec\,(A/\delta(Q)A). 
\]
By 
\[
D_+(x) = \Spec\,A[y, z]/(a+yz) \subset X, 
\]
we have $a \not\in \m^2$, as otherwise $D_+(x)$ would not be regular. 
Hence $\Delta_f = \{ a=0\}$ is regular at $0$. 
Thus (2) holds. 
This completes the proof of Step \ref{s2-sm-red-fib}. 
\end{proof}
Step \ref{s1-sm-red-fib} and 
Step \ref{s2-sm-red-fib} complete the proof of Proposition \ref{p-sm-red-fib}. 
\end{proof}

\begin{thm}\label{t-Delta-nonreg}
Let $f: X \to T$ be a generically smooth conic bundle, where $X$ and $T$ are regular noetherian integral schemes. 
Set 
\[
(\Delta_f)_{{\rm non}\text{-}{\rm reg}} := \{ t \in \Delta_f\,|\, \MO_{\Delta_f, t}
\text{ is not a regular local ring}\}. 
\]
Then the set-theoretic equality  
\[
(\Delta_f)_{{\rm non}\text{-}{\rm reg}} = \Sigma_f
\]
holds. 
\end{thm}

\begin{proof}
By Theorem \ref{t-disc-non-sm} or Proposition \ref{p-Sigma-Delta}, 
both sides are subsets of $\Delta_f$ . 

Let us show the inclusion $(\Delta_f)_{{\rm non}\text{-}{\rm reg}} \supset \Sigma_f$. 
Pick a point $t \in \Delta_f \setminus (\Delta_f)_{{\rm non}\text{-}{\rm reg}}$, i.e., 
$t$ is a point of $\Delta_f$ at which $\Delta_f$ is regular. 
By Proposition \ref{p-Delta-reg-red}, the fibre $X_t$ over $t$ is geometrically reduced. 
It follows from Theorem \ref{t-disc-non-sm} that $t \not\in \Sigma_f$. 
This completes the proof of the inclusion $(\Delta_f)_{{\rm non}\text{-}{\rm reg}} \supset \Sigma_f$. 

Let us prove the opposite inclusion $(\Delta_f)_{{\rm non}\text{-}{\rm reg}} \subset \Sigma_f$. 
Pick a point $t \in \Delta_f  \setminus \Sigma_f$. 
Then $X_t$ is not smooth but geometrically reduced (Theorem \ref{t-disc-non-sm}). 
By Proposition \ref{p-sm-red-fib}, $\Delta_f$ is regular at $t$, i.e., $t \not\in (\Delta_f)_{{\rm non}\text{-}{\rm reg}}$. 
Therefore, we get $(\Delta_f)_{{\rm non}\text{-}{\rm reg}} \subset \Sigma_f$. 
\end{proof}

\begin{thm}\label{t-sm-red-fib}
Let $f: X \to T$ be a generically smooth conic bundle, where $T$ is a regular noetherian integral scheme. 
Then the following are equivalent. 
\begin{enumerate}
\item $\Delta_f$ is regular. 
\item $X$ is regular and any fibre of $f$ is geometrically reduced. 
\end{enumerate}
\end{thm}

\begin{proof}
By Proposition \ref{p-sm-red-fib}, (2) implies (1). 
Assume (1), i.e., $\Delta_f$ is regular. 
By Proposition \ref{p-Delta-reg-red}, any fibre of $f$ is geometrically reduced. 
Fix a point ${\cred P} \in X$ and set $0:=f({\cred P})$. 
Hence it is enough to show that $X$ is regular at ${\cred P}$. 
Taking the strict henselisation of the local ring $\MO_{T, 0}$, 
we may assume that $T= \Spec\,A$, $(A, \m, \kappa)$ is a strictly henselian regular local ring, 
and $0$ is the closed point of $T$. 

Since $X_0$ is geometrically reduced, we can write 
\[
X =\Proj\, A[x, y, z]/(ax^2 + yz) 
\]
for some $a \in A$ (Corollary \ref{c-char0-local2}(2), Proposition \ref{p-char2-local}(2)). 
Then we obtain  $\delta(Q) = -a$ for $Q:=ax^2+yz$ (Definition \ref{d-delta}). 
As $\Delta_f = \Spec\,(A/\delta(Q)A)$ is regular, we have $a \not\in \m^2$. 
It suffices to show that $A[y, z]/(a + yz)$ is regular. 
We may assume that $a \in \m$. 
Since $(A/\m)[y, z]/(yz) = \kappa[y, z]/(yz)$ is smooth outside the origin $(0, 0)$, 
we may assume that ${\cred P} =(0, 0)$. 
Then the prime ideal ${\cred P}$ is the image of the maximal ideal $\n := \m A[y, z] +(y, z)A[y, z]$ of $A[y, z]$ 
by the natural surjection $A[y, z] \to A[y, z]/(a + yz)$. 
By $a \not\in \m^2$, we obtain $a + yz \not\in \n^2$, and hence $A[y, z]/(a+yz)$ is regular 
at ${\cred P}$. 
Thus $X$ is regular. 
\end{proof}

\begin{prop}\label{p-geom-red-cano}
Let $k$ be an algebraically closed field. 
Let $f: X \to T$ be a conic bundle, where 
$T$ is a smooth variety over  $k$. 
Assume that 
\begin{enumerate}
\item $f$ is generically smooth, and 
\item any fibre of $f$ is geometrically reduced. 
\end{enumerate}
Then $X$ is canonical in the sense of \cite[Definition 2.8]{Kol13}. 
If $k$ is characteristic $p>0$, then $X$ is strongly $F$-regular. 
\end{prop}

In the proof below, we use the following fact: if $X' \to X$ is an \'etale surjective morphism, then $X$ is canonical if and only if so is $X'$ 
 \cite[Proposition 2.15]{Kol13}. 

\begin{proof}
Let us prove the assertion under the assumption that $k$ is  of characteristic $p>0$. 
We prove the assertion by induction on $\dim T$. 

We first treat the case when $\dim T=1$. 
By 
{\cred Corollary \ref{c-char0-local2}(2)} 
and Proposition \ref{p-char2-local}(2), 
we may assume, after taking suitable \'etale cover of $X$, that 
\[
X = \Proj\,A[x, y, z]/(xy+cz^2), \qquad T = \Spec\,A, \qquad c \in A. 
\]
The singular locus of $X$ is contained in 
\[
D_+(z) = \Spec\,A[x, y]/(xy + c). 
\]
This is an $A_n$-singularity, and hence canonical and strongly $F$-regular 
(note that these singularities are toric).  

Assume that $\dim T \geq 2$ and 
that the assertion holds for the lower dimensional case. 
By (1), $\Delta_f$ is a nonzero effective Cartier divisor on $T$. 
Fix a closed point $0 \in T$ around which we will work. 
Take a general smooth prime divisor $T'$ on $T$ passing through $0$. 
Set $X' := X \times_T T'$. 
Then $f':X' \to T'$ is a conic bundle over $T'$. 
In particular, $X'$ is strongly $F$-regular by the induction hypothesis. 
Since $\dim(X' \cap {\rm Sing}\,X) < \dim X'-1$ (where $\dim \emptyset := -\infty$), 
we have 
\[
(K_X+X')|_{X'} =K_{X'}, 
\]
i.e., the  different $\Diff_{X'}(0)$ is equal to zero (for the definition of $\Diff_{X'}(0)$, see \cite[Subsection 4.1]{Kol13}). 
Since $(X', \Diff_{X'}(0) =0)$ is strongly $F$-regular, 
it follows from inversion of adjunction \cite[Theorem 4.1]{Das15} that $(X, X')$ is purely $F$-regular around $X'$. 
In particular, $X$ is strongly $F$-regular and $(X, X')$ is log canonical 
\cite[Theorem 3.3]{HW02}. 
Since $X'$ is a nonzero effective Cartier divisor, $X$ is canonical. 
This completes the proof for the case when $k$ is  of characteristic $p>0$. 

If $k$ is of characteristic zero, then the above argument works 
by using the inversion of adjunction for log canonicity \cite[Theorem]{Kaw07}.
\end{proof}

For a later use, 
we establish the following Jacobian criterion. 

\begin{prop}[Jacobian criterion]\label{p-Jac-criterion}
Let $k$ be a field. 
Set $A:=k[t_1, ..., t_d]$ and let 
\[
X := \Proj\,A[x, y, z]/(Q)
\]
be a conic bundle over $A$ for $Q := Q(t_1, ..., t_d, x, y, z) \in A[x, y, z]$. 
Fix a $k$-rational point 
\[
w := (t_{10}, ..., t_{d0}) \times [x_0:y_0:z_0] \in X \subset \mathbb P^2_A, 
\]
where $t_{10}, ..., t_{d0}, x_0, y_0, z_0 \in k$. 
Then the following are equivalent. 
\begin{enumerate}
\item $X$ is not smooth at $w$. 
\item The equation 
\[
\partial_{x}Q(w) =\partial_{y}Q(w) =\partial_{z}Q(w)=\partial_{t_1}Q(w) = \cdots =\partial_{t_d}Q(w) =0,
\]
holds, where $\partial_{\bullet}$ denotes the partial differential with respect to $\bullet$. 
\end{enumerate}
\end{prop}

\begin{proof}
Taking the base change to the algebraic closure, 
the problem is reduced to the case when $k$ is algebraically closed. 
If $k$ is of characteristic $\neq 2$, then the assertion is well known. 
We may assume that $k$ is of characteristic two. 

We have 
\[
Q = ax^2 + by^2 + cz^2 + \alpha yz + \beta zx + \gamma xy, \qquad a, b, c, \alpha, \beta, \gamma \in A=k[t_1, ..., t_d]. 
\]
By symmetry, we may assume that $x_0 =1$. 
Set 
\[
Q' := Q'(t_1, ..., t_d, y, z) := Q(t_1, ..., t_d, 1, y, z) = 
 a + by^2 + cz^2 + \alpha yz + \beta z + \gamma y, 
\]
which corresponds to the affine open subset $X' :=D_+(x)$ of $X$ defined by $x \neq 0$. 
Consider the following conditions. 
\begin{enumerate}
\item[(1)'] $X'$ is not smooth at $w$. 
\item[(2)'] 
\[
\partial_{y}Q'(w) =\partial_{z}Q'(w)=\partial_{t_1}Q'(w) = \cdots =\partial_{t_d}Q'(w) =0. 
\]
\end{enumerate}
It is clear that (1) $\Leftrightarrow$ (1)' $\Leftrightarrow$ (2)'. 
Therefore it suffices to show that (2) $\Leftrightarrow$ (2)'. 
It is clear that (2) $\Rightarrow$ (2)', because 
the partial differentials 
$\partial_{y}, \partial_{z}, \partial_{t_1}, ..., \partial_{t_d}$ 
commute with the substitution $x=1$. 

Therefore, it is enough to show that (2)' $\Rightarrow$ (2). 
Assume (2)'. 
Then we have 
\[
\partial_{y}Q(w) =\partial_{z}Q(w)=\partial_{t_1}Q(w) = \cdots =\partial_{t_d}Q(w) =0. 
\]
It suffices to show $\partial_{x}Q(w)=0$. 
By $\partial_{y}Q(w) =\partial_{z}Q(w)=0$, 
we have 
\[
\alpha z_0 + \gamma = 0, \qquad \alpha y_0 +\beta =0. 
\]
Therefore, we get  
\[
\partial_{x}Q(w)=\beta z_0 + \gamma y_0 =(\alpha y_0 )z_0 + (\alpha z_0)y_0=0, 
\]
as required. 
\end{proof}

\begin{rem}\label{r-Jac-criterion}
Proposition \ref{p-Jac-criterion} holds even when $A=k[[t_1, ..., t_d]]$. 
In this case, the $k$-rational point $w$ is lying over $\m =(t_1, ..., t_d)$, 
and hence $t_{10} = \cdots = t_{d0}=0$. 
\end{rem}

\section{The Mori--Mukai formula}\label{s-MM-formula}


In this section, we prove the Mori--Mukai formula (Theorem \ref{t-cbf} in Subsection \ref{ss-MM-proof}): 
\begin{equation}\label{e0-s-MM-formula}
\Delta_f \equiv -f_*(K_{X/T}^2)
\end{equation}
for an arbitrary generically smooth conic bundle $f: X \to T$ between smooth projective varieties. 
As explained in Subsection \ref{ss-intro-overview}, 
the problem is reduced to the case when $f$ coincides with  $\Univ_{\P_T(E)/T}^{\theta} \to \Hilb_{\P_T(E)/T}^{\theta}$, 
where $T$ is a smooth projective curve, $E$ is a locally free sheaf on {\cred $T$ of} rank $3$, 
and $\Univ_{\P_T(E)/T}^{\theta}$ denotes the universal family associated with 
the relative Hilbert scheme $\Hilb_{\P_T(E)/T}^{\theta}$ parametrising conics. 
We start by studying the absolute Hilbert scheme $\Univ_{\P^2_{\Z}/\Z}^{\theta} \to \Hilb_{\P^2_{\Z}/{\Z}}^{\theta}$ parametrising conics (Subsection \ref{ss-univ-Hilb-abs}), 
because the absolute case determines the local structure of the relative one  $\Univ_{\P_T(E)/T}^{\theta} \to \Hilb_{\P_T(E)/T}^{\theta}$. 
More specifically, 
if $T = \bigcup_{i \in I} T_i$ is an open cover trivialising $E$, 
then $\Univ_{\P_{T_i}(E|_{T_i})/T_i}^{\theta} \to \Hilb_{\P_{T_i}(E|_{T_i})/T_i}^{\theta}$ coincides with 
$\Univ_{\P^2_{\Z}/\Z}^{\theta} \times_{\Z} T_i \to \Hilb_{\P^2_{\Z}/{\Z}}^{\theta} \times_{\Z} T_i$. 
We then summarise properties on the relative version 
$\Univ_{\P_T(E)/T}^{\theta} \to \Hilb_{\P_T(E)/T}^{\theta}$ in Subsection \ref{ss-univ-Hilb-rel}, most of which are immediate conclusions from the absolute case. 

In order to include generically non-smooth conic bundles, we shall introduce an invertible sheaf $\Delta_f^{\bdl}$, which we shall call the discriminant bundle (Subsection \ref{s-disc-bdl}), 
for an arbitrary conic bundle $f: X \to S$ with $S$ integral. 
When $f$ is generically smooth, it satisfies $\Delta_f^{\bdl} \simeq \MO_T(\Delta_f)$ 
(Remark \ref{r-disc-bdl-div}) and 
the above Mori--Mukai formula (\ref{e0-s-MM-formula}) is generalised as follows: 
\begin{equation}\label{e2-s-MM-formula}
\Delta_f^{\bdl} \equiv -f_*(K_{X/T}^2). 
\end{equation}
Moreover, the proof of (\ref{e0-s-MM-formula}) becomes simpler by using the notion of  discriminant bundles, as otherwise we would  need to be careful with the 
{\cred generically} smooth condition. 

\subsection{The universal family}\label{ss-univ-Hilb-abs}

In this subsection, we study the conic bundle $\Univ_{\P^2_{\Z}/\Z}^{\theta} \to \Hilb_{\P^2_{\Z}/{\Z}}^{\theta}$, 
where $\Hilb_{\P^2_{\Z}/{\Z}}^{\theta}$ denotes the Hilbert scheme of conics and 
$\Univ_{\P^2_{\Z}/\Z}^{\theta}$ is its universal family. 
As explained above, 
this plays an essential role in the proof of the Mori--Mukai formula (\ref{e0-s-MM-formula}). 
It is also worth studying this conic bundle $\Univ_{\P^2_{\Z}/\Z}^{\theta} \to \Hilb_{\P^2_{\Z}/{\Z}}^{\theta}$ 
as a concrete example. 
For foundations on Hilbert schemes, we refer to \cite{FGI05}. 

\begin{nota}\label{n-Hilb-Univ}
\begin{enumerate}
\item Set $\theta(m) :=2m+1 \in \Z[m]$, 
which is nothing but the Hilbert polynomial of conics. 
In other words, if $\kappa$ is a field and $C$ is a conic on $\P^2_{\kappa}$, 
then the following holds for any $m \in \Z$: 
\[
\theta(m) =\chi(C, \MO_{\P^2}(m)|_C) = 2m+1. 
\]
\item We have 
\[
\Hilb^{\theta}_{\P^2_{\Z}/\Z} \simeq \P^5_{\Z} = \Proj\,\Z[a, b, c, \alpha, \beta, \gamma],  
\]
where $\Z[a, b, c, \alpha, \beta, \gamma]$ denotes the polynomial ring over $\Z$ with six variables $a, b, c, \alpha, \beta, \gamma$. 
We identify $\Hilb^{\theta}_{\P^2_{\Z}/\Z}$ with 
$\P^5_{\Z} = \Proj\,\Z[a, b, c, \alpha, \beta, \gamma]$. 
\item Let $\Univ^{\theta}_{\P^2_{\Z}/\Z} \subset \Hilb^{\theta}_{\P^2_{\Z}/\Z} \times_{\Z} \P^2_{\Z}$ be the universal closed subscheme: 
\[
\Univ^{\theta}_{\P^2_{\Z}/\Z} = \{ ax^2 + by^2 + cz^2 + \alpha yz + \beta zx + \gamma xy =0\} \subset  
\P^2_{\Z}  \times_{\Z} \Hilb^{\theta}_{\P^2_{\Z}/\Z} 
\]
\item 
The induced composite morphism 
\[
f_{\univ} :   \Univ^{\theta}_{\P^2_{\Z}/\Z} \hookrightarrow \P^2_{\Z}  \times_{\Z} \Hilb^{\theta}_{\P^2_{\Z}/\Z} \xrightarrow{ {\rm pr}_2} \Hilb^{\theta}_{\P^2_{\Z}/\Z}
\]
is a generically smooth conic bundle. 
Set $\Delta_{\univ} := \Delta_{f_{\univ}} \subset \Hilb^{\theta}_{\P^2_{\Z}/\Z}$, 
which satisfies the following (Definition \ref{d-delta} and Definition \ref{d-disc}): 
\[
\Delta_{\univ}  =\{ 4abc + \alpha \beta \gamma - a\alpha^2 -b \beta^2 -c \gamma^2 =0\} 
\subset \Hilb^{\theta}_{\P^2_{\Z}/\Z} = \P^5_{\Z}= \Proj\,\Z[a, b, c, \alpha, \beta, \gamma]. 
\]
Similarly, we set $\Sigma_{\univ} := \Sigma_{f_{\univ}} \subset  \Hilb^{\theta}_{\P^2_{\Z}/\Z}$. 
\item 
For a field $K$ and the morphism $\Spec\,K \to \Spec\,\Z$, we set 
\[
\Delta_{\univ, K} := \Delta_{\univ} \times_{\Z} K \quad {\rm and} \quad 
\Sigma_{\univ, K} := \Sigma_{\univ} \times_{\Z} K. 
\]
\end{enumerate}
\end{nota}

\begin{prop}\label{p-univ-sm}
We use Notation \ref{n-Hilb-Univ}. Then $\Univ^{\theta}_{\P^2_{\Z}/\Z}$ is smooth over $\Z$. 
In particular, 
$\Univ^{\theta}_{\P^2_{\Z}/\Z}, \Univ^{\theta}_{\P^2_{\Q}/\Q},$ and $\Univ^{\theta}_{\P^2_{\F_p}/\F_p}$ 
are regular for  any prime number $p$. 
\end{prop}

\begin{proof}
Both $f_{\univ} : \Univ^{\theta}_{\P^2_{\Z}/\Z} \to \Hilb^{\theta}_{\P^2_{\Z}/\Z}$ and 
$\Hilb^{\theta}_{\P^2_{\Z}/\Z} \to \Spec\,\Z$ are flat by $\Hilb^{\theta}_{\P^2_{\Z}/\Z} \simeq \P^5_{\Z}$. 
Hence also $\Univ^{\theta}_{\P^2_{\Z}/\Z} \to \Spec\,\Z$ is flat. 
Fix a prime number $p$. 
It suffices to show that 
\[
\Univ^{\theta}_{\P^2_{\overline{\F}_p}/\overline{\F}_p} = \{ 
Q :=ax^2 + by^2+cz^2 + \alpha yz + \beta zx + \gamma xy=0 \} \subset \P^2_{\overline{\F}_p} \times_{\overline{\F}_p} \P^5_{\overline{\F}_p}
\]
is smooth, where $\overline{\F}_p$ denotes the algebraic closure of $\F_p$. 
We shall apply the Jacobian criterion. 
Suppose that 
$\Univ^{\theta}_{\P^2_{\overline{\F}_p}/\overline{\F}_p}$ 
is not smooth at a closed point 
\[
w := [x_0:y_0:z_0] \times [a_0:b_0:c_0:\alpha_0:\beta_0:\gamma_0] \in \Univ^{\theta}_{\P^2_{\overline{\F}_p}/\overline{\F}_p} \subset 
\P^2_{\overline{\F}_p} \times_{\overline{\F}_p} \P^5_{\overline{\F}_p}
\]
Then $[a_0:b_0:c_0:\alpha_0:\beta_0:\gamma_0] \in \P^5_{\overline{\F}_p}$ 
is its image. 
By symmetry, we may assume that $a_0 \neq 0$ or $\alpha_0 \neq 0$.

Assume $\alpha_0\neq 0$. 
By the Jacobian criterion (Proposition \ref{p-Jac-criterion}), we have $\partial_aQ(w)=\partial_bQ(w) =\partial_cQ(w)=0$, which concludes that 
$x_0 = y_0=z_0=0$. This contradicts $[x_0 : y_0 : z_0] \in \P^2_{\overline{\F}_p}$. 

Assume $a_0 \neq 0$. 
For $\mathbb A^5_{\overline{\F}_p} = D_+(a)$ and $A :=\overline{\F}_p[b, c, \alpha, \beta, \gamma]$, 
we have 
\[
f_{\univ}^{-1}(D_+(a)) \simeq \{ x^2 + by^2+cz^2 + \alpha yz + \beta zx + \gamma xy=0 \} 
\subset \P^2_A. 
\]
By Proposition \ref{p-Jac-criterion}, 
our singular point satisfies $\partial_b = \partial_c=0$, which implies $[x_0:y_0:z_0] = [1:0:0]$. 
However, such a point does not lie on $\Univ^{\theta}_{\P^2_{\overline{\F}_p}/\overline{\F}_p}$, 
because $[x_0:y_0:z_0] = [1:0:0]$ and $a_0 \neq 0$ imply
\[
a_0x_0^2 + b_0y_0^2+c_0z_0^2 + \alpha_0 y_0z_0 + \beta_0 z_0x_0 + \gamma_0 x_0y_0
= a_0x_0^2 \neq 0. 
\]
\end{proof}

\begin{thm}\label{t-univ-Delta-Sigma}
We use Notation \ref{n-Hilb-Univ}. 
Then the following set-theoretic equality holds: 
\[
\Sigma_{\univ} = 
\{ v \in \Delta_{\univ} \,|\, \Delta_{\univ} \to \Spec\,\Z \text{ is not smooth at }v\}. 
\]
\end{thm}

\begin{proof}
Set 
\[
\overline{\Sigma}_{\univ} := \{ v \in \Delta_{\univ} \,|\, \Delta_{\univ} \to \Spec\,\Z \text{ is not smooth at }v\}, 
\]
which we equip with the reduced scheme structure. 
Fix a prime number $p$. 
It is enough to show the set-theoretic equality $\Sigma_{\univ} \times_{\Z} \F_p = \overline{\Sigma}_{\univ} \times_{\Z} \F_p$. 
Since $\Univ^{\theta}_{\P^2_{\F_p}/\F_p} \to \Hilb^{\theta}_{\P^2_{\F_p}/\F_p} $ is a conic bundle 
of regular schemes (Proposition \ref{p-univ-sm}), 
the assertion follows from Theorem \ref{t-disc-non-sm} and Proposition \ref{p-sm-red-fib}. 
\end{proof}

\begin{prop}\label{p-Delta-normal}
We use Notation \ref{n-Hilb-Univ}. 
Let $K$ be a field. 
Then $\Delta_{\univ, K}$ is geometrically integral over $K$ and geometrically normal over $K$. 
In particular, $\Delta_{\univ, K}$ is a normal prime divisor on $\P^5_K$
\end{prop}

\begin{proof}
We may assume that $K$ is an algebraically closed field by taking the base change 
{\cred (Remark \ref{r-disc-bc})}. 
Set $p:= {\rm char}\,K$, which is possibly zero. 
Recall that $\Delta_{\univ, K}$ is an effective Cartier divisor on 
$\P^5_K$. 
It suffices to show that $\Delta_{\univ, K}$ is a normal prime divisor.

\setcounter{step}{0}


\begin{step}\label{s1-Delta-normal}
$\Delta_{\univ, K}$ is smooth outside some closed subset of dimension two. 
\end{step}

\begin{proof}[Proof of Step \ref{s1-Delta-normal}]
We first treat the case when $p \neq 2$ and $p \neq 3$. 
Consider the hyperplane $\P^4_{K}$ defined by $\gamma=0$: 
\[
\P^4_{K} = \{ \gamma=0 \} \subset \P^5_{K}. 
\]
We then have 
\[
\Delta_{K} \cap \P^4_{K}  =\{\delta' := 4abc - a\alpha^2 -b \beta^2 =0\} 
\subset \P^4_{K} =\Proj\,K[a, b, c, \alpha, \beta]. 
\]
Then the singular locus 
$(\Delta_{\univ, K} \cap \P^4_{K})_{{\rm sing}}$ 
of $\Delta_{\univ, K} \cap \P^4_{K}$ 
is given by 
\[
\partial_a \delta' = 4bc -\alpha^2=0, \qquad \partial_b\delta' =4ac -\beta^2=0, \qquad \partial_c\delta' =4ab=0,
\]
\[
\partial_{\alpha} \delta' = -2a\alpha=0, \qquad \partial_{\beta}\delta' =-2b\beta=0. 
\]
On the open subset $\{ a \neq 0\}$, these equations become $b=\alpha=4ac -\beta^2 =0$, which is one-dimensional. 
By symmetry, these equations define a one-dimensional locus also on the open subset $\{ b \neq 0\}$. 
Finally, on the remaining closed subset $\{ a=b=0\}$, the equations become $\{ a=b=\alpha =\beta =0\}$, which is zero-dimensional. 
To summarise, $\Delta_{\univ, K} \cap \P^4_{K}$ is smooth outside a one-dimensional closed subset. 
Therefore, $\Delta_{\univ, K}$ is smooth outside a two-dimensional closed subset when $p \neq 2$ and $p \neq 3$. 
If $p=3$, then we can not use the Jacobian criterion for homogeneous polynomials. 
However, we can still apply a similar argument to the above after taking the standard affine cover of $\P^4_K$.

Assume $p=2$. 
We then have 
\[
\Delta_{\univ, K} = \{ \delta := \alpha \beta \gamma + a\alpha^2 +b \beta^2 +c \gamma^2 =0\} 
\subset \P^5_{K} 
= \Proj\,K[a, b, c, \alpha, \beta, \gamma]. 
\]
We have 
\[
\partial_a \delta =\alpha^2, \qquad \partial_b\delta  =\beta^2, \qquad \partial_c\delta  =\gamma^2,
\]
\[
\partial_{\alpha} \delta =\beta\gamma, \qquad \partial_{\beta}\delta  =\alpha\gamma, \qquad \partial_{\gamma}\delta  =\alpha\beta. 
\]
Therefore, the singular locus of $\Delta_{\univ, K}$ is given by 
$\{\alpha = \beta = \gamma =0\}$, which is two-dimensional. 
This completes the proof of Step \ref{s1-Delta-normal}.
\end{proof}

\begin{step}\label{s2-Delta-normal}
$\Delta_{\univ, K}$ is a normal prime divisor on $\P^5_K$. 
\end{step}

\begin{proof}[Proof of Step \ref{s2-Delta-normal}]
Since 
$\Delta_{\univ, K}$ is an effective Cartier divisor on $\P^5_{K}$, 
$\Delta_{\univ, K}$ is Cohen--Macaulay. 
By Step \ref{s1-Delta-normal}, $\Delta_{\univ, K}$ is normal by Serre's criterion. 
If $\Delta_{\univ, K}$ is not irreducible, then $\Delta_{\univ, K}$ would be non-normal, because two distinct prime divisors on $\P^5_{K}$ 
intersect along a three-dimensional non-empty closed subset. 
Therefore, $\Delta_{\univ, K}$ is a normal prime divisor. 
This completes the proof of Step \ref{s2-Delta-normal}.
\end{proof}
Step \ref{s2-Delta-normal} completes the proof of Proposition \ref{p-Delta-normal}.
\end{proof}

\begin{cor}
$\Delta_{\univ}$ is flat over $\Spec\,\Z$. 
\end{cor}

\begin{proof}
The assertion follows from the fact that 
$\Delta_{\univ} \to \Spec\,\Z$ is a equi-dimensional morphism (Proposition \ref{p-Delta-normal})
from a Cohen--Macaulay scheme to a regular scheme. 
\end{proof}

\subsection{Relative universal families}\label{ss-univ-Hilb-rel}

\begin{nota}\label{n-Hilb-Univ2}
\begin{enumerate}
\item Let $k$ be an algebraically closed field, 
let $C$ be a smooth projective curve over $k$, and 
let $E$ be a locally free sheaf on $C$ of rank $3$. 
\item Set $\theta(m) :=2m+1 \in \Z[m]$, 
which is the Hilbert polynomial of conics. 
\item We have the Hilbert scheme $\Hilb^{\theta}_{\P_C(E)/C}$ and 
let $\Univ^{\theta}_{\P_C(E)/C} \subset \P_C(E) \times_C \Hilb^{\theta}_{\P_C(E)/C}$ be the universal closed subscheme. 
\item 
The induced morphism 
\[
f_{\univ, C, E} :   \Univ^{\theta}_{\P_C(E)/C} \hookrightarrow \Hilb^{\theta}_{\P_C(E)/C} \times_{C} \P_C(E) \xrightarrow{ {\rm pr}_1} \Hilb^{\theta}_{\P_C(E)/C}
\]
is a generically smooth conic bundle. 
Set 
\[
\Delta_{\univ, C, E}  := \Delta_{f_{\univ, C, E}}\qquad {\rm and}\qquad 
\Sigma_{\univ, C, E}  := \Sigma_{f_{\univ, C, E}}. 
\]
\end{enumerate}
\end{nota}

\begin{thm}\label{t-rel-univ}
We use Notation \ref{n-Hilb-Univ2}. 
Then the following hold. 
\begin{enumerate}
\item The induced morphism $\pi : \Hilb^{\theta}_{\P_C(E)/C} \to C$ is a $\P^5$-bundle, i.e., 
there exists an open cover $C = \bigcup_{i \in I} C_i$ such that 
the induced morphism $\Hilb^{\theta}_{\P_C(E)/C} \times_C C_i$ is isomorphic to $\P^5 \times_k C_i$ for every $i \in I$. 
\item $\Hilb^{\theta}_{\P_C(E)/C}$ is a $6$-dimensional smooth projective variety. 
\item $\Univ^{\theta}_{\P_C(E)/C}$ is a $7$-dimensional smooth projective variety. 
\item $\Delta_{\univ, C, E}$ is a reduced normal divisor on $\Hilb^{\theta}_{\P_C(E)/C}$. 
\item For the singular locus $(\Delta_{\univ, C, E})_{{\rm sing}}$ of $\Delta_{\univ, C, E}$, we have the set-theoretic equality 
\[
(\Delta_{\univ, C, E})_{{\rm sing}} = \Sigma_{\univ, C, E}. 
\]
Furthermore, it holds that $\dim \Sigma_{\univ, C, E} \leq 3$. 
\end{enumerate}
\end{thm}

\begin{proof}
Fix a non-empty open subset $C'$ of $C$ such that $E' := E|_{C'} \simeq \MO_{C'}^{\oplus 3}$. 
Then we have the following diagrams in which all the squares are cartesian: 
\[
\begin{tikzcd}
\Univ^{\theta}_{\P_{C'}(E')/C'} \arrow[r, hook] \arrow[d, "f'"] & \Univ^{\theta}_{\P_{C}(E)/C} \arrow[d, "f_{\univ, C, E}"]\\
\Hilb^{\theta}_{\P_{C'}(E')/C'} \arrow[r, hook]  \arrow[d]& \Hilb^{\theta}_{\P_{C}(E)/C} \arrow[d]\\
C' \arrow[r, hook] & C. 
\end{tikzcd}\qquad 
\begin{tikzcd}
\Univ^{\theta}_{\P_{C'}(E')/C'} \arrow[r] \arrow[d, "f'"] & \Univ^{\theta}_{\P^2_k/k} \arrow[d, "f_{\univ} \times_{\Z} k"]\\
\Hilb^{\theta}_{\P_{C'}(E')/C'} \arrow[r]  \arrow[d]& \Hilb^{\theta}_{\P^2_k/k} \simeq \P^5_k \arrow[d]\\
C' \arrow[r] & \Spec\,k. 
\end{tikzcd}
\]
Then the assertion (1) holds by $\Hilb^{\theta}_{\P^2_k/k} \simeq \P^5_k$ (cf. Notation \ref{n-Hilb-Univ}). 

Let us show (2) and (3). 
By (1) and Proposition \ref{p-univ-sm}, both $\Hilb^{\theta}_{\P_C(E)/C} \to C$ and $\Univ^{\theta}_{\P_C(E)/C} \to C$ are smooth morphisms. 
Hence it is enough to show that {\cred $\Hilb^{\theta}_{\P_C(E)/C}$} and {\cred $\Univ^{\theta}_{\P_C(E)/C}$} are connected. 
Since $C$ and all the fibres of $\pi : \Hilb^{\theta}_{\P_C(E)/C} \to C$ are connected, 
also  $\Hilb^{\theta}_{\P_C(E)/C}$ is connected. 
Similarly, ${\cred \Univ^{\theta}_{\P_C(E)/C}}$ is connected, because so are 
$\Hilb^{\theta}_{\P_C(E)/C}$ and all the fibres of $\Univ^{\theta}_{\P_C(E)/C} \to \Hilb^{\theta}_{\P_C(E)/C}$. 
This completes the proof of (2) and (3).

The assertion (4) follows from 
Remark \ref{r-disc-bc} and  Proposition \ref{p-Delta-normal}. 
Let us show (5). 
Since $\Hilb^{\theta}_{\P_C(E)/C}$ and $\Univ^{\theta}_{\P_C(E)/C}$ are smooth varieties by (2) and (3), 
the set-theoretic equality $(\Delta_{\univ, C, E})_{{\rm sing}} = \Sigma_{\univ, C, E}$ follows from 
Theorem \ref{t-Delta-nonreg}. 
Since $\Delta_{\univ, C, E}$ is a reduced normal divisor on $\Hilb^{\theta}_{\P_C(E)/C}$, 
we obtain 
\[
\dim \Sigma_{\univ, C, E} = \dim (\Delta_{\univ, C, E})_{{\rm sing}} \leq 
\dim \Delta_{\univ, C, E} -2 \leq (\dim \Hilb^{\theta}_{\P_C(E)/C}-1) -2 =3. 
\]
Thus (5) holds. 
\qedhere


\end{proof}

\begin{lem}\label{l-rel-univ-N}
We use Notation \ref{n-Hilb-Univ2}. 
Let $\Gamma$ and $L$ be curves on $\Hilb^{\theta}_{\P_C(E)/C}$ 
such that $\pi(\Gamma) =C$ and $\pi(L)$ is a point. 
Then 
\[
N_1(\Hilb^{\theta}_{\P_C(E)/C}) \otimes_{\Z} \Q = \Q [\Gamma] + \Q[L],  
\]
where $N_1(\Hilb^{\theta}_{\P_C(E)/C}) := Z_1(\Hilb^{\theta}_{\P_C(E)/C})/\equiv$, 
$Z_1(\Hilb^{\theta}_{\P_C(E)/C}) := \bigoplus_{B:\text{curve}} \Z B$, 
$\equiv$ denotes the numerical equivalence, 
and $[\Gamma]$ and $[L]$ are the numerical equivalence classes. 
\end{lem}

\begin{proof}
For an ample Cartier divisor $H$ on $C$, 
we obtain $\pi^*H \cdot \Gamma >0$ and $\pi^*H \cdot L=0$, 
which imply that $[\Gamma]$ and $[L]$ are linearly independent over $\Q$. 
Hence it suffices to show that the following sequence is exact, because it implies $\rho(\Hilb^{\theta}_{\P_C(E)/C})=2$: 
\[
0 \to \Pic\,C \xrightarrow{\pi^*} \Pic\,(\Hilb^{\theta}_{\P_C(E)/C}) \xrightarrow{\cdot L } \Z. 
\]
By $\pi_*\MO_{\Hilb^{\theta}_{\P_C(E)/C}} = \MO_C$ (Theorem \ref{t-rel-univ}(1)), 
$\pi^* : \Pic\,C \to \Pic\,(\Hilb^{\theta}_{\P_C(E)/C})$ is injective. 
Pick $M \in \Pic\,(\Hilb^{\theta}_{\P_C(E)/C})$ satisfying  $M \cdot L=0$. By using Theorem \ref{t-rel-univ}(1), 
it is easy to see that $\pi_*M \in \Pic\,C$ and $\pi^*\pi_*M \to M$ is an isomorphism. 
\end{proof}

\subsection{Discriminant bundles}\label{s-disc-bdl}

Recall that the discriminant scheme $\Delta_f$ for a conic bundle $f: X \to T$ 
satisfies the following properties. 
\begin{enumerate}
    \item[(I)] For a point $t \in T$, $X_t$ is not smooth if and only if $t \not\in \Delta_f$ (Theorem \ref{t-disc-non-sm}). 
    \item[(II)] Let 
\[
\begin{CD}
X' @>\alpha >> X\\
@VVf'V @VVf V\\
T' @>\beta >> T
\end{CD}
\]
be a cartesian diagram of noetherian schemes such that $f$ is a conic bundle. 
In particular, also $f'$ is a conic bundle. 
Then the equality 
\[
\Delta^{\bdl}_{f'} = \beta^{-1}\Delta^{\bdl}_f
\]
of closed subschemes holds (Remark \ref{r-disc-bc}). 
\end{enumerate}

The purpose of this subsection is to introduce an invertible sheaf $\Delta_f^{\bdl}$ on $T$ 
that satisfies the following two properties. 
\begin{enumerate}
    \item[(I)'] Let $f : X \to T$ be a generically smooth conic bundle 
    such that $T$ is a noetherian integral scheme. 
    Then $\Delta^{\bdl}_f \simeq \MO_T(\Delta_f)$ (Remark \ref{r-disc-bdl-div}). 
    \item[(II)'] Let 
\[
\begin{CD}
X' @>\alpha >> X\\
@VVf'V @VVf V\\
T' @>\beta >> T
\end{CD}
\]
be a cartesian diagram, 
where $T$ and $T'$ are notherian integral schemes, 
$f$ is a conic bundle, and $\beta$ is of finite type. 
In particular, also $f'$ is a conic bundle. 
Then $\Delta^{\bdl}_{f'} \simeq \beta^*\Delta^{\bdl}_f$ (Theorem \ref{t-disc-bdl-bc}).  
\end{enumerate}
Although our definition (Definition \ref{d-disc-bdl}) might look unnatural at first sight, 
it is designed in order that $\Delta^{\bdl}_f$ satisfies (I)' and (II)'. 
Indeed, it is easy to see that there is {\cred the} unique way, if it exists, to define invertible sheaves $\Delta_f^{\bdl}$ 
satisfying (I)'and (II)' when  $f:X \to T$ is a conic bundle and $T$ is a regular noetherian integral scheme. 





\begin{dfn}\label{d-disc-bdl}
Let $f: X \to T$ be a conic bundle, 
where $T$ is a noetherian integral scheme. 
We then have the following cartesian diagram: 
\[
\begin{CD}
X @>>> \Univ^{\theta}_{\P(f_*\omega_{X/T}^{-1})/T}\\
@VVf V @VVgV\\
T @>j >> \Hilb^{\theta}_{\P(f_*\omega_{X/T}^{-1})/T},
\end{CD}
\]
where $\theta$ is the Hilbert polynomial of conics (Subsection \ref{ss:notation}(\ref{ss-nota-conic})) and 
$j$ denotes the morphism induced by 
the closed immersion $X \hookrightarrow \P(f_*\omega_{X/T}^{-1})$ 
(Proposition \ref{p-conic-emb}(3)). 
We set 
\[
\Delta^{\bdl}_{f}:= j^*\MO_T(\Delta_{g}), 
\]
which we call the {\em discriminant bundle} of $f$. 
Note that $\Delta_{g}$ is an effective Cartier divisor on $T$ (Lemma \ref{l-disc-bdl}) and 
$\MO_T(\Delta_{g})$ denotes the invertible sheaf corresponding to $\Delta_g$. 
In particular, $\Delta^{\bdl}_{f}$ is an invertible sheaf on $T$. 
\end{dfn}
\begin{lem}\label{l-disc-bdl}
We use the same notation as in Definition \ref{d-disc-bdl}. 
Then $\Delta_g$ is an effective Cartier divisor. 
\end{lem}

\begin{proof}
By the same argument as in Theorem \ref{t-rel-univ}(1)(2), we see that 
$\Hilb^{\theta}_{\P(f_*\omega_{X/T}^{-1})/T} \to T$ 
is a $\P^5$-bundle and $\Hilb^{\theta}_{\P(f_*\omega_{X/T}^{-1})/T}$ is a noetherian integral scheme. 
Furthermore, $\Univ^{\theta}_{\P(f_*\omega_{X/T}^{-1})/T} \to \Hilb^{\theta}_{\P(f_*\omega_{X/T}^{-1})/T}$ is a generically smooth conic bundle. 
Thus $\Delta_g$ is an effective Cartier divisor. 
\end{proof}

\begin{rem}\label{r-disc-bdl-div}
Let $f: X \to T$ be a conic bundle, 
where $T$ is a noetherian integral scheme. 
If $f$ is generically smooth, then 
it follows from Definition \ref{d-disc-bdl} that  $\Delta^{\bdl}_f \simeq \MO_T(\Delta_f)$. 
On the other hand, if $f$ is not generically smooth, 
then $\Delta_f = T$, which is no longer a Cartier divisor. 
\end{rem}

\begin{lem}\label{l-Hilb-bc}
Let $\beta : T' \to T$ be a morphism of noetherian schemes. 
Let $E$ be a locally free sheaf of rank $3$ on $T$. 
Set $E' := \beta^*E$. 
Then 
\[
\Hilb^{\theta}_{\P_{T'}(E')/T'} \simeq \Hilb^{\theta}_{\P_T(E)/T} \times_T T', 
\]
where $\theta := 2m+1 \in \Z[m]$, which is the Hilbert polynomial of conics. 
\end{lem}

\begin{proof}
The assertion follows from $\P_{T'}(E') \simeq \P_T(E) \times_T T'$. 
\end{proof}

\begin{thm}\label{t-disc-bdl-bc}
Let 
\[
\begin{CD}
X' @>\alpha >> X\\
@VVf'V @VVf V\\
T' @>\beta >> T
\end{CD}
\]
be a cartesian diagram of noetherian schemes, where $T$ and $T'$ are noetherian integral schemes and $f$ is a conic bundle. 
Assume that one of the following holds. 
\begin{enumerate}   
\item[(A)] $\beta$ is flat. 
    \item[(B)] $\beta$ is a closed immersion such that $\beta(T')$ is an effective Cartier divisor on $T$. 
    \item[(C)] Both $T$ and $T'$ are regular and $\beta$ is of finite type.  
\end{enumerate}
Then the following hold. 
\begin{enumerate}
    \item The induced homomorphism 
    \[
    \theta : \beta^*f_*\omega_{X/T}^{-1} \to  f'_*\alpha^*\omega^{-1}_{X/T}
    \]
    is an isomorphism. 
    \item 
    An isomorphism $\Delta^{\bdl}_{X'/T'} \simeq \beta^*\Delta^{\bdl}_{X/T}$ holds. 
\end{enumerate}
\end{thm}

\begin{proof}
By Remark \ref{r-disc-bc} and Lemma \ref{l-Hilb-bc}, (1) implies (2). 
Hence it suffices to prove (1). We may assume that $T$ and $T'$ are affine. 
If (A) holds, then the assertion (1) immediately follows from the flat base change theorem. 

Assume (B). 
Recall that 
$\beta^*f_*\omega_{X/T}^{-1}$ is a locally free sheaf of rank $3$ (Proposition \ref{p-conic-emb}). 
By the same argument as in Proposition \ref{p-conic-emb}, 
also $f'_*\alpha^*\omega^{-1}_{X/T}$ is a locally free sheaf of rank $3$. 
Then it is enough to show that 
$\theta : \beta^*f_*\omega_{X/T}^{-1} \to f'_*\alpha^*\omega^{-1}_{X/T}$ is surjective. 
Hence the problem is reduced to the surjectivity of 
\[
H^0(X, \omega_{X/T}^{-1}) \to H^0(X', \alpha^*\omega_{X/T}^{-1}). 
\]
We have an exact sequence: 
\[
0 \to \omega^{-1}_{X/T} \otimes I_{X'} \to \omega^{-1}_{X/T} \to \alpha_*\alpha^*\omega^{-1}_{X/T} \to 0. 
\]
By $I_{X'} =\MO_{X}(-X') = f^*\MO_T(-T')$, we obtain 
\[
R^1f_*(\omega^{-1}_{X/T} \otimes I_{X'}) 
\simeq 
R^1f_*(\omega^{-1}_{X/T} \otimes f^*\MO_T(-T')) 
\simeq 
R^1f_*(\omega^{-1}_{X/T}) \otimes \MO_T(-T') \overset{{\rm (*)}}{=}0, 
\]
where ${\rm (*)}$ follows from Proposition \ref{p-conic-emb}. 
This completes the proof of (1) when (B) holds. 

Assume (C). 
Since the problem is local on $T$ and $T'$, 
the problem is reduced to the case when 
$\beta : T'\to T$ is factored as follows: 
\[
T =:T_0 \hookrightarrow T_1 \hookrightarrow \cdots \hookrightarrow 
T_r := T' \times \mathbb A^n \xrightarrow{{\rm pr}_1} T', 
\]
where each $T_i$ is a regular affine noetherian integral scheme and 
each $T_i \hookrightarrow T_{i+1}$ is a closed immersion with $\dim T_i = \dim T_{i+1} -1$ \cite[Theorem 14.2 and Theorem 21.2]{Mat86}. 
Therefore, the problem is reduced to the case when (A) or (B) holds. 
This completes the proof of  (1).  
\qedhere


\end{proof}


\subsection{Proof of the Mori--Mukai formula}\label{ss-MM-proof}

\begin{lem}\label{l-sm-surface}
We work over an algebraically closed field $k$. 
Let $f: X \to T$ be a conic bundle, 
where $X$ is a smooth surface and $T$ is a smooth curve. 
Then $f$ is generically smooth and 
every singular fibre of $f$ is reduced and consists of exactly two $(-1)$-curves.  
\end{lem}

\begin{proof}
Let $X'$ be a relatively minimal model of $f$, so that we obtain the induced morphisms 
\[
f : X \xrightarrow{\varphi} X' \xrightarrow{f'} T. 
\]
Then the assertion follows from the fact that 
$f' : X' \to T$ is a $\P^1$-bundle and $\varphi : X \to X'$ is a sequence of blowups. 
\end{proof}

{\cred 
We are ready to prove the Mori-Mukai formula in arbitrary characteristics, 
which  generalises the case of characteristic zero 
\cite[Proposition 6.2]{MM83}.}

\begin{thm}\label{t-cbf}
We work over an algebraically closed field $k$. 
Let $f: X\to T$ be a 
conic bundle, 
where $T$ is a smooth projective variety. 
Then the equality 
\begin{equation}\label{e1-cbf}
\Delta^{\bdl}_f \cdot C = - \omega_{X/T} \cdot \omega_{X/T} \cdot f^{-1}(C) 
\end{equation}
holds for every curve $C$ on $T$. 
In particular,  when also $X$ is smooth over $k$, the numerical equivalence 
\[
\Delta^{\bdl}_f \equiv -f_*(K_{X/T}^2)
\]
holds for a Cartier divisor $K_{X/T}$ on $X$ satisfying $\omega_{X/T} \simeq \MO_X(K_{X/T})$. 
\end{thm}

\begin{proof}
The in-particular part follows from the projection formula. 
Hence it suffices to show (\ref{e1-cbf}). 
In what follows, we write $\mathcal L^2 := \mathcal L \cdot \mathcal L$ 
(note that $\mathcal L^2$ differs from $\mathcal L^{\otimes 2}$). 

\setcounter{step}{0}
\begin{step}\label{s1-cbf}
Let $C$ be a projective curve on $T$. 
Let $C^N \to C$ be the normalisation of $C$ and let $g : X \times_T C^N \to C^N$ 
be the base change of $f$. 
Then the following hold. 
\begin{enumerate}
    \item $\Delta^{\bdl}_f \cdot C = \deg \Delta^{\bdl}_g$. 
    \item $\omega_{X/T}^{2} \cdot f^{-1}(C) = \omega^2_{X \times_T C^N/C^N}$. 
\end{enumerate}
\end{step}

\begin{proof}[Proof of Step \ref{s1-cbf}]
We have the following cartesian diagram: 
\[
\begin{tikzcd}
X \times_T C^N \arrow[d, "g"] \arrow[r, "\alpha"] & X \arrow[d, "f"]\\
C^N \arrow[r, "\beta"] & T.  
\end{tikzcd}
\]
The assertion (1) holds by  
\[
\Delta^{\bdl}_f \cdot C = \deg (\beta^*\Delta^{\bdl}_f) = \deg \Delta^{\bdl}_g, 
\]
where the latter equality follows from 
$\beta^*\Delta^{\bdl}_f \simeq  \Delta^{\bdl}_g$ (Theorem \ref{t-disc-bdl-bc}(2)). 
We obtain  (2)  by the following: 
\[
\omega_{X/T}^{2} \cdot f^{-1}(C) 
\overset{{\cred (\star)}}{=}(\omega_{X/T}|_{f^{-1}(C)})^2 
=(\alpha^*\omega_{X/T})^2 
\overset{{\cred (\star\star)}}{=}
\omega_{X \times_T C^N/C^N}^2, 
\]
where 
{\cred $(\star)$ holds by \cite[Lemma 1.10]{Bad01} and $(\star\star)$} 
follows from 
$\alpha^*\omega_{X/T} \simeq \omega_{X \times_T C^N/C^N}$ \cite[Theorem 3.6.1]{Con00}. 
{\cred 
For the definition and some properties of the intersection number $(\omega_{X/T}|_{f^{-1}(C)})^2$, 
we refer to \cite[Section 1]{Bad01}.}  
This completes the proof of Step \ref{s1-cbf}. 
\end{proof}

\begin{step}\label{s2-cbf}
The assertion of Theorem \ref{t-cbf} holds 
when $T$ is a smooth projective curve, $f$ is generically smooth, and $\Delta_f$ is reduced. 
\end{step}

\begin{proof}[Proof of Step \ref{s2-cbf}]
By Theorem \ref{t-sm-red-fib}, $X$ is smooth over $k$ and any fibre of $f$ is reduced, i.e., each singular fibre consists of two $(-1)$-curves (Lemma \ref{l-sm-surface}). 
Let 
\[
f : X \xrightarrow{\mu} X' \xrightarrow{f'} T
\]
be a relatively minimal model. 
Then any fibre of $f' : X' \to T$ is $\P^1$. 
Set $e$ to be the number of the singular fibres of $f: X \to T$. 
Then we have 
\[
K_X^2 =K_{X'}^2 -e. 
\]
For the genus $g(T)$ {\cred of $T$}, it holds that 
\[
K_{X'/T}^2 = (K_{X'} -f'^*K_T)^2 = K_{X'}^2 -2 K_{X'} \cdot f'^*K_T 
\]
\[
=8(1-g(T)) -2 \cdot (-2) \cdot (2g(T)-2)=0. 
\]
By $\Delta_f^{\bdl} \simeq \MO_T(\Delta_f)$ (Remark \ref{r-disc-bdl-div}), we obtain 
\[
\deg \Delta^{\bdl}_f =\deg \Delta_f = e = K_{X'}^2 - K_X^2 = K_{X'/T}^2 -K_{X/T}^2 = -K_{X/T}^2. 
\]
This completes the proof of Step \ref{s2-cbf}. 
\end{proof}

\begin{step}\label{s3-cbf}
The assertion of Theorem \ref{t-cbf} holds when  $f: X \to T$ coincides with 
the induced morphism 
\[
f_{\univ, C, E} :   \Univ^{\theta}_{\P_C(E)/C} \to \Hilb^{\theta}_{\P_C(E)/C}
\]
for some smooth projective curve $C$ and locally free sheaf $E$ on $C$ of rank $3$. 
\end{step}

\begin{proof}[Proof of Step \ref{s3-cbf}]
We have the projection
\[
\pi : \Hilb^{\theta}_{\P_C(E)/C} \to C. 
\]
Recall that there exists an open cover $C = \bigcup_{i \in I} C_i$ such that $\pi^{-1}(C_i) \simeq \P^5 \times C_i$ 
(Theorem \ref{t-rel-univ}). 
Take a general line $L$ contained in  a fibre of $\pi$. 
Fix a smooth projective curve $T'$ obtained as a complete intersection of general hyperplane sections of $T = \Hilb^{\theta}_{\P_C(E)/C}$. 
It follows from Lemma \ref{l-rel-univ-N} that 
\[
N_1(T) \otimes_{\Z} \Q = \Q [T'] + \Q[L]. 
\]
Therefore, it is enough to show 
\[
\Delta^{\bdl}_f \cdot T'=
-\omega_{X/T}^2 \cdot f^{-1}(T')  \qquad \text{and}\qquad 
\Delta^{\bdl}_f \cdot L = -\omega_{X/T}^2 \cdot f^{-1}(L). 
\]

We now prove $\Delta^{\bdl}_f \cdot T'=
-\omega_{X/T}^2 \cdot f^{-1}(T')$. 
Recall that $\Delta_{\univ, C, E}$ is a reduced divisor on 
a $6$-dimensional smooth projective variety $T=\Hilb^{\theta}_{\P_C(E)/C}$ (Theorem \ref{t-rel-univ}). 
Hence the scheme-theoretic intersection $T' \cap \Delta_{\univ, C, E}$  is  smooth, 
because the curve $T'$ is defined as a complete intersection of general hyperplane sections. 
For the base change $f' : X \times_T T' \to T'$ of $f$, we have $\Delta_{f'} = T' \cap \Delta_f = T' \cap \Delta_{\univ, C, E}$ (Remark \ref{r-disc-bc}), and hence $\Delta_{f'}$ is smooth. 
Then 
\[
\Delta^{\bdl}_f \cdot T' \overset{{\rm (i)}}{=} \deg \Delta_{f'} \overset{{\rm (ii)}}{=}  -\omega_{X \times_T T'/T'}^2 \overset{{\rm (iii)}}{=} -\omega_{X/T}^2 \cdot f^{-1}(T'), 
\]
where (i) and (iii) hold by Step \ref{s1-cbf}(1) and Step \ref{s1-cbf}(2) respectively, and 
(ii) follows from Step \ref{s2-cbf}. 
This completes the proof of $\Delta^{\bdl}_f \cdot T'=
-\omega_{X/T}^2 \cdot f^{-1}(T')$. 

It suffices to show $\Delta^{\bdl}_f \cdot L =-\omega_{X/T}^2 \cdot f^{-1}(L)$. 
By the same argument as in the previous paragraph, it is enough to prove that 
the scheme-theoretic intersection $L \cap \Delta_{\univ, C, E}$  is  smooth. 
Let $F (\simeq \P^5_k)$ be the fibre of $\pi : \Hilb^{\theta}_{\P_C(E)/C} \to C$ containing $L$. 
There is $i \in I$ such that $F \hookrightarrow \pi^{-1}(C_i)$. 
Consider the following composite isomorphism:  
\[
\iota: F \hookrightarrow \pi^{-1}(C_i) \simeq \P^5 \times C_i \xrightarrow{{\rm pr}_1} \P^5. 
\]
We have $\iota^*\Delta_{\univ, k} = \Delta_{\univ, C, E}|_F$, because  
$\Delta_{\univ, C, E}|_F$ is the pullback of $\Delta_{\univ, C, E}|_{\pi^{-1}(C_i)}$. 
It follows from Proposition \ref{p-Delta-normal} that the scheme-theoretic intersection $M \cap \Delta_{\univ, k}$ 
is smooth for a general line $M$ on $\P^5$. 
As $L$ is chosen to be a general line on $F \simeq \P^5$, also $L \cap \Delta_{\univ, C, E}$ is smooth. 
This completes the proof of Step \ref{s3-cbf}. 
\end{proof}

\begin{step}\label{s4-cbf}
The assertion of Theorem \ref{t-cbf} holds when  $T$ is a smooth projective curve. 
\end{step}

\begin{proof}[Proof of Step \ref{s4-cbf}]
Set $E := f_*(\omega_{X/T}^{-1})$. 
Then $E$ is a locally free sheaf of rank $3$ (Proposition \ref{p-conic-emb}(2)) and 
$f : X \hookrightarrow \P_T(E) \to T$ is a flat family of conics on $\P_T(E)$ relatively over $T$. 
Therefore, we have the following cartesian diagram 
\[
\begin{tikzcd}[xscale=2]
X \arrow[d, "f"]\arrow[r, "\alpha"] & \Univ^{\theta}_{\P_T(E)/T}  \arrow[d, "f_{\univ, T, E}"]\\
T\arrow[r, "\beta"] & \Hilb^{\theta}_{\P_T(E)/T}. 
\end{tikzcd}
\]
Note that $\beta : T \to \Hilb^{\theta}_{\P_T(E)/T}$ is a section of the projection 
$\Hilb^{\theta}_{\P_T(E)/T} \to T$. 
In particular, $\beta : T \to \Hilb^{\theta}_{\P_T(E)/T}$ is a closed immersion. 
We then obtain 
\[
\deg \Delta^{\bdl}_f \overset{{\rm (i)'}}{=} \Delta^{\bdl}_{f_{\univ, T, E}} \cdot \beta(T) \overset{{\rm (ii)'}}{=} 
- \omega^2_{\Univ^{\theta}_{\P_T(E)/T}/\Hilb^{\theta}_{\P_T(E)/T}} \cdot 
f_{\univ, T, E}^{-1}(\beta(T))
\overset{{\rm (iii)'}}{=} -\omega_{X/T}^2, 
\]
where ${\rm (i)'}$ and ${\rm (iii)'}$ hold by Step \ref{s1-cbf}(1) and Step \ref{s1-cbf}(2) respectively, and 
${\rm (ii)'}$ follows from Step \ref{s3-cbf}. 
This completes the proof of Step \ref{s4-cbf}. 
\end{proof}

\begin{step}\label{s5-cbf}
Theorem \ref{t-cbf} holds without any additional assumptions. 
\end{step}

\begin{proof}[Proof of Step \ref{s5-cbf}]
It suffices to show that $\Delta_f \cdot C= -f_*\omega_{X/T}^2 \cdot C$ for any curve $C$ on $T$. 
This follows from  
\[
\Delta^{\bdl}_f \cdot C\overset{{\rm (i)''}}{=}  \deg \Delta^{\bdl}_{h} 
\overset{{\rm (ii)''}}{=}  - \omega_{X \times_T C^N/C^N}^2
\overset{{\rm (iii)''}}{=}  - \omega_{X/T}^2 \cdot f^{-1}(C), 
\]
where $h: X \times_T C^N \to C^N$ denotes the base change of $f :X \to T$, 
${\rm (i)''}$ and ${\rm (iii)''}$ hold by Step \ref{s1-cbf}(1) and Step \ref{s1-cbf}(2) respectively, and 
${\rm (ii)''}$ follows from Step \ref{s4-cbf}. 
This completes the proof of Step \ref{s5-cbf}. 
\end{proof}
Step \ref{s5-cbf} completes the proof of Theorem \ref{t-cbf}. 
\qedhere



\end{proof}

\section{Surfaces}\label{s-surface}

In this section, we focus on surface conic bundles $f : X \to T$. 
More precisely, we treat the case when $X$ has at worst canonical singularities, 
since there is nothing to do for the case when $X$ is smooth over $k$ (Lemma \ref{l-sm-surface}). 
The motivation is to seek a relation between $\Delta_f$ and the singularities of $X$. 
For this purpose, we assume that $f: X \to T$ has {\cred the} unique singular fibre $X_0$, 
where $0 \in T$ is a closed point. 
The primary goal of this section is to establish the following results (Theorem \ref{t-disc-vs-sing} in Subsection \ref{ss-disc-vs-sing}). 
\begin{enumerate}
    \item $\deg \Delta_f = m-1$, where $m$ denotes the number of the {\cred irreducible}  components of the fibre $Y_0$, 
    where $Y \to X$ is the minimal resolution of $X$. 
    \item 
    If $X_0$ is reduced, then 
    \begin{itemize}
        \item $X$ has {\cred the} unique singularity $x$, and 
        \item $x$ is of type $A_n$, where $n := m-2 = \deg \Delta_f -1$. 
    \end{itemize}
\item 
If $X_0$ is not reduced, then one and only one of the following holds. 
\begin{enumerate}
    \item $X$ has exactly two  singularities $x$ and $x'$. Moreover, both $x$ and $x'$ are of type $A_1$. In this case, $\deg \Delta_f = 2$. 
    \item 
    $X$ has {\cred the} unique singularity $x$. Moreover, 
    $x$ is of type $D_n$ with $n \geq 3$, where $n= m-1 = \deg \Delta_f$. 
\end{enumerate}
\end{enumerate}
For example, $X$ has {\cred the} unique singularity of type $D_4$
if $\deg \Delta_f =4$ and $X_0$ is not reduced (note that each of these conditions can be checked, as far as $f: X \to T$ is explicitly given). 
In order to establish the above properties (1)--(3), 
we first classify the dual graphs of the exceptional loci of the minimal resolutions of 
surface conic bundles (Subsection \ref{ss-surface-char0}). 
In Subsection \ref{ss-2dim-classify}, 
we  shall summarise further classification results for surface conic bundles and exhibit several examples.

\subsection{The classification of dual graphs}\label{ss-surface-char0}

In this subsection, we classify the singularities of surface conic bundles for the case when 
the total space $X$ has at worst canonical singularities 
(Proposition \ref{p-dim2-red-sing}, Proposition \ref{p-dim2-nonred-sing}). 
The results in this subsection should be well known to experts. 
We include the proofs for the reader's convenience. 
{\cred We refer to \cite[\S 4.1, cf. 4.19]{KM98} for some foundational results on canonical surface singularities.} 



\begin{prop}\label{p-dim2-red-sing}
We work over an algebraically closed field $k$. 
Let $f: X \to T$ be a conic bundle, 
where $X$ is a canonical surface and $T$ is a smooth curve. 
Fix a closed point $0 \in T$. 
Assume that $X_0$ is reduced and $X_0$ is {\cred the} unique singular fibre. 
Let $\varphi : Y \to X$ be the minimal resolution of $X$. 
Then the dual graph of the singular fibre of the composition $Y \xrightarrow{\varphi} X \xrightarrow{f} T$ is as follows: 
\[
(-1) - (-2) - (-2) - \cdots - (-2) - (-1). 
\]
In particular, $X$ has {\cred the} unique singularity $x$ and $x$ is of type $A_n$ for some $n \in \Z_{>0}$. 
\end{prop}


\begin{proof}
We run a $K_Y$-MMP over $T$: 
\[
\psi : Y =: Z_m \xrightarrow{\psi_{m-1}} Z_{m-1}  \xrightarrow{\psi_{m-2}} 
\cdots 
\xrightarrow{\psi_2} 
Z_2 
\xrightarrow{\psi_1}  Z_1 =: Z, 
\]
i.e., each $\psi_i : Z_{i+1} \to Z_i$ is a contraction of a $(-1)$-curve contained in the fibre over $0 \in T$. 
After possibly replacing $T$ by an open neighbourhood of $0 \in T$, 
we may assume that $Z_1 =Z = \P^1 \times T$. 
Then the singular fibre $(Z_2)_0$ consists of exactly two $(-1)$-curves $E_1+E_2$.


It is enough to prove the following three properties (1)-(3). 
\begin{enumerate}
    \item The blowup centre of $\psi_i  : Z_{i+1} \to Z_i$ is a smooth point of the fibre $(Z_i)_0$ for every $i$. 
    \item The fibre $(Z_i)_0$ over $0$ is reduced for every  $i$. 
    \item The blowup centre of $\psi_i:Z_{i+1} \to Z_i$ 
    is {\cred disjoint} from     any $(-2)$-curve contained in  $(Z_i)_0$. 
\end{enumerate}
Indeed, if the singular fibre $(Z_i)_0$ is a chain 
\[
(-1) - (-2) - (-2) - \cdots - (-2) - (-1), 
\]
which is reduced, then the next blowup centre lies on one of the leftmost and rightmost $(-1)$-curves and avoids  the $(-2)$-curves. 

Let us show (1). 
Suppose that $(Z_i)_0$ is reduced and the blowup centre of $\psi_i : Z_{i+1} \to Z_i$ 
is a singular point of $(Z_i)_0$, i.e., it is an intersection of two irreducible components. 
Then $(Z_{i+1})_0$ is not reduced along the resulting  $(-1)$-curve $F_{i+1}$. 
Then we can show, by induction on $j$, that 
there is a $(-1)$-curve $F_j$ on $Z_j$ such that $(Z_j)_0$ is non-reduced along $F_j$. 
This {\cred implies} that $Y=Z_m$ contains a $(-1)$-curve $F$ inside the singular fibre $Y_0$ 
such that $Y_0$ is not reduced along $F$. 
Then $X_0$ is non-reduced along $\varphi(F)$. 
This contradicts our assumption. 
Thus (1) holds. 
By induction on $i$, the assertion (2) follows from (1). 

Let us show (3). 
Note that the number of the irreducible components of $(Z_i)_0$ is equal to $i$. 
Let $\nu_i$ be the number of $(-2)$-curves inside the fibre $(Z_i)_0$, e.g., 
we have $\nu_1 = \nu_2 = 0$ and $\nu_3 = 1$. 
For $i \geq 2$, let us show that $\nu_i = i-2$ by induction on $i$. 
Fix $i \geq 2$ and assume $\nu_i = i-2$. 
By (1), we obtain 
\[
\nu_{i+1}  \leq \nu_i +1. 
\]
Furthermore, the equality holds if and only if the blowup centre lies on a $(-1)$-curve. We obtain 
\[
m-2 \leq \nu_m \leq \nu_2 +m-2 = m-2, 
\]
where the first inequality follows from the fact that $X$ has 
two irreducible components and 
$\Ex(\varphi)$ consists of $(m-2)$ irreducible components. 
Therefore, 
we obtain 
$\nu_{i+1}  = \nu_i +1$. 
This completes the proof of $\nu_i =i-2$ for $i \geq 2$. 
Then we get $\nu_{i+1} = \nu_i +1$, i.e., 
the blowup centre of $\psi_i : Z_{i+1} \to Z_i$ lies on a $(-1)$-curve. 
Thus (3) holds. 
\qedhere



\end{proof}



\begin{prop}\label{p-dim2-nonred-sing}
We work over an algebraically closed field $k$. 
Let $f: X \to T$ be a conic bundle, 
where $X$ is a canonical surface and $T$ is a smooth curve. 
Fix a closed point $0 \in T$. 
Assume that $X_0$ is non-reduced and $X_0$ is {\cred the} unique singular fibre. 
Let $\varphi : Y \to X$ be the minimal resolution of $X$. 
Then the dual graph of the singular fibre of the composition $Y \xrightarrow{\varphi} 
X \xrightarrow{f} T$ is either a chain 
\[
(-2) - (-1) - (-2)
\]
or
\[
D_n - (-1) 
\]
with $n \geq 3$, where the latter case means that $(-1)$-curve intersects only with the long tail. 
\end{prop}

\begin{proof}
We run a $K_Y$-MMP over $T$: 
\[
\psi : Y =: Z_m \xrightarrow{\psi_{m-1}} Z_{m-1}  \xrightarrow{\psi_{m-2}} 
\cdots 
\xrightarrow{\psi_2} 
Z_2 
\xrightarrow{\psi_1}  Z_1 =: Z, 
\]
i.e., each $\psi_i : Z_{i+1} \to Z_i$ is a contraction of a $(-1)$-curve contained in the fibre over $0 \in T$. 
After possibly replacing $T$ by an open neighbourhood of $0 \in T$, 
we may assume that $Z_1 =Z = \P^1 \times T$. 
It is clear that $(Z_i)_0$ contains at least one $(-1)$-curve for every $i \geq 2$. 
Since 
$X_0$ is irreducible and $\Ex(\varphi : Y=Z_m \to X)$ consists of $(-2)$-curves, 
the following holds: 
\begin{enumerate}
\item We have the irreducible decomposition $((Z_m)_0)_{\red} = A_1 \cup \cdots \cup A_{m-1} \cup B$, 
where each $A_i$ is a $(-2)$-curve and $B$ is a $(-1)$-curve. 
\end{enumerate}
We see that 
\begin{enumerate}
\item[(2)] $(Z_i)_0$ contains no prime divisor $C$ satisfying $C^2 \leq -3$,  
\end{enumerate}
as otherwise its proper transform $C' \subset (Z_m)_0$ would satisfy $C'^2 \leq -3$, which contradicts (1). 

We now show that the blowup centre of $\psi_2 : Z_3 \to Z_2$ is the singular point $E_1 \cap E_2$ for {\cred the irreducible decomposition} $(Z_2)_0 = E_1 \cup E_2$. 
Otherwise, $(Z_3)_0$ contains two $(-1)$-curves $F$ and $F'$ which are mutually disjoint. 
This property is stable under taking a blowup. 
Hence also $(Z_m)_0$ has at least two $(-1)$-curves, which contradicts (1). 
Therefore, the dual graph of $(Z_3)_0$ is given by 
\[
(-2) - (-1) - (-2). 
\]

The blowup centre of the next blowup $Z_4 \to Z_3$ is disjoint from the two $(-2)$-curves by (2). 
Then $(Z_4)_0$ consists of three $(-2)$-curves and {\cred the} unique $(-1)$-curve. 
The blowup centre of $\psi_4 : Z_5 \to Z_4$ is disjoint from these $(-2)$-curves. 
Repeating this procedure, we see that the resulting dual graph of $Y_0 = (Z_m)_0$ is as in the statement.
\qedhere


\end{proof}


\subsection{Discriminants vs singularities}\label{ss-disc-vs-sing}

We are ready to prove a main theorem of this section. 

\begin{thm}\label{t-disc-vs-sing}
We work over an algebraically closed field $k$. 
Let $f: X \to T$ be a conic bundle, 
where $X$ is a canonical surface and $T$ is a smooth curve. 
Assume that there exists a closed point $0 \in T$ such that 
$X_0$ is {\cred the} unique singular fibre. 
Let $\varphi : Y \to X$ be the minimal resolution of $X$ and 
set $m$ to be the number of the irreducible components of the fibre $Y_0$ over $0 \in T$. 
Then the following hold. 
\begin{enumerate}
\item $\deg \Delta_f = m-1$. 
\item 
Assume that $X_0$ is reduced. 
Then $X$ has {\cred the} unique singularity $x$ and $x$ is of type $A_n$, 
where $n= m-2 = \deg \Delta_f -1$. 
\item 
Assume that $X_0$ is not reduced. 
Then one and only one of the following holds. 
\begin{enumerate}
    \item $X$ has exactly two  singularities $x$ and $x'$. Moreover,  both $x$ and $x'$ are of type $A_1$. In this case, $\deg \Delta_f = 2$. 
    \item 
    $X$ has {\cred the} unique singularity $x$ and $x$ is of type $D_n$ with $n \geq 3$, 
where $n= m-1 = \deg \Delta_f$. 
\end{enumerate}
\end{enumerate}
\end{thm}

\begin{proof}
Let us show (1). 
By compactifying $T$ suitably, we may assume that $T$ is a smooth projective curve. 
We run a $K_Y$-MMP: $\psi: Y \to Z$ over $T$, so that $f_Z : Z \to T$ is a $\P^1$-bundle. 
Let $f_Y : Y \to T$ be the induced morphism. 
By $\deg \Delta_f = -K_{X/T}^2$ ({\cred Remark \ref{r-disc-bdl-div},} Theorem \ref{t-cbf}), we obtain  
\[
\deg \Delta_f = -K_{X/T}^2 = -(K_X-f^*K_T)^2 = -(K_Y -f_Y^*K_T)^2
\]
\[
= -(K_Z-f_Z^*K_T)^2 + (m-1) = m-1. 
\]
Thus (1) holds. 
By (1), the assertions (2) and (3) follow from 
Proposition \ref{p-dim2-red-sing} and Proposition \ref{p-dim2-nonred-sing}, respectively. 
\end{proof}


\subsection{Description and examples}\label{ss-2dim-classify}
In this subsection, we study generically smooth conic bundles over $k[[t]]$. 
We shall classify such conic bundles 
for the case when either 
\begin{itemize}
\item 
${\rm char}\,k \neq 2$ (Proposition \ref{p-surface-char0-red}, Proposition \ref{p-surface-char0-nonred}), or 
\item ${\rm char}\,k = 2$ and the fibre $X_0$ over the closed point is reduced
(Proposition \ref{p-surface-char2-red}). 
\end{itemize}
The remaining case is when ${\rm char}\,k = 2$ and the fibre $X_0$ is non-reduced. 
Although we will not give the complete classification in this case, 
we shall check, by exhibiting several examples, that 
all the canonical (RDP) singularities of type $A$ and $D$ in Artin's list \cite{Art77}  actually appear.  


\subsubsection{${\rm char}\,k \neq 2$}
In Proposition \ref{p-surface-char0-red} (resp. Proposition \ref{p-surface-char0-nonred}), 
we treat the case when the singular fibre is reduced (resp. non-reduced).




\begin{prop}\label{p-surface-char0-red}
Let $k$ be an algebraically closed field of characteristic $\neq 2$. 
Set $A := k[[t]]$ and let $0$ be the closed point of $\Spec\,A$. 
Let $f: X 
\to \Spec\,A$ be a conic bundle, 
where $X$ is canonical and $X_0$ is not smooth but {\cred is} reduced. 
Then 
the following hold. 
\begin{enumerate}
    \item We have an $A$-isomorphism 
    \[
    X \simeq \Proj\,A[x, y, z]/(x^2 + y^2 +t^{n+1}z^2)
    \]
    for some $n \geq 0$. 
        \item 
    $\deg \Delta_f = n+1$. 
    \item 
    If $n=0$, then $X$ is regular. 
    If $n\geq 1$, then $X$ has {\cred the} unique singularity and it is of type $A_{n}$. 
\end{enumerate}
\end{prop}

\begin{proof}
Let us show (1). 
By Corollary \ref{c-char0-local2}(2), we can write 
\[
X \simeq \Proj\,A[x, y, z]/(x^2 + y^2 +cz^2)
\]
for some $c \in A = k[[t]]$. 
Then we have $c = dt^m$ for some $m \in \Z_{\geq 0}$ and $d \in k[[t]]^{\times}$. 
Since $X_0$ is not smooth, we get $m >0$. 
By Hensel's lemma, we have $\sqrt{d} \in k[[t]]$. 
Thus (1) holds. 

The assertion (2) follows from Definition \ref{d-delta} and Definition \ref{d-disc}. 
The assertion (3) holds by either direct computation, or 
(2) and Theorem \ref{t-disc-vs-sing}. 
\end{proof}




\begin{prop}\label{p-surface-char0-nonred}
Let $k$ be an algebraically closed field of characteristic $\neq 2$. 
Set $A := k[[t]]$ and let $0$ be the closed point of $\Spec\,A$. 
Let $f: X 
\to \Spec\,A$ be a conic bundle, 
where $X$ is canonical and $X_0$ is not  reduced. 
Then 
the following hold. 
\begin{enumerate}
    \item We have an $A$-isomorphism 
    \[
    X \simeq \Proj\,A[x, y, z]/(x^2 + ty^2 +t^{n-1}z^2)
    \]
    for some $n \geq 2$. 
        \item 
    $\deg \Delta_f = n$. 
        \item 
    If $n=2$, then $X$ has exactly two  singular points and both of them are of type $A_1$. 
    If $n \geq 3$, then $X$ has {\cred the} unique {\cred singular} point and it is of type $D_n$. 
\end{enumerate}
\end{prop}

\begin{proof}
Let us show (1). 
By Corollary \ref{c-char0-local2}(1), we can write 
\[
X = \Proj\,A[x, y, z]/(x^2 + by^2+cz^2+\alpha yz)
\]
for some $b, c, \alpha \in \m := tk[[t]]$.

We now reduce the problem to the case when $\alpha =0$. 
For a suitable $m \in \Z_{>0}$, we have 
\[
by^2+cz^2+\alpha yz = t^m(b'y^2 +c'z^2 + \alpha' yz), 
\]
where $b', c', \alpha' \in A$ and one of $b', c', \alpha'$ is a unit of $A$. 
If $b' \in A^{\times}$, then we may assume that $b' =1$ (Hensel's lemma) and 
$\alpha'=0$ by completing a square: 
$y^2 +c'z^2 + \alpha' yz = (y + \frac{\alpha'}{2}z)^2 + (c' - \frac{\alpha'^2}{4})z^2$. 
If $c' \in A^{\times}$, then the problem is reduced to the case when $b' \in A^{\times}$ by switching $y$ and $z$. 
If $b' \in \m$, $c' \in \m$, and $\alpha' \in A^{\times}$, then
we may assume that $b' \in A^{\times}$ by applying a linear transform 
$(x, y, z) \mapsto (x, y, y+z)$. 
In any case,  the problem is reduced to the case when $\alpha=0$. 

After possibly switching $y$ and $z$, Hensel's lemma enables us to write 
\[
X \simeq  \Proj\,A[x, y, z]/(x^2 + t^r(y^2+t^sz^2))
\]
for some $r \geq 0$ and $s \geq 0$.

We now show that $r=1$. Suppose $r \geq 2$. 
It suffices to conclude that $X$ is not normal. 
Take the affine open subset defined by $z \neq 0$: 
\[
X' := D_+(z) = \Spec\,k[[t]][x, y]/(x^2 + t^r(y^2 + t^s) ). 
\]
By $r \geq 2$, the line defined by $x=t=0$ is contained in the singular locus of $X'$. 
Hence $X'$ is not normal. 
This completes the proof of $r=1$.

We then obtain 
\[
X \simeq  \Proj\,A[x, y, z]/(x^2 + ty^2+t^{s+1}z^2)
\]
with $s \geq 0$. Thus (1) holds.

The assertion (2) follows from Definition \ref{d-delta} and Definition \ref{d-disc}. 
The assertion (3) holds by either direct computation, or 
(2) and Theorem \ref{t-disc-vs-sing}. 
\qedhere



\end{proof}


\subsubsection{Reduced-fibre case with ${\rm char}\,k = 2$}\label{ss-surface-char2}


\begin{prop}\label{p-surface-char2-red}
Let $k$ be an algebraically closed field of characteristic two. 
Set $A := k[[t]]$ and let $0$ be the closed point of $\Spec\,A$. 
Let $f: X 
\to \Spec\,A$ be a conic bundle, 
where $X$ is canonical and $X_0$ is not smooth but {\cred is} reduced. 
Then 
the following hold. 
\begin{enumerate}
    \item  We have an $A$-isomorphism 
    \[
    X \simeq \Proj\,A[x, y, z]/(xy +t^{n+1}z^2)
    \]
    for some $n \geq 0$. 
        \item 
    $\deg \Delta_f = n+1$. 
    \item 
    If $n=0$, then $X$ is regular. 
    If $n\geq 1$, then $X$ has {\cred the} unique singularity and it is of type $A_{n}$. 
\end{enumerate}
\end{prop}

\begin{proof}
Let us show (1). 
By Proposition \ref{p-char2-local}(2), we can write 
\[
X \simeq \Proj\,A[x, y, z]/(xy +cz^2)
\]
for some $c \in \m := tk[[t]]$. 
Then we can write $ c = c't^m$ for some $c' \in A^{\times}$ and $m \in \Z_{\geq 0}$. 
By the following equality of ideals of $A[x, y, z]$
\[
(xy +cz^2) = (xy +c't^m z^2) =  ( (c'^{-1}x)y + t^mz^2), 
\]
we may assume that $c=t^m$. Since $X_0$ is not smooth, we get $m >0$. 
Thus (1) holds. 

The assertion (2) follows from Definition \ref{d-delta} and Definition \ref{d-disc}. 
The assertion (3) holds by either direct computation, or 
(2) and Theorem \ref{t-disc-vs-sing}. 
\end{proof}

\subsubsection{Non-reduced-fibre case with ${\rm char}\,k = 2$}


\begin{lem}\label{l-2dim-char2-nonred1}
Let $k$ be an algebraically closed field of characteristic two. 
Set $A := k[[t]]$ and let $0$ be the closed point of $\Spec\,A$. 
Let $f: X \to \Spec\,A$ be a conic bundle, 
where $X$ is canonical and $X_0$ is not reduced. 
Then 
an $A$-isomorphism 
\[
X \simeq \Proj\,A[x, y, z]/(ax^2 + by^2+cz^2 + t^n xy)
\]
holds for some $a, b, c \in A$ and $n \geq 1$. 
\end{lem}

\begin{proof}
We have 
\[
X = \Proj\,A[x, y, z]/(ax^2 + by^2 + cz^2 + \alpha yz + \beta zx + \gamma xy)
\]
for some $a, b, c \in A = k[[t]]$ and $\alpha, \beta, \gamma \in \m = tk[[t]]$. 
We may assume that $v_t(\alpha) \geq v_t(\beta) \geq v_t(\gamma) =:n$, 
where $v_t : k[[t]] \to \Z \cup \{\infty\}$ denotes the discrete valuation with $v_t(t)=1$. 
Then we can write 
\[
\alpha = t^n \alpha', \quad \beta = t^n \beta', \quad \gamma = t^n\gamma'
\]
for some $\alpha', \beta' \in A$ and $\gamma' \in A^{\times}$.  
Taking the multiple with $\gamma'^{-1}$, we get 
\[
X \simeq \Proj\,A[x, y, z]/(a'x^2 + b'y^2+c'z^2 + t^n( \alpha' yz + \beta' zx + xy))
\]
for some $a', b', c' \in A = k[[t]]$. 
It holds that 
\begin{eqnarray*}
\alpha' yz + \beta' zx + xy =
(x + \alpha' z)(y+\beta'z) - \alpha'\beta'z^2. 
\end{eqnarray*}
Applying the coordinate change $(x + \alpha' z, y+\beta'z, z) \mapsto (x, y, z)$, we obtain 
\[
X \simeq \Proj\,A[x, y, z]/(a''x^2 + b''y^2+c''z^2 +t^n xy)
\]
for some $a'', b'', c'' \in A$ and $n \in \Z_{\geq 0}$. 
Since $X_0$ is not reduced, we obtain $n>0$. 
\end{proof}

\begin{nota}\label{n-2dim-char2-nonred}
Let $k$ be an algebraically closed field of characteristic two. 
Set $A := k[[t]]$, $\m := tk[[t]]$, and $\kappa := A/\m (\simeq k)$. 
Let $0$ be the closed point of $\Spec\,A$. 
Let 
\[
f: X = \Proj\,A[x, y, z]/(Q) \to \Spec\,A
\] 
be a conic bundle, 
where $X$ is canonical and $X_0$ is not reduced. 
Assume that 
\[
Q = ax^2 + by^2 + cz^2 + t^n xy
\]
for some $a, b, c \in A$ and $n \geq 1$. 
We often write $a(t) :=a, b(t) :=b, c(t) :=c$. 
For $a(t) = a_0 + a_1 t + a_2 t^2+ \cdots$, 
it holds that $a(0) = a_0, a'(t) = a_1 + 2a_2t + 3a_3t^2 + \cdots$, and $a'(0) = a_1$. 
\end{nota}

\begin{lem}\label{l-2dim-char2-nonred2}
We use Notation \ref{n-2dim-char2-nonred}. 
Assume $n\geq 2$. Then the following hold. 
\begin{enumerate}
    \item 
    The singular locus of $X$ is given by the intersection of the following two lines on $\P^2_{\kappa}$: 
    \[
    \sqrt{a(0)}x + \sqrt{b(0)}y+\sqrt{c(0)}z = \sqrt{a'(0)} x + \sqrt{b'(0)}y +\sqrt{c'(0)}z=0.
    \]
    In particular, these lines are distinct and $X$ has {\cred the} unique singular point. 
    \item 
    We may assume that the singular point is 
    either $[0:1:0] \in \P^2_{\kappa}$ or $ [0:0:1]\in \P^2_{\kappa}$. 
    More precisely, there exist $\widetilde{a}, \widetilde{b}, \widetilde{c} \in A$ such that an $A$-isomorphism 
    \[
    X \simeq \Proj\,A[x, y, z]/(\widetilde{a}x^2 + \widetilde{b}y^2+\widetilde{c}z^2 + t^n xy)
    \] 
    holds and the singular point of 
    $\Proj\,A[x, y, z]/(\widetilde{a}x^2 + \widetilde{b}y^2+\widetilde{c}z^2 + t^n xy)$, which is unique by (1), 
    is either $[0:1:0] \in \P^2_{\kappa}$ or $ [0:0:1]\in \P^2_{\kappa}$. 
\end{enumerate}
\end{lem}

\begin{proof}
Let us show (1). 
Let $[x_0:y_0:z_0] \in \P^2_{\kappa}$ 
{\cred be} 
a singular point of $X$. 
By Proposition \ref{p-Jac-criterion}, 
this singular point  $(t, [x:y:z]) = (0, [x_0:y_0:z_0])$ is a solution of 
\[
Q = \partial_tQ = \partial_x Q =\partial_y Q = \partial_z Q =0. 
\]
By $t=0$, the following equalities automatically hold: 
\[
\partial_x Q =\partial_y Q = \partial_z Q =0. 
\]
By $n \geq 2$, the remaining conditions $Q = \partial_t Q=0$ can be written as  
\[
a(0)x^2 + b(0)y^2+c(0)z^2 = a'(0) x^2 + b'(0)y^2 +c'(0)z^2=0.  
\]
Since both are double lines, the singular locus is either a line or a point. 
If the singular locus is a line, then $X$ would not be normal. 
Thus (1) holds. 

Let us show (2). Note that the unique singular point can be written 
as $[\alpha : \beta: \gamma] \in \P^2_{\kappa}$. 
If $[\alpha : \beta :\gamma ]=[0:0:1]$, then there is nothing to show. 
By symmetry, we may assume that $\beta \neq 0$, and hence $\beta =1$. 
Applying $[x : y : z] \mapsto [x + \alpha y : y : z + \gamma y]$, 
we may assume that $\alpha = \gamma =0$. 
Thus (2) holds. 
\end{proof}

\begin{prop}\label{p-2dim-char2-nonred1}
We use Notation \ref{n-2dim-char2-nonred}. 
Assume that $n \geq 2$ and $[0:1:0]$ is the unique singularity of $X$. 
Then 
one and only one of (I) and (II) holds. 
\begin{enumerate}
    \item[(I)] 
    \begin{enumerate}
    \renewcommand{\labelenumii}{(\arabic{enumii})}
    \item  We have an $A$-isomorphism 
    \[
    X \simeq \Proj\,A[x, y, z]/(tx^2 +z^2 + t^n xy)
    \]
    for some $n > 0$. 
        \item 
    $\deg \Delta_f = 2n$. 
    \item 
\begin{itemize}
\item If $n=1$, then $X$ has exactly two singular points and both of them are of type $A_1$. 
\item If $n \geq 2$, then     $X$ has {\cred the} unique singular point and it is of type $D^0_{2n}$ in the sense of \cite[Page 16]{Art77}. 
    \end{itemize}
    \end{enumerate}
    \item[(II)] 
    \begin{enumerate}
    \renewcommand{\labelenumii}{(\arabic{enumii})}
    \item  We have an $A$-isomorphism 
    \[
    X \simeq \Proj\,A[x, y, z]/(x^2 +tz^2 + t^n xy)
    \]
    for some $n >0$. 
        \item 
    $\deg \Delta_f = 2n+1$. 
    \item 
    $X$ has {\cred the} unique singular point and it is of type $D^0_{2n+1}$ in the sense of \cite[Page 16]{Art77}. 
\end{enumerate}
\end{enumerate}

\end{prop}

\begin{proof}
Let us show (1). For $\gamma :=t^n$, 
we have 
\[
X = \Proj\,A[x, y, z]/(ax^2 + by^2 + cz^2 + \gamma xy). 
\]
Note that the condition $\gamma =t^n$ is not stable under the following argument.  
Since $[0:1:0]$ is a singular point of $X$, 
we obtain $b(0) = b'(0)=0$ (Lemma \ref{l-2dim-char2-nonred2}(1)). 
Recall that the singular point $[0:1:0]$ coincides with the solution of the following equation 
(Lemma \ref{l-2dim-char2-nonred2}): 
\begin{equation}\label{e1-2dim-char2-nonred1}
    \sqrt{a(0)}x +\sqrt{c(0)}z = \sqrt{a'(0)} x +\sqrt{c'(0)}z=0.
\end{equation}
In particular, 
\[
{\rm (I)}\,\,c(0) \neq 0 \qquad\qquad \text{or}\qquad\qquad {\rm (II)}\,\, a(0) \neq 0. 
\]
In what follows, we only treat 
(I), since the proofs are very similar. 


Assume (I), i.e., $c(0) \neq 0$. 
We get $c \in A^{\times}$. 
Taking the multiplication with $c^{-1}$, the defining equation of $X$ becomes 
\[
ax^2 + by^2 + z^2 + \gamma xy =0. 
\]
By the fact that $[0:1:0]$ is {\cred the} unique solution of (\ref{e1-2dim-char2-nonred1}), 
we obtain $a'(0) \neq 0$. 
Applying $z \mapsto z +\sqrt{a(0)} x$, we may assume that $a(0)=0$. 
Hence we can write 
\[
a =a(t)= a_1t +a_2t^2+ \cdots \in k[[t]], 
\]
where $a_1, a_2, ... \in k$ and $a_1 \neq 0$.  

We now erase the term $by^2$. 
Let $b_mt^m$ be the leading term of $b$, 
i.e.,  $b = b_mt^m + b_{m+1}t^{m+1} + \cdots \in k[[t]]$ with $b_m, b_{m+1}, ... \in k$ and $b_m \neq 0$. 
If $m = 2\ell$ for some $\ell \in \Z_{\geq 0}$ (i.e., $m$ is even), 
then we can erase $b_{2\ell}t^{2\ell}$ by applying $z \mapsto z+ \lambda y$ for suitable $\lambda \in k$. 
Hence we may assume that $m = 2\ell +1$ for some $\ell \in \Z_{\geq 0}$ (i.e., $m$ is odd). 
We can write 
\[
\gamma = \gamma_st^s + \gamma_{s+1}t^{s+1} + \cdots
\]
for some $\gamma_s, \gamma_{s+1}, ... \in k$ with $\gamma_s \neq 0$. 
Since $X_0$ is non-reduced, we have $s \geq 1$. 
We treat the following two cases separately. 
\begin{enumerate}
\renewcommand{\labelenumi}{(\roman{enumi})}
\item $2\ell + 1 \leq s +\ell$. 
\item $2\ell + 1 > s + \ell$. 
\end{enumerate}
Assume (i). 
In this case, we apply the $A$-linear transform 
$(x, y, z) \mapsto (x+\mu t^{\ell}y, y, z)$ for some $\mu \in k$,
so that the defining polynomial of $X$ becomes 
\[
a(x+\mu t^{\ell} y)^2 + by^2 + z^2 + \gamma(x+\mu t^{\ell} y)y 
= ax^2 + (a\mu^2t^{2\ell} + b + \gamma \mu t^{\ell} )y^2 + z^2 + \gamma xy.  
\]
If $\mu \neq 0$, then the leading term of $\gamma \mu t^{\ell}$ is equal to 
$\gamma_s \mu t^{\ell+s}$, which is of degree $\geq 2 \ell +1$ 
by the assumption $2\ell + 1 \leq s +\ell$. 
Therefore, the leading term of $a\mu^2t^{2\ell} + b + \gamma \mu t^{\ell}$  is of degree $\geq 2\ell +1$. 
It is enough to find $\mu \in k$ that makes this inequality strict. 
The coefficient of $t^{2\ell+1}$ in $a\mu^2t^{2\ell} + b + \gamma \mu t^{\ell} \in k[[t]]$ is equal to 
\[
a_1 \mu^2 + b_{2\ell +1} + \gamma_{\ell +1} \mu, 
\]
where we set $\gamma_{\ell+1} :=0$ when $\ell +1 < s$. 
By $a_1 \neq 0$, there is a solution $\mu \in k$ of the equation $a_1 \mu^2 + b_{2\ell +1} + \gamma_{\ell +1} \mu =0$. 
This completes the case when (i). 
Hence we may assume (ii). 
In this case, we apply the  $A$-linear transform 
$(x, y, z) \mapsto (x+\nu t^{2\ell+1 -s}y, y, z)$ for some $\nu \in k$, 
so that the defining polynomial of $X$ becomes 
\[
a(x+\nu t^{2\ell+1 -s} y)^2 + by^2 + z^2 + \gamma(x+\nu t^{2\ell+1 -s} y)y 
\]
\[
= ax^2 + (a\nu^2t^{4\ell+2 -2s} + b + \gamma\nu t^{2\ell+1 -s}  )y^2 + z^2 + \gamma xy.  
\]
The leading term of $a\nu^2t^{4\ell+2 -2s}$ is of degree $1 +4\ell+2 -2s$, 
which is larger than $2\ell +1$ by the assumption $2\ell + 1 > s + \ell$. 
Hence the coefficient of $t^{2\ell+1}$ in 
$a\nu^2t^{4\ell+2 -2s} + b + \gamma\nu t^{2\ell+1 -s} \in k[[t]]$ is equal to 
$b_{2\ell+1} + \gamma_s \nu$. 
Hence it is enough to set $\nu := \gamma_s^{-1} b_{2\ell+1}$. 
Therefore, we may assume that $b=0$. 

Then the defining equation of $X$ becomes 
\[
ax^2 + z^2 + \gamma xy =0. 
\]
By using $z^2$, we may assume that $a$ only has odd terms, i.e., $a = a_1t + a_3t^3 + a_5t^5+ \cdots \in k[[t]]$. 
We can write 
\[
ax^2 = (a_1t+ a_3t^3 + a_5t^5+ \cdots)x^2 = t ((\sqrt{a_1}+ \sqrt{a_3}t + \sqrt{a_5}t^2+ \cdots)x)^2, 
\]
and hence we may assume that $a=t$. 
Finally, we can assume that $\gamma = t^n$ after replacining $y$ with $uy$ for some unit $u \in A^{\times}$. 
Since $X_0$ is not reduced, we get $n >0$. 
Thus (1) holds. 

The assertion (2) follows from Definition \ref{d-delta} and Definition \ref{d-disc}. 
The assertion (3) holds by Theorem \ref{t-disc-vs-sing} and \cite[page 16]{Art77}. 
\qedhere



\end{proof}

\begin{ex}
Let $k$ be an algebraically closed field of characteristic two. 
Set $A :=k[[t]]$ and let $0 \in \Spec\,A$ be the closed point. 
Fix $r \in \Z_{>0}$ and $s \in \Z_{>0}$. 
\begin{enumerate}
\item[(I)]  
Consider a conic bundle: 
\[
f:  X := \Proj\,A[x, y, z]/(x^2 +ty^2 +t^{2r}z^2 + t^s xy) \to \Spec\,A.  
\]
Then the following hold. 
\begin{enumerate}
\renewcommand{\labelenumii}{(\arabic{enumii})}
\item $\deg \Delta_f = 2r+2s$. 
\item $[0:0:1] \in X_0$ is {\cred the} unique singularity of $X$ and it is of type 
    $D_{2r+2s}^r$ in the sense of \cite[Page 16]{Art77}. 
\end{enumerate} 
\item[(II)]  
Consider a conic bundle: 
\[
f:  X := \Proj\,A[x, y, z]/(x^2 +ty^2 +t^{2r+1}z^2 + t^s xy) \to \Spec\,A.  
\]
Then the following hold. 
\begin{enumerate}
\renewcommand{\labelenumii}{(\arabic{enumii})}
\item $\deg \Delta_f = 2r+2s+1$. 
\item $[0:0:1] \in X_0$ is {\cred the} unique singularity of $X$ and it is of type 
    $D_{2r+2s+1}^r$ in the sense of \cite[Page 16]{Art77}. 
\end{enumerate} 
\end{enumerate}
\end{ex}

\begin{proof}
We only prove (I). 
The assertion (1)  follows from Definition \ref{d-delta} and Definition \ref{d-disc}. 
Let us show  (2). 
We now prove that $[0:0:1] \in X_0$ is {\cred the} unique singularity of $X$. 
Take a singular point $[x_0 : y_0:z_0] \times t_0 \in \P^2_{x, y, z} \times \A^1_t$ 
of $\widetilde X := \{x^2 +ty^2 +t^{2r}z^2 + t^s xy=0\} \subset \P^2_{x, y, z} \times \A^1_t$. 
By the Jacobian criterion (Proposition \ref{p-Jac-criterion}), we obtain the following equation: 
\[
x_0^2 +t_0y_0^2 +t_0^{2r}z_0^2 + t_0^s x_0y_0=0, 
\]
\[
t^s_0 y_0=0, \qquad t_0^sx_0=0,\qquad y_0^2 + st_0^{s-1}x_0y_0=0.
\]
If $t_0\neq 0$, then we would get $x_0=y_0=0$, which leads to a contradiction $z_0=0$. 
Hence $t_0=0$. Then we get $x_0=y_0=0$. 
Thus $[0:0:1] \times \{0\} \in \P^2_{x, y, z} \times \A^1_t$ is {\cred the} unique singularity 
of $\widetilde X := \{x^2 +ty^2 +t^{2r}z^2 + t^s xy=0\}$. 
Thus $[0:0:1] \in X_0$ is {\cred the} unique singularity of $X$.


The singularity is the origin of 
\[
x^2 + ty^2 + t^{2r} + t^s xy \in k[[t]][x, y]. 
\]
Applying the $k[[t]]$-linear transform $(x , y) \mapsto (x+t^r, y)$, 
the equation becomes 
\[
x^2 + ty^2 + t^s xy + t^{r+s}y \in k[[t]][x, y]. 
\]
This is of type $D_{2r+2s}^r$ in the sense of \cite[Page 16]{Art77}. 
Thus (2) holds. 
\end{proof}

\section{Special phenomena in characteristic two}\label{s-char2-ex}



Although $\Delta_f$ can be non-reduced under our definition (e.g., Theorem \ref{t-disc-vs-sing}), 
$\Delta_f$ was classically defined as its reduced structure $(\Delta_f)_{\red}$ \cite{Bea77}, \cite{MM83}. 
In this section, we shall observe that 
both $\Delta_f$ and $(\Delta_f)_{\red}$ can behave worse in characteristic two than  characteristic $\neq 2$ 
(Subsection \ref{ss-ex-char2}). 
We start by establishing some results which hold in characteristic $\neq 2$ (Subsection \ref{ss-results-char0}).

\subsection{Results in characteristic $\neq 2$}\label{ss-results-char0}

\begin{prop}\label{p-sm-nonred-fib1}
Let $f: X \to T$ be a conic bundle of noetherian regular $\Z[\frac{1}{2}]$-schemes. 
Assume that $\dim T \leq 1$. 
Then $\Delta_f$ is regular and any fibre of $f$ is geometrically reduced. 
\end{prop}

\begin{proof}
If $\dim T=0$, then there is nothing to show. 
In what follows, we assume that $\dim T=1$ and $T$ is an integral scheme. 
Note that $f$ is generically smooth (Lemma \ref{l-char0-gene-sm}). 
By Proposition \ref{p-sm-red-fib}, it is enough to show that  any fibre of $f$ is geometrically reduced. 
Suppose that there exists a point $0 \in T$ such that $X_0$ is not geometrically reduced. Let us derive a contradiction. 
Taking the strict henselisation of $\MO_{T, 0}$, we may assume that $T=\Spec\,A$, 
where $(A, \m, \kappa)$ is a strictly henselian local ring with $\dim A=1$. 
In particular, $A$ is a discrete valuation ring. 
Let $t$ be a uniformiser, i.e., $\m =tA$. 
By Corollary \ref{c-char0-local2}, we can write 
\[
X = \Proj\,A[x, y, z]/(Q), \qquad  Q = x^2 + by^2 + cz^2 + \alpha yz,
\]
for some $b, c, \alpha \in \m$. 
Applying a suitable linear transform, we may assume that $b \neq 0$ 
(cf. the proof of Proposition \ref{p-surface-char0-nonred}(1)). 
We can write $b = t^{\ell} \widetilde{b}$ for some $\ell \in \Z_{>0}$ and $\widetilde b \in A^{\times}$. 
By Hensel's lemma, we may assume $\widetilde b =1$: 
\[
Q = x^2 + t^{\ell} y^2 + c z^2 + \alpha yz. 
\]
Since the open subset $D_+(y) = \{ x^2 + t^{\ell} + c z^2 + \alpha z =0\}$ of $X$ is regular, 
we obtain $\ell =1$. 
We have $\alpha = t \widetilde{\alpha}$ for some $\widetilde \alpha \in A$. 
By completing a square: 
\[
t^{\ell} y^2  + \alpha yz = t y^2 + t \widetilde \alpha yz = 
t \left(y+\frac{\widetilde \alpha}{2}z\right)^2 - \frac{t\widetilde \alpha^2}{4}z^2, 
\] 
we may assume that $\alpha =0$, i.e., 
\[
Q = x^2 + t y^2 + c z^2. 
\]
By the same argument as above, we may assume that $c=t$: 
\[
Q = x^2 + t y^2 + t z^2. 
\]
Then $X$ is not regular at the closed point $[0:1: \sqrt{-1}] \in \P^2_{\kappa}$ over $0 \in T$, which is absurd. 
\end{proof}

\begin{prop}\label{p-sm-nonred-fib}
Let $f: X \to T$ be a conic bundle of noetherian regular $\Z[\frac{1}{2}]$-schemes. 
Then the following hold. 
\begin{enumerate}
    \item $\Delta_f$ is reduced. 
    \item 
    If $T$ is a smooth surface over an algebraically closed field $k$ of 
    {\cred characteristic} $\neq 2$ (in particular, $X$ is a smooth threefold over $k$), then $\Delta_f$ is normal crossing.  
\end{enumerate}
\end{prop}
\begin{proof}
The assertion (1) holds by applying Proposition \ref{p-sm-nonred-fib1} 
after taking the base change $\Spec\,\MO_{T, \xi} \to T$ 
for a generic point {\cred $\xi$} of $\Delta_f$.  
The assertion (2) holds by (1) and the proof of \cite[Ch. I, Proposition 1.2(iii)]{Bea77}. 
\end{proof}

\subsection{Examples in characteristic two}\label{ss-ex-char2}

The following example shows that singularities of $\Delta_f$ and $(\Delta_f)_{\red}$ can be arbitrarily bad even if $X$ is a smooth threefold and $T$ is a smooth surface (cf. Proposition \ref{p-sm-nonred-fib}(2)).

\begin{ex}\label{e-3fold-char2}
We work over an algebraically closed field $k$ of characteristic two. 
Set 
\[
T := \mathbb A^2_{u, v}  = \Spec\,k[u, v]\qquad\text{and}\qquad \mathbb P^2_{x, y, z} := {\rm Proj}\,k[x, y, z].
\]
Fix $\varphi(u, v) \in k[u, v] \setminus \{0\}$. 
Set 
\[
X := \{ x^2 + uy^2 + vz^2 + \varphi(u, v)^2 yz\} \subset 
\mathbb A^2_{u, v} \times \mathbb P^2_{x, y, z}= T \times \mathbb P^2_{x, y, z}. 
\]
Then the induced morphism $f :X \to T$ is a conic bundle. 
The discriminant scheme $\Delta_f \subset T$ of $f$ is given as follows 
(Definition \ref{d-delta} and Definition \ref{d-disc}): 
\[
\Delta_f = \{ \varphi(u, v)^4 =0 \} \subset T = \mathbb A^2_{u, v}.  
\]
In particular, 
the conic bundle $f:X \to T$ is generically smooth by $\varphi(u, v) \neq 0$.

Let us show, by using a Jacobian criterion (Proposition \ref{p-Jac-criterion}),  that $X$ is smooth over $k$. 
Suppose that a closed point $(u_0, v_0) \times [x_0:y_0:z_0] \in \mathbb A^2_{u, v} \times \mathbb P^2_{x, y, z}$ is a solution of 
\[
Q= \partial_uQ = \partial_vQ =\partial_xQ = \partial_yQ = \partial_zQ=0,  
\]
\[
\text{where} \qquad 
Q := x^2 + uy^2 + vz^2 + \varphi(u, v)^2 yz.
\]
By $\partial_u Q=0$ and $\partial_v Q=0$, we have $y_0 = z_0=0$. 
Then $[x_0: y_0:z_0]=[1:0:0]$.  However, this is not a solution of $Q=0$, 
which is absurd. 
Hence $X$ is smooth. 
This conic bundle $f:X \to T$ satisfies the following properties. 
\begin{enumerate}
\item For any point $t \in \Delta_f$, its fibre $X_t$ is not geometrically reduced. 
\item The discriminant divisor $\Delta_f \subset T$ is set-theoretically given by $\{\varphi(u, v)=0\}$
In particular, the following phenomena happen in characteristic two, 
each of which does not occur in any other characteristic. 
\begin{itemize}
\item Assume that $\varphi(u, v)=u$. 
In this case, $(\Delta_f)_{\red}$ is a smooth curve, and $X_t$ is non-reduced for any $t \in (\Delta_f)_{\red}$ (cf. Theorem \ref{t-sm-red-fib}, Proposition \ref{p-sm-nonred-fib}). 
\item Assume that $\varphi(u, v)=uv(u+v)$. 
In this case, $(\Delta_f)_{\red}$ is not normal crossing 
 (cf. Proposition \ref{p-sm-nonred-fib}). 
\end{itemize} 
\end{enumerate}
As $\varphi(u, v)$ is chosen to be an arbitrary nonzero element of $k[u, v]$, 
the singularities of $\Delta_f$ and $(\Delta_f)_{\red}$ can be arbitrarily bad. 
\end{ex}


\begin{ex}\label{e-surface-char2}
We use the same notation as in Example \ref{e-3fold-char2}. 
Set $\varphi(u, v):=u$. We have  
\[
X = \{ x^2 + uy^2 + vz^2 + u^2 yz\} \subset T \times \mathbb P^2_{x, y, z}. 
\]
Take the localisation $T'$ of $T = \Spec\,k[u, v]$ at the generic point of $\{ u =0\}$, i.e., 
\[
T' := \Spec\,k[u, v]_{(u)} = \Spec\,k(v)[u]_{(u)}. 
\]
Then $k(v)[u]_{(u)}$ is a 
discrete valuation ring. 
Take the base change $X' := X \times_T T'$ and let $f' : X' \to T'$ be the induced 
conic bundle. 
Then the following hold.
\begin{enumerate}
    \item $X'$ and $T'$ are regular. 
    \item $\dim T' = 1$. 
    \item $f' : X' \to T'$ is generically smooth. 
    \item $\Delta_{f'} = \Spec\,(k(v)[u]_{(u)}/(u^4))$ is not reduced. 
\end{enumerate}
In particular, Proposition \ref{p-sm-nonred-fib1} does not hold in characteristic two even if we impose the generically smooth assumption. 
\end{ex}

\begin{ex}[Fano threefolds with non-reduced discriminants]\label{e-Fano3-char2}
We work over an algebraically closed field $k$ of characteristic two. 
Set  
\[
X := \{ s x^2 + ty^2 + uz^2 + s yz\} \subset \mathbb P^2_{x,y,z} \times \mathbb P^2_{s, t,u}. 
\]
and let $f : X \to T :=\P^2_{s, t, u}$ be the induced morphism, which is a conic bundle. 
By Proposition \ref{p-Jac-criterion}, 
we can check that $X$ is smooth over $k$. 
It follows from $X \in |\MO_{\P^2 \times \P^2}(2, 1)|$ that $X$ is a smooth Fano threefold. 
We have 
\[
\Delta_f = \Proj\,k[s, t, u]/(s^3),  
\]
which is  non-reduced. 
Furthermore, $(\Delta_f)_{\red}$ is a smooth rational curve. 
\end{ex}

\begin{rem}[Failure of Bertini]\label{r-Bertini-fail}
We use the same notation as in Example \ref{e-Fano3-char2}. 
For any smooth curve $T'$ on $T = \P^2_{s, t, u}$, 
its inverse image $X' := f^{-1}(T')$ is not smooth.  
Indeed, if $X'$ is  smooth, 
then the resulting conic bundle $f': X' \to T'$ is a generically smooth conic bundle from a smooth 
surface $X'$ to a smooth curve $T'$, which {\cred implies} 
that $\Delta_{f'}$ is reduced 
(Proposition \ref{p-sm-red-fib}, Lemma \ref{l-sm-surface}). 
However, this is absurd, because $\Delta_f$ is non-reduced and hence so is $\Delta_{f'}$ (Remark \ref{r-disc-bc}). 

Note that the original proof of the Mori-Mukai formula (Theorem \ref{intro-MM}) in characteristic zero depends on the Bertini theorem for base point free divisors \cite[Proposition 6.2]{MM83}. 
This example shows that the same argument does not work in characteristic two. 
\end{rem}

\begin{rem}[Wild conic bundles]
We work over an algebraically closed field $k$ of characteristic two. 
If $f: X \to T$ is a wild conic bundle, then it is hard from its definition (Definition \ref{d-disc-bdl}) to know what is $\Delta_f^{\bdl}$. 
When $f: X \to T$ is a wild conic bundle 
from a smooth projective threefold $X$ to a smooth projective surface $T$, 
the following holds: 
\[
\Delta^{\bdl}_f \overset{{\rm (i)}}{\equiv} -f_*(K_{X/T}^2) \overset{{\rm (ii)}}{\equiv} -K_T, 
\]
where (i) follows from Theorem \ref{t-cbf} and (ii) holds by \cite[Corollary 4]{MS03}. 
\end{rem}


\bibliographystyle{skalpha}
\bibliography{library.bib}

\end{document}